\RequirePackage{fix-cm}
\documentclass{svjour3}                     % onecolumn (standard format)
\smartqed  % flush right qed marks, e.g. at end of proof
\usepackage{mathptmx}      % use Times fonts if available on your TeX system
\usepackage{amsmath,amssymb,amsfonts}%
\usepackage{graphicx}
\usepackage{subfigure}
\usepackage{stmaryrd}
\usepackage[T1]{fontenc}
\usepackage{latexsym,mathrsfs}
\usepackage{float}
\usepackage{xcolor}
\usepackage{algorithm, algorithmic}
\usepackage{array,bm}
\usepackage[colorlinks,linkcolor=blue, anchorcolor=blue, citecolor=blue]{hyperref}
\usepackage{color,caption,cases}
 \usepackage{tabularx,booktabs }
\usepackage{indentfirst}
\usepackage{appendix}
\usepackage[numbers,sort&compress]{natbib}
\numberwithin{equation}{section}
\newcommand{\malA}{\mathcal{A}}
\newcommand{\malB}{\mathcal{B}}

\newcommand{\malH}{\mathcal{H}}
\newcommand{\malO}{\mathcal{O}}
\newcommand{\malS}{\mathcal{S}}

\journalname{Arxiv }

\begin{document}
\title{Fractional-step High-order and Bound-preserving Method for Convection Diffusion Equations 
}	
\titlerunning{ Fractional-step High-order and Bound-preserving Method for Convection Diffusion Equations }        % if too long for running head
			
\author{
	Baolin Kuang    \and
	Hongfei Fu \and 
	Shusen Xie$^*$ 
}
\authorrunning{ Baolin Kuang \and Hongfei Fu \and Shusen Xie } % if too long for running head

\institute{ 
	Shusen Xie \at
	School of Mathematical Sciences \& Laboratory of Marine Mathematics, Ocean University of China, Qingdao  266100, China\\
	\email{shusenxie@ouc.edu.cn}
	\and
	Baolin Kuang  \at
	School of Mathematical Sciences, Ocean University of China, Qingdao  266100, China\\
	\email{blkuang@stu.ouc.edu.cn} 
	\and 
  Hongfei Fu		\at
	Corresponding author. School of Mathematical Sciences \& Laboratory of Marine Mathematics, Ocean University of China, Qingdao  266100, China\\
	\email{fhf@ouc.edu.cn} 
}

\date{Received: date / Accepted: date}
% The correct dates will be entered by the editor

\maketitle
 
\begin{abstract}
In this paper, we derive two bound-preserving and mass-conserving schemes based on the fractional-step method and high-order compact (HOC) finite difference method for nonlinear convection-dominated diffusion equations.  We split the one-dimensional equation into three stages, and employ appropriate temporal and spatial discrete schemes respectively.
We show that our scheme is weakly monotonic and that the bound-preserving property can be achieved using the limiter of Li et al. (SIAM J Numer Anal 56: 3308–3345, 2018) under some mild step constraints.
By employing the alternating direction implicit (ADI) method, we extend the scheme to two-dimensional problems, further reducing computational cost. We also provide various numerical experiments to verify our theoretical results.
%The transport part is solved by explicit strong stability preserving Runge-Kutta substeppings combined with HOC finite difference scheme. And the diffusion part is discretized by the implicit Crank-Nicolson and HOC scheme.
% \PACS{PACS code1 \and PACS code2 \and more}

\end{abstract}

% \begin{keywords}
%  {HOC finite difference method \and bound-preserving \and mass conservation \and fractional-step method \and ADI methods}
% \end{keywords}

%%%%%%%%%%%%%%%%%%%%%%%%%%%%%%%%%%%%%%%%%
\section{Introduction}
%%%%%%%%%%%%%%%%%%%%%%%%%%%%%%%%%%%%%%%%%

 \qquad Consider the nonlinear convection diffusion equation with periodic boundary condition
\begin{equation}
\begin{aligned}
& u_t + \nabla \cdot \mathbf{F}(u)- \nu \Delta u = S (\mathbf{x},t),
\quad \mathbf{x} \in \Omega, t \in(0, T], \label{model:all}\\
&    u(\mathbf{x},0)= u_{0}(\mathbf{x}),
\end{aligned}
\end{equation}
where $\Omega$ is a bounded domain in $\mathbb{R}$ or $\mathbb{R}^2$. The symbols $\Delta$ and $\nabla\cdot$ denote the Laplacian and divergence operators, respectively. The diffusion coefficient $\nu$ is typically a small positive constant in convection-dominated problems.
Here, $\mathbf{F}(u)$ represents the flux vector, and $S(\mathbf{x}, t)$ denotes the source term, both of which are well-defined smooth functions.

convection diffusion equations are fundamental mathematical physics models that are widely used to describe various physical, chemical, and biological phenomena, such as the transport of groundwater in soil \cite{parlange1980water}, the long-distance transport of pollutants in the atmosphere \cite{zlatev1984implementation}, and the dispersion of tracers in anisotropic porous media \cite{fattah1985dispersion}, etc. Of particular importance is the case where convection dominates,
i.e., {when the diffusion coefficient is small compared to other scales}, which is often used in models of two phase flow in oil reservoirs\cite{peaceman2000fundamentals}.
However,  convection-dominant diffusion problems exhibit strong hyperbolic properties, making them challenging for standard finite difference or finite element methods. These methods often face issues such as computational complexities and non-physical numerical oscillations, especially when standard methods struggle to resolve steep gradients and ensure mass conservation, both of which are crucial for realistic applications.
On the other hand, in some physical processes, the physical quantity $u$ represents the density, concentration, or pressure of a certain substance or population. Therefore, the solution must be non-negative and often distributed in the interval $[0,1]$. In such cases, negative values are physically meaningless, such as in radionuclide transport calculations \cite{genty2011maximum}, chemotaxis problems \cite{zhang2022positivity}, streamer discharges simulations \cite{zhuang2014positivity}.
In addition, mass conservation is an important physical property. Under the periodic boundary condition, the solution of problem \eqref{model:all} satisfies
\begin{equation} \label{model:mass}
 \int_\Omega u(\mathbf{x}, t) d \mathbf{x} = \int_\Omega u^0(\mathbf{x}) d \mathbf{x}+ \int_0^t \int_\Omega S(\mathbf{x},s) d \mathbf{x} ds,
\end{equation}
which indicates that \eqref{model:all} inherently conserves mass and it is essential for numerical schemes to also conserve mass in the discrete sense.
% Therefore, while ensuring high order accuracy, constructing an efficient numerical method that preserves physical properties is of great research significance and scientific value. 
%  Otherwise, negative values may result in
% ill-posedness and instability for certain nonlinear equations \cite{srinivasan2018positivity,van2019positivity}.
% Moreover, the computational cost for enforcing these properties should not be significant if it is needed for each time step.
Hence, it is worthwhile to study an effective and high-order bound-preserving numerical scheme for nonlinear convection diffusion equations.

% the fractional step high-order numerical approximation of
As we knew, high-order schemes can more effectively retain the high order information of the nonlinear problem itself. Moreover, under the same error accuracy requirements, high-order schemes can use a coarser subdivision grid, which significantly improves the calculation efficiency. Therefore, it is meaningful to construct high-order schemes for the aforementioned model problems.
Although general high-order methods provide high accuracy, they often suffer from poor stability and wide grid templates, which degrade the matrix's sparse structure corresponding to the numerical schemes. To overcome these shortcomings and improve computational efficiency, the high-order compact (HOC) difference method was proposed. This method achieves as high as possible with the smallest grid template.
 Since compact schemes have a spectral-like resolution and extremely low numerical dissipation, the HOC method has aroused extensive research interest over the past few decades, and various effective HOC difference schemes have been developed successively \cite{wang1d,Zhang2002HOC,qiao2008,li2018A,mohebbi2010high}.
%  For instance, Chu et al. \cite{chu1998three} developed a new three-point combined compact difference scheme with sixth-order accuracy at periodic boundaries and fifth-order accuracy at non-periodic boundaries. 
% used the Steklov average operator to construct a fourth-order HOC difference scheme for the three-dimensional convection diffusion equation with constant coefficients and provided optimal error analysis using the Bramble-Hilbert lemma. Mohebbi et al. \cite{mohebbi2010high} combined the time cubic spline configuration method to propose a fourth-order HOC difference scheme in both time and space for the one-dimensional heat equation and convection diffusion equation. Wang \cite{wang2015fourth} proposed an effective Picard monotone iterative fourth-order HOC algorithm for the two-dimensional quasi-linear elliptic boundary value problem and proved the sufficient condition for the unique existence of the numerical solution and the optimal error estimate using the upper and lower solution methods.
%
% Although HOC difference schemes provide highly accurate numerical solutions on relatively coarser grids, numerical experiments presented in \cite{2002A} showed high cost in using linear sparse solvers to invert the resulting matrix at each time step when solving time-dependent problems.
% To further reduce computational costs, Karaa and Zhang \cite{karaa2004high} proposed a high order alternating direction implicit (ADI) method for solving 2D convection diffusion problems with constant coefficients, which is fourth-order in space and second-order in
% time.
ADI methods have proven to be very efficient in the approximation of the solutions of multi-dimensional evolution problems by transferring the numerical solutions of a high-dimensional problem into numerical solutions of a successive of independent one-dimensional problems. 
To further reduce computational costs, Karaa and Zhang \cite{karaa2004high} proposed a compact alternating direction implicit (ADI) method, which only require
solving one-dimensional implicit problems for each time step. 
% Tian and Ge \cite{tian2007fourth} derived a compact ADI scheme by using a spatial discrete exponential fourth-order compact difference formula and the Crank-Nicolson (C-N) method for the time discretization.
% Tian \cite{tian2011rational} also proposed another unconditionally stable rational compact ADI difference method.
% This method is unconditionally stable and, compared to \cite{karaa2004high}, it has a smaller dissipation error and better resolution properties.

On the other hand, preserving physical properties is a crucial reference point for the construction of numerical schemes. In recent years, scholars have developed many structure-preserving numerical methods for the convection diffusion equation.
 Bertolazzi et al. \cite{berikelashvili2007convergence} proposed a finite volume scheme based on the second-order maximum-preserving principle.
 Rui et al. \cite{rui2010mass} developed a mass-conserving characteristic finite element scheme.
 Shu et al. \cite{zhang2012maximum} proposed a high order finite volume WENO scheme that satisfies the maximum principle.
 Xiong et al. \cite{xiong2015high} studied the discontinuous Galerkin method based on the maximum value preservation principle.
In fact, the implementation of high-order structure-preserving difference methods is challenging. Typically, weak monotonicity does not apply to general finite difference schemes. However, Li et al. \cite{li2018A} demonstrated that certain high-order compact finite difference schemes satisfy such properties. This implies that a simple limiter can ensure the upper and lower bounds of the numerical solution without reducing the order and while maintaining mass conservation.
Although many studies have extensively explored high-order compact finite difference methods, this is the first time that weak monotonicity in compact finite difference approximations has been discussed, and it is also the first instance where high-order finite difference schemes have established weak monotonicity. However, applying limiters, especially when handling high-dimensional problems, inevitably incurs significant computational costs. Furthermore, the step size limitation in \cite{li2018A} is relatively stringent, suggesting that there is room for further improvement in the algorithm.
In conclusion, while various methods have been proposed to tackle the challenges posed by convection-dominated diffusion equations, the construction of an effective and high-order bound-preserving numerical scheme remains a crucial and ongoing area of research.

Fractional step methods, also known as operator splitting methods, are widely used for many complex time-dependent models. A notable example is the well-known Strang method \cite{strang1968}, which is the splitting method applied in this paper. The basic idea is to divide the original complex problem into a series of smaller problems, allowing us to solve simpler systems instead of tackling the larger system \cite{JIA2011387}. 
% The splitting technique is generally used in one of two ways \cite{hundsdorfer2003numerical}:
% \begin{itemize}
%  \item Coordinate Axis Splitting: Each split system only involves information along one of the coordinate axes.
%  \item Physical Phenomenon Splitting: The difference operator is split into several parts, each representing a particular physical phenomenon.
% \end{itemize}
%
The splitting method takes into account the different properties of the components, allowing the use of appropriate numerical methods for each part. This approach transforms a system of nonlinear difference equations into a linear system and a set of ordinary difference equations, making implementation easier and reducing computational cost \cite{einkemmer2016overcoming}. Additionally, the method facilitates parallel implementation in many cases and enables the construction of numerical methods with better geometric properties. For these reasons, the operator splitting method is well-suited for the numerical approximation of models describing complex phenomena.
To reduce the constraints of using the bound-preserving limiter method and decrease the computational cost, we consider combining the Strang method to improve the overall approach. This involves strongly decoupling the convection and diffusion terms and handling them separately. For high-dimensional problems, incorporating the ADI method further reduces the multidimensional problem to the computational complexity of several one-dimensional problems.

The remainder of this paper is organized as follows.
In Section 2, we obtain an improved HOC difference scheme of one dimension, and show that under some reasonable step constraints, the limiter can enforce the bound-preserving property;
In Section 3, we generalize the scheme to two-dimension, and also give the restrictions on bound-preserving bound;
Numerical examples are provided in Section 4.

%%%%%%%%%%%%%%%%%%%%%%%%%%%%%%%%%%%%%%%%%%%%%%%%%%%
 \section{ One-dimensional problems } \label{sec:1d}
%%%%%%%%%%%%%%%%%%%%%%%%%%%%%%%%%%%%%%%%%

\qquad In this section, we aim to construct a Strang splitting method based on HOC difference schemes for the one-dimensional convection diffusion equation. Under some mild step constraints and by incorporating the bound-preserving limiter (BP limiter) \cite{li2018A}, we will demonstrate that the proposed scheme is weakly monotonic and has the bound-preserving property.

% by using bound-preserving limiter.
%on the. satisfying the bound-preserving property \eqref{1d:MP} through from effectively.

We consider the convection diffusion equation with periodic boundary condition:
 \begin{equation}
\begin{aligned}
& u_t + f(u)_x -\nu u_{xx} = S(x,t),
 \quad 0\leq x \leq L,\ 0\leq t \leq T \label{model:1d} \\ 
& u(x, 0) =u^{0} (x),
\end{aligned}
\end{equation}
 where $\nu $ is a small positive constant compared with $|f'(u)|$.

%==========================================================
\subsection{Strang splitting temporal discretization} \label{sec_1d_time}
%===========================================================

\qquad For a positive integer $N_t$, the total number of time steps, let $\tau = {T}/{N_t}$ be the time step size, and $t_n = n \tau(n=0,1,\ldots,N_t)$.
We denote the exact solution operator associated with the nonlinear hyperbolic equation by $\malS_{N}$:
\begin{align}
 {u}_t+ f(u)_x =S(x,t), \label{1d_hyper}
\end{align}
and the exact solution operator associated with the linear parabolic equation by $\malS_{L}$ :
\begin{align}
 {u}_t= \nu u_{xx}. \label{1d_para}
\end{align}
Then, the second-order Strang splitting approximation of \eqref{model:1d} involves three substeps for approximating the solution at time $t_{n+1}$ from input solution $u(x,t_n)$\cite{strang1968}:
\begin{align}
 {u}({x}, t_n+\tau)=\malS_{L}(\tau / 2) \malS_{N}(\tau) \malS_{L}(\tau / 2) {u}({x}, t_n) + \malO (\tau^3). \label{strang}
\end{align}
That is, on each time interval $\left(t_n, t_{n+1}\right)$, Strang splitting approximation is reduced to the following process:
\begin{numcases}{}
 u_{t} -\nu u_{xx} = 0, & $t\in [t_n, t_{n+1/2}]; \quad u^{*}({x,t_n})=u^n(x)$, \label {1d:Strang:diffusion1} \\
 u_{t} + f(u)_x = S(x, t), & $t\in [t_n, t_{n+1}]; \qquad u^{**}(x,t_n)=u^{*}(x,t_{n+{1}/{2}})$, \label {1d:Strang:hyperbolic} \\
 u_{t} -\nu u_{xx} = 0, &  $t\in [t_{n+1/2}, t_{n+1}]; \quad u^{***}(x,t_{n+{1}/{2}})=u^{**}(x,t_{n+1})$. \label {1d:Strang:diffusion2}
\end{numcases} 
% In addition, as shown in [\cite{hundsdorfer2003numerical}, Chapter IV], the above scheme \eqref{1d:Strang:diffusion1}--\eqref{1d:Strang:diffusion2} bears a splitting error of $O\left(\tau^2\right)$.
% \begin{lemma} \label{lem:strang}
% Suppose that $f(\mathbf{u})(\boldsymbol{x}, t) \in C^2\left((0, T) ;\left[C^2(\Omega)\right]^d\right)$ and the exact solution to convection diffusion equation $u \in C^3\left((0, T) ; C^2(\Omega)\right)$. Let $\left\{u^n\right\}_{n=0}^N$ be the numerical solution produced by the second-order splitting scheme
% \eqref{scheme:1d:semi_t:1}--\eqref{scheme:1d:semi_t:4}. There exists a constant $C>0$ independent of $\tau$, such that
%  \begin{align*}
%   \left\|u\left(\cdot, t_n\right)-u^n\right\|_{\infty} \leq C \tau^2, \quad 1 \leq n \leq N .
%  \end{align*}
% \end{lemma}

Due to the different natures of hyperbolic and parabolic equations, subproblems \eqref{1d_hyper} and \eqref{1d_para} can be solved by different temporal and spatial discretization methods, which is one of the main advantages of operator splitting technique.

For the diffusion process \eqref{1d:Strang:diffusion1}, we employ the classical Crank-Nicolson (C-N) implicit time discretization method:
\begin{align} \label{1d:schem:Un1}
 \frac{U^{n,1}-U^{n}}{{\tau}/{2}} - \nu \frac{U_{xx}^{n,1} + U_{xx}^{n}}{2} =0,
 \quad \textrm {on} (t_n, t_{n+\frac{1}{2}}).
\end{align}

As for the hyperbolic equation \eqref {1d:Strang:hyperbolic}, we apply the second-order explicit strong stability preserving (SSP) Runge-Kutta scheme \cite{shu1988efficient} which consists of two forward Euler steps on $\left(t_n, t_{n+1}\right)$ :
\begin{numcases}{}
 \frac{U^{n,2} -U^{n,1}}{\tau} + f\left(U^{n,1} \right)_x = S(x,t^n),
 \label{1d:schem:Un2} \\
 \frac{{U}^{n,3} -U^{n,1}}{\tau} + \frac{1}{2} \big(f\left(U^{n,1}\right)_x+f\left(U^{n,2} \right)_x\big)
 = \frac{1}{2} \left(S(x,t^n) + S(x,t^{n+1}) \right).
 \label{1d:schem:Un3}
\end{numcases}
The explicit temporal discretization inevitably brings stability time step restriction \cite{gottlieb2009high}, but this is within an acceptable range, see Theorem \ref{th:1d:BP:all}.

For the diffusion equation \eqref {1d:Strang:diffusion2}, similar to the equation \eqref{1d:schem:Un1}, we give the semi-discrete scheme as:
\begin{align} \label{1d:schem:Unplus1}
 \frac{U^{n+1}-{U}^{n,3}}{{\tau}/{2}} -\nu \frac{U_{xx}^{n+1} +{U}_{xx}^{n,3}}{2} = 0,
 \quad \textrm {on} (t_{n+\frac{1}{2}}, t_{n+1}),
\end{align}
where $U^{n,1}$, $U^{n,2}$, and $U^{n,3}$ are intermediate variables, $U^{n+1}$ denotes the semi-discrete numerical solution at $t=t_{n+1}$.

Thus, combing the above approximations together, we summarize the temporal semi-discretization scheme as follows:
%{\small
\begin{numcases}{}
 \label{scheme:1d:semi_t:1}
 \dfrac{U^{n,1}-U^{n}}{{\tau}/{2}} - \nu \frac{U_{xx}^{n,1} + U_{xx}^{n}}{2} =0,
 & $\textrm {on} (t_n, t_{n+\frac{1}{2}}) $, \\
   \label{scheme:1d:semi_t:2}
  \dfrac{U^{n,2} -U^{n,1}}{\tau} + f(U^{n,1})_x = S(x,t^n),
  & $\textrm {on} (t_n, t_{n+1}) $ , \\
   \label{scheme:1d:semi_t:3}
  \dfrac{{U}^{n,3} -U^{n,1}}{\tau} + \frac{1}{2}  (f(U^{n,1})_x+f(U^{n,2})_x)
  = \dfrac{1}{2} \left(S(x,t^n) + S(x,t^{n+1}) \right),
  & $\textrm {on} (t_n, t_{n+1}) $,
    \label{scheme:1d:semi_t:4}
  % \end{cases}
 \\
 \frac{U^{n+1}-{U}^{n,3}}{{\tau}/{2}} -\nu \frac{U_{xx}^{n+1} +{U}_{xx}^{n,3}}{2} = 0,
 & $\textrm {on} (t_{n+\frac{1}{2}}, t_{n+1}) $.
\end{numcases}
This completes the one time step description of Strang splitting method.
In general, if all the solutions involved in the three-step splitting algorithm \eqref {scheme:1d:semi_t:1}--\eqref {scheme:1d:semi_t:4} are smooth, the Strang splitting method is third order at each time step and second order accurate when applied to advance the solution of \eqref{model:1d} from the initial time $t=0$ to the final time $T$ \cite{hundsdorfer2003numerical}.
%As shown in [\cite{hundsdorfer2003numerical}, Chapter IV], the above scheme (\ref{scheme:1d:full:1})--(\ref{scheme:1d:full:4}) bears a splitting error of $O(\tau^2) $ and we state the result in Lemma \ref{strang}.
%\begin{lemma} \label{strang}
% \cite{hundsdorfer2003numerical} Let $\left\{u^n\right\}_{n=0}^N$ be the numerical solution produced by the second-order splitting scheme (\ref {scheme:1d:semi_t:1})--(\ref {scheme:1d:semi_t:4}). There exists a constant $C>0$ independent of $\tau$, such that
% \begin{align*}
 %  \left\|u(\cdot, t_n\right)-u^n\right\|_{\infty} \leq C \tau^2,
 %  \quad 1 \leq n \leq N .
 % \end{align*}
%\end{lemma}

%===========================================================
\subsection{HOC spatial discretization}\label{sec_1d_space}
%===========================================================

Given a positive integer $N_x$, let $h_{x} = {L} / {N_x}$ and $ x_i= i h_x$ and let the domain $[0,\ L]$ be discretized into grids that are
described by the set $\Omega_h=\left\{x|x=\left\{x_{i}\right\}, 1 \leq i \leq N_x \right\}$. The space of periodic grid functions defined on $\Omega_h$ is denoted by $\mathcal{V}_h = \left\{u| u = \left\{u_{i}\right\}, u_{i+N_x} = u_i \right\}$.
For grid function $v \in \mathcal{V}_h $, we introduce the following difference operators:
\begin{align*}
  &  D_{\hat{x}} v_i =\frac{v_{i+1}-v_{i-1}}{2 h_x}, \quad
  \delta_x^2 v_i = \frac{v_{i+1}-2 v_i + v_{i-1}}{h_x^2}, \\
  & \malA_{x} v_{i }=
  \Big(I+\frac{h_x^2}{12} \delta_x^2\Big) v_{i }=\frac{1}{12}(v_{i-1}+10 v_{i}+v_{i+1}),\\
  & \malB_{x} v_{i}=
  \Big(I+\frac{h_x^2}{6} \delta_x^2\Big) v_{i }=\frac{1}{6}(v_{i-1}+4 v_{i}+v_{i+1}).
\end{align*}
By Taylor expansion, it's easy to get
\begin{equation} \label{HOC_AB}
 \malA_{x} v_{x x} = \delta_x^2 v+O\left(h_x^{4} \right), \quad 
  \malB_{x} v_{x} = D_{\hat{x}} v+O\left(h_x^{4} \right).
\end{equation}

Below, we use $u^n=\{u_i^n\}\in \mathcal{V}_h$ to denote the fully-discrete numerical solution at $t=t^n$, and $u^{n,1}=\{u_i^{n,1}\}$, $u^{n,2}=\{u_i^{n,2}\}$ and $u^{n,3}=\{u_i^{n,3}\}$ to represent the fully-discrete intermediate solutions.  For simplicity, denote the finite difference operator  
\begin{equation}
 \malH_{1x} = \malA_{x} -\frac{\tau}{4} \nu\delta_x^2, 
 \quad 
 \malH_{2x} = \malA_{x} +\frac{\tau}{4} \nu \delta_x^2. \label{scheme:HOC} 
\end{equation}
Substituting \eqref{HOC_AB} to \eqref{scheme:1d:semi_t:1}--\eqref{scheme:1d:semi_t:4}, a fully-discrete Strang splitting fourth-order compact finite difference scheme, named the Strang-HOC difference scheme, can be proposed as
\begin{numcases}{}
 \malH_{1x} u^{n,1}_i = \malH_{2x} u^{n}_i, \label{scheme:1d:full:1}   \\
 \malB_{x} u^{n,2}_i = \malB_{x} u^{n,1}_i -\tau D_{\hat{x}} f^{n,1}_i + \tau \malB_{x} S^{n}_i, \label{scheme:1d:full:2}   \\
 \malB_{x}{u}^{n,3}_i = \malB_{x} u^{n,1}_i -\frac{1}{2} \tau D_{\hat{x}}
 \Big(f^{n,1}_i + f^{n,2}_i \Big) + \frac{1}{2} \tau \malB_{x} \Big(S^{n}_i + S^{n+1}_i \Big), \label{scheme:1d:full:3}   \\
 \malH_{1x} u^{n+1}_i = \malH_{2x}{u}^{n,3}_i, \label{scheme:1d:full:4}
\end{numcases}
where $ f^{n,s}_i = f(u_i^{n ,s})$ for $s=1,2$, and $S^n_i =S(x_i, t^n)$. 

\begin{remark}\rm
 The algebraic systems resulting from the proposed scheme \eqref{scheme:1d:full:1}--\eqref{scheme:1d:full:4} are symmetric positive definite, well-posed and with only cyclic tridiagonal constant-coefficient matrices, which only needs to be generated once. Consequently, the proposed scheme can be implemented with a computational complexity of just $O(N_x)$ per time step. This efficient implementation makes the scheme highly suitable for practical applications and numerical simulations.
\end{remark}

\begin{remark}
    To have nonlinear stability and eliminate oscillations for shocks, a TVBM (total variation bounded in the means) limiter was introduced for the compact finite difference scheme solving scalar convection equations in \cite{shu1994TVB} and modified by \cite{li2003TVB}. Therefore, for problems such as large gradients, we will use the TVB limiter to further reduce numerical oscillations, see Example 2 in Section \ref{sec:tests}.  
\end{remark}
%============================================
\subsection{Discrete mass conservation}
%============================================

We define the discrete $L^2$ inner product as:
$$(w,\ v)= \sum_{i=1}^{N_x} h_x w_i v_i,\quad v,w \in \mathcal{V}_h.$$
It is easy to verify that the following results are valid.
\begin{lemma}\label{lem:int}
For any grid function $v\in \mathcal{V}_h $, we have
\begin{align*}
& (\mathcal{A}_x v, \ 1)=(v,\ 1),\quad (\mathcal{B}_x v,\ 1)=(v,\ 1), \\
& (\delta^2_x v,\ 1)=0,\quad (D_{\hat{x}} v,\ 1)=0.
\end{align*}
\end{lemma}

Next, we prove that the solution of the Strang-HOC difference scheme \eqref{scheme:1d:full:1}--\eqref{scheme:1d:full:4} possesses the discrete version of mass conservation.
  \begin{theorem}(Discrete Mass Conservation) \label{thm:mass:1d}
  Let $u^n\in \mathcal{V}_h$ be the solution of the Splitting-HOC difference scheme \eqref{scheme:1d:full:1}--\eqref{scheme:1d:full:4}. Then there holds
  \begin{align}\label{MassCons:D}
   (u^{n+1},\ 1) =({u}^{o} , \ 1)+ \frac{\tau}{2} \sum_{k=0}^{n} (S^{k} + S^{k+1} ,\ 1), \quad n=0,\cdots,N_t-1.
  \end{align}
  \end{theorem}
 \begin{proof}
 Taking the inner product on both sides of \eqref{scheme:1d:full:1} with 1 and using Lemma \ref{lem:int}, we obtain
  \begin{align}\label{MassCons:D1}
    (u^{n,1},\ 1) = (u^{n},\ 1).
  \end{align}
 Similarly, \eqref{scheme:1d:full:4} yields
 \begin{align}\label{MassCons:D2}
   (u^{n+1},\ 1) =({u}^{n,3} , \ 1).
 \end{align}
Taking the inner product on both sides of \eqref{scheme:1d:full:3} with 1 and using Lemma \ref{lem:int}, we have
   \begin{align}\label{MassCons:D3}
    ({u}^{n,3} , 1) = (u^{n,1}, 1) + \frac{\tau}{2} \left(S^{n} + S^{n+1} , 1\right).
  \end{align}
  Collecting \eqref{MassCons:D1}--\eqref{MassCons:D3} together, we obtain
   \begin{align}\label{MassCons:D4}
    (u^{n+1},\ 1) =({u}^{n} ,\ 1)+ \frac{\tau}{2} \left(S^{n} + S^{n+1} , 1\right),
   \end{align}
  which implies the conclusion.
\end{proof}

%=====================================
\subsection{Bound-preserving analysis}
%===========================================
In this section, we assume the exact solution satisfies
\begin{equation}
     \mathop{\min}_{x} u^{o} (x)= m \leq u(x, t) \leq M=\mathop{\max}_{x} u^{o} (x), \quad \forall t \geq 0.
\end{equation}
For simplicity, we assume the source term $S(x,t) =0$.

For convenience, we denote $\lambda=\frac{\tau}{h_x}, \mu = \frac{\tau}{h_x^{2} }, \boldsymbol u = (u_1,...,u_{N_x})^{T}$.
Under periodic boundary condition, the operator $\malH_{1x}$
 can be written correspondingly in the following form:
% $$
% \mathbf A_{x} = \frac{1}{12}
% \left(\begin{array}{ccccc}
%  10  & 1  &  &  & 1 \\
%  1  & 10  & 1  &  & \\
%     & \ddots & \ddots & \ddots & \\
%     &  & 1  & 10  & 1 \\
%  1  &  &  & 1  & 10
% \end{array}\right), %_{N_x \times N_x},
% \quad \mathbf B_{x} = \frac{1}{6}
% \left(\begin{array}{ccccc}
%  4 & 1  &  &  & 1 \\
%  1 & 4  & 1  &  & \\
%   & \ddots & \ddots & \ddots & \\
%   &  & 1  & 4  & 1 \\
%  1 &  &  & 1  & 4
% \end{array}\right), %_{N_x \times N_x}
% $$
% $$
% \mathbf D_{\hat{x}} = \frac{1}{2h_x}
% \left(\begin{array}{ccccc}
%  0  & 1  &  &  & -1 \\
%  -1  & 0  & 1  &  & \\
%    & \ddots & \ddots & \ddots & \\
%    &  & 0  & 0  & 1 \\
%  1  &  &  & -1  & 10
% \end{array}\right), %_{N_x \times N_x},
% \quad \mathbf\delta_x^2 = \frac{1}{h_x^2}
% \left(\begin{array}{ccccc}
%  -2  & 1  &  &  & 1 \\
%  1  & -2  & 1  &  & \\
%     & \ddots & \ddots & \ddots & \\
%     &  & 1  & -2  & 1 \\
%  1  &  &  & 1  & -2
% \end{array}\right), %_{N_x \times N_x},
% $$
\begin{equation}\label{def:H1x}
\mathbf H_{1x} =\frac{1}{12}
\left(\begin{array}{ccccc}
 10+6 \nu\mu & 1-3 \nu\mu &    &    & 1-3 \nu\mu \\
  1-3 \nu\mu & 10+6\nu\mu & 1-3 \nu\mu &    &    \\
     & \ddots  & \ddots   & \ddots  &    \\
      &    & 1-3 \nu\mu & 10+6\nu\mu & 1-3 \nu\mu \\
 1-3 \nu\mu &    &    & 1 -3\nu\mu & 10+6\nu\mu
\end{array}\right). %_{N_x \times N_x}.
\end{equation}
% $$
% \mathbf H_{2x} =\frac{1}{12}
% \left(\begin{array}{cccccc}
%  10-6 \nu\mu & 1+3 \nu\mu &    &    & 1+3 \nu\mu \\
%  1+3 \nu\mu & 10-6\nu\mu & 1+3 \nu\mu &    &    \\
%     & \ddots  & \ddots  & \ddots  &    \\
%      &    & 1+3 \nu\mu & 10-6\nu\mu & 1-3 \nu\mu \\
%  1+3 \nu\mu &    &    & 1+3 \nu\mu & 10-6\nu\mu
% \end{array} \right).
% $$

% On the other hand, substituting \eqref{scheme:HOC} into \eqref{scheme:1d:full:1} and \eqref{scheme:1d:full:4}, we can get the rewritten numerical scheme:
% \begin{numcases}{}
% \mathbf H_{1x} \boldsymbol u^{n,1} = \mathbf H_{2x} \boldsymbol u^{n}, \label{1d:scheme:mat:1} \\
% \mathbf B_{x} \boldsymbol u^{n,2} = \mathbf B_{x} \boldsymbol u^{n,1} -\tau \mathbf D_{\hat{x}} \boldsymbol f^{n,1} + \tau \mathbf B_{x} \boldsymbol S^{n,1},
% \label{1d:scheme:mat:2} \\
% \mathbf B_{x} u^{n,3} = \mathbf B_{x} \boldsymbol u^{n,1} -\frac{1}{2} \tau
%     \mathbf D_{\hat{x}} \left(\boldsymbol f^{n,1} + \boldsymbol f^{n,2} \right) +
%     frac{1}{2} \tau \mathbf B_{x} \left(\boldsymbol S^{n,1} + \boldsymbol S^{n,2} \right),
%  \label{1d:scheme:mat:3}  \\
% \mathbf H_{1x} \boldsymbol u^{n+1} = \mathbf H_{2x} \boldsymbol {u}^{n,3}. \label{1d:scheme:mat:4}
% \end{numcases}

Next, we will discuss the property of bound-preserving by making full use of  the coefficient matrix and the BP limiter in Algorithm \ref{alg:1d:xie} \cite{li2018A}.
 There are many equivalent definitions or characterizations
of $\mathbf{M}$-matrix, see \cite{plemmons1977m}. One convenient, sufficient, but not necessary characterization of
nonsingular $\mathbf{M}$-matrix is as follows.
\begin{lemma}\label{lem:Mmatrix}
\cite{poole1974survey,shen2021discrete} Suppose a real square matrix $\mathbf A = (a_{i,j})_{N_x \times N_x}$ is $\mathbf L$ matrix, i.e, $\mathbf A$ satisfies that $a_{ii} \geq 0 $ for each $i$ and
$a_{i,j} \leq 0$ whenever $ i \neq j$. If all the row sums of $\mathbf A$ are
non-negative and at least one row sum is positive, then $\mathbf A$ is nonsingular $\mathbf{M}$-matrix, which means $\mathbf A^{-1}$ is non-negative.
\end{lemma}

\begin{lemma} \label{lem:M-BP}
 Assume $\mathbf A \boldsymbol u^{n+1} = {\boldsymbol {u}}^{n}$, where $\mathbf A = (a_{i,j})_{N_x \times N_x}$ is a $\mathbf{M}$-matrix and tridiagonal with the sum of per row (column) equals to one. If $ {\boldsymbol {u}}^{n} \in [m,M] $, then $\boldsymbol u^{n+1}$ is bound-preserving, i.e. $\boldsymbol u^{n+1} \in [m,M] $ for $n=0, \cdots,N_t-1$.
\end{lemma}
\begin{proof}

Firstly, notice that the matrix only modifies the immediate neighbors of $\boldsymbol u_i $ and $\mathbf A^{-1}$ is a non-negative matrix from Lemma \ref{lem:Mmatrix}.
Next we will show $\boldsymbol u^{n+1}$ is bound-preserving.
Let $\mathcal C$ denotes the corresponding operator of $\mathbf A$.
 Without loss of generality, we assume
\begin{equation}
  \mathcal C u_{i}^{n+1}
= - \alpha u_{i-1}^{n+1} + (1+\alpha + \beta)u_{i}^{n+1} -\beta u_{i+1}^{n+1}
= {u}_{i}^n, \quad i =1,\cdots,N_x,
\end{equation}
where $\alpha$ and $\beta $ both are non-negative by the definition of $\mathcal C$, and $\alpha + \beta \leq 1$.

If $u_{i}^{n+1} \equiv constant$, the result is clearly true.
Otherwise, we assume exist $J$, such that $\min \limits_{i} \{u_i^{n+1} \} = u_{J}^{n+1}$ for all $n$.
Then
\begin{align*}
\begin{aligned}
 u_{J}^{n+1} & = - \alpha u_{J}^{n+1} + (1+\alpha + \beta)u_{J}^{n+1}
      - \beta u_{J}^{n+1} \\
    & \geq - \alpha u_{J-1}^{n+1} + (1+\alpha + \beta)u_{J}^{n+1}
      -\beta u_{J+1}^{n+1} \\
    & = u_{J}^{n} \geq m.
\end{aligned}
\end{align*}
Similarly, we might as well assume exist $K$, such that $\max \limits_{i} \{u_i^{n+1} \} = u_{K}^{n+1}$ for all $n$.
Then
\begin{align*}
\begin{aligned}
u^{n+1}_{K} & = - \alpha u^{n+1}_{K} + (1+\alpha + \beta)u^{n+1}_{K} -\beta u^{n+1}_{K} \\
   & \leq - \alpha u^{n+1}_{K-1} + (1+\alpha + \beta)u^{n+1}_{K}- \beta u^{n+1}_{K+1}\\
   & = { u^{n}_{K} } \leq M.
\end{aligned}
\end{align*}
 So, we can conclude that if $\boldsymbol {u}^{n} \in [m,M] $, then $\boldsymbol u^{n+1}$ is bound-preserving, i.e. $\boldsymbol u^{n+1} \in [m,M] $.

\end{proof}

%\begin{small}
\begin{algorithm} 
% \cite[Theorem 2. 7]{2018xie}
 \caption{ BP limiter for periodic data $u_{i}$ satisfying $\bar{u}_{i} \in {[m, M]}$\cite{li2018A}.}
 \label{alg:1d:xie}
 \begin{algorithmic} [1]
  
  \REQUIRE $u_{i}$ satisfies $\bar{u}_{i} =\frac{1}{c+2} \left(u_{i-1} + c u_{i}
     +u_{i+1} \right) \in[m, M], c \geq 2 $. Let $u_{0}, $ $u_{N+1}$ denote $u_{N_x}, u_{1}$, respectively.
  
  \ENSURE The output satisfies $v_{i}\in[m, M], i=1, \ldots, N_x$ and $\sum_{i=1}^{N_x}
    v_{i} =\sum_{i=1}^{N_x} u_{i}$.
  
  \STATE \textbf{Step 0}: First set $v_{i} =u_{i}, i=1, \ldots, N_x$.
    Let $v_{0}, v_{{N_x}+1}$ denote $v_{{N_x}}, v_{1}$, respectively.
  
  \STATE $\mathbf{Step I}$ : Find all the sets of class I $S_{1}, \ldots,S_{K}$
    (all local sawtooth profiles) and all the sets of class II $T_{1},\ldots, T_{K}$(consists of point values between $S_i$ and $S_{i+1}$ and two boundary
     points $u_{n_i}$ and $u_{m_{i+1} }$).
  
  \STATE \textbf{Step II}: For each $T_{j} (j=1, \ldots, K)$,
  
  \FORALL {index $i$ in $T_{j}$}
  \IF {$u_{i} <m$}
  \STATE $\quad v_{i-1} \leftarrow v_{i-1} -\frac{\left(u_{i-1} -m\right)_{+} }
    {\left(u_{i-1} -m\right)++\left(u_{i} +1-m\right)_{+} } \left(m-u_{i} \right)_{+}$
  
  \STATE $\quad v_{i+1} \leftarrow v_{i+1} -\frac{\left(u_{i+1} -m\right)_{+}}
     {\left(u_{i-1} -m\right)_{+} +\left(u_{i+1} -m\right)_{+} } \left(m-u_{i} \right)_{+}$
   
  \STATE $\quad v_{i} \leftarrow $ m
  
  \ENDIF
  \IF {$u_{i} >M $}
  \STATE $\quad v_{i-1} \leftarrow v_{i-1} +\frac{\left(M-u_{i-1} \right)_{+} }{\left(M-u_{i-1} \right)_{+} +\left(M-u_{i+1} \right)_{+} } \left(u_{i} -M\right)_{+}$
  
  \STATE $\quad v_{i+1} \leftarrow v_{i+1} +\frac{\left(M-u_{i+1} \right)_{+} }{\left(M-u_{i-1} \right)_{+} +\left(M-u_{i+1} \right)_{+} } \left(u_{i} -M\right)_{+}$
  
  \STATE $\quad v_{i} \leftarrow $ m
  
  \ENDIF
  \ENDFOR
  
  \STATE \textbf{Step III}: for each sawtooth profile $S_{j} =\big\{u_{m_{j} }, \ldots, u_{n_{j} } \big\} (j=1, \ldots, K)$, let $N_{0}$ and $N_{1}$ be the numbers of undershoot and overshoot points in $S_{j}$ respectively.
  
  \STATE Set $V_{j} =N_{1} M+N_{0} m+v_{m_{j} } +v_{n_{j} }$.
  
  \STATE Set $A_{j} =v_{m_{j} } +v_{n_{j} } +N_{1} M-\left(N_{1} +2\right) m, B_{j} =\left(N_{0} +2\right) M-v_{m_{j} } -v_{n_{j} } -N_{0} m$.
  
  \IF{$V_{j} -U_{j} >0$}
  \FOR{ $i=m_{j}, \ldots, n_{j}$}
  \STATE $\quad v_{i} \leftarrow v_{i} -\frac{v_{i} -m}{A_{j} } \left(V_{j} -U_{j} \right)$
  
  \ENDFOR
  \ELSE
  \FOR {$i=m_{j}, \ldots, n_{j}$}
  \STATE $\quad v_{i} \leftarrow v_{i} +\frac{M-v_{i} }{B_{j} } \left(U_{j} -V_{j} \right)$
  
  \ENDFOR
  \ENDIF
  
 \end{algorithmic}
\end{algorithm}
%\end{small}

The definition of (weakly) monotonic scheme is given below.
\begin{definition}
 A finite difference scheme $v^{n+1}_j = H(v^{n}_{j-s},v^{n}_{j-s+1},\dots,v^{n}_{j+s}), s>0$ is monotone if $H$ is a monotone nondecreasing function of each of its $2s +1 $ arguments. Furthermore, if the scheme $\bar{v}^{n+1} := L v^{n+1}$ and satisfies $\bar v^{n}_j$ can be written as $ H(v^{n}_{j-s},v^{n}_{j-s+1},\dots,v^{n}_{j+s})$, where $L = (l_{i,j})_{N_x\times N_x}$ is a matrix with each row (or column) summing to one and $l_{i,j}\ge 0$ for all $i,j$, then we called the scheme is weakly monotonic.
\end{definition}

 \begin{lemma} \cite{li2018A} \label{lem:xie:1d:BPlimiter}
 Assume periodic data $u_{i} (i=1, \ldots, N_x)$ satisfies 
 $$\bar{u}_{i} =\frac{1}{c+2} \left(u_{i-1} + c u_{i} +u_{i+1} \right) \in[m, M], \quad c \geq 2$$ 
 for all $i=1, \ldots, N_x$ with $u_{0} :=u_{N_x}$ and $u_{N_x+1} :=u_{1}$; then the output of Algorithm \ref{alg:1d:xie} satisfies $\sum_{i=1}^{N_x} v_{i} =\sum_{i=1}^{N_x} u_{i}$ and $v_{i} \in[m, M]$ for all $i$.
\end{lemma}

\begin{remark} \label{rem:order}
 % The limiter described in Algorithm \ref{alg:1d:xie} are high order accurate limiters in the same sense \cite{li2018A}.
 As stated in \cite{li2018A}, the BP limiter modifies $u_i $ by $ O(h_x^4)$, indicating that the limiter does not affect the accuracy of our scheme. Furthermore, according to Lemma \ref{lem:xie:1d:BPlimiter}, the limiter preserves the mass conservation of the original numerical solution.
\end{remark}

%%%%%%--------------------------------

The proof of the bound-preserving property of the Strang-HOC difference scheme \eqref{scheme:1d:full:1}--\eqref{scheme:1d:full:4} is divided into the following lemmas. Firstly, we consider the diffusion process \eqref{scheme:1d:full:1} in the first half-time step $[t_n,t_{n+1/2}]$.
%----------------------------
\begin{lemma} \label{lem:1d:diffusion}
 Under the constraint $\nu\mu \leq \frac{5}{3}$, if $u^{n} \in[m, M]$, then $ {u}^{n,1 }$ computed by the scheme \eqref{scheme:1d:full:1} with the BP limiter satisfies $ m \le {u}^{n,1 } \le M$.
\end{lemma}
\begin{proof}
  For the scheme \eqref{scheme:1d:full:1}, incorporating the definition of $\malH_{1x}$, the following weak monotonicity holds under the condition $\nu \mu \leq \frac{5}{3}$ :
 \begin{align*}
  \malH_{1x} u_i^{n}
   & = \frac{1}{12} \big((1+3 \nu \mu) u^{n}_{i-1} + (10-6 \nu \mu)u^{n}_{i} +
        (1+3 \nu \mu)u^n_{i+1} \big) \notag \\
   & = \hat K\left(u_{i-1}^{n}, u_{i}^{n}, u_{i+1}^{n} \right)
   =\hat K(\uparrow, \uparrow, \uparrow).
 \end{align*}
 where $\uparrow$ denotes that the partial derivative with respect to the corresponding argument is non-negative \cite{li2018A}.
 % and $\malH_{1x}$ denote represents the operator form of matrix $\mathbf H_{1x}$, i.e, $$\malH_{1x} v_i = \frac{1}{12} \left((1+3 \nu \mu) v_{i-1} + (10-6 \nu \mu)v_{i} +
 %       (1+3 \nu \mu)v_{i+1} \right). $$
Therefore, $m \leq u_{i}^{n} \leq M$ implies $m=\hat K(m, m, m) \leq \malH_{1x} u_i^{n,1} \leq$ $\hat K(M, M, M)=M$ for all $ 1\leq i \leq N_{x}$.
 
% To insure $ m \leq u^{n,1} \leq M$, on the one hand, we can guarantee
% that $\mathbf H_{1x}$ is an $\mathrm{M}$ matrix by applying {Lemma \ref{lem:M-BP}}.
% Clearly, $\mathbf H_{1x}$ is a symmetric and strictly diagonally dominant matrixwith positive diagonal elements, making it positive definite. Thus, according to Lemma \ref{lem:Mmatrix}, we only need to guarantee $\mathbf H_{1x}$ is an $\mathrm L$ matrix, which requires all off-diagonal elements to be nonpositive. Since the diffusion coefficient $\nu$ is positive, this condition converts to 
% $1-3 \nu \mu \leq 0, i.e. \nu \mu \geq \frac{1}{3}$.

% On the other hand, if $\nu \mu < \frac{1}{3}$, although $\mathbf H_{1x}$ is not an $\mathbf{M}$-matrix, it's still a convex combination of $\boldsymbol u^{n ,1}$ and satisfies $\malH_{1x} u^{n,1 }_i =\frac{1}{c+2} \left(u_{i-1}^n + c u_{i}^n +u_{i+1}^n \right) \in[m, M] $, where $c= \frac{10+6 \nu \mu}{1-3 \nu \mu}\geq 2$.
% Through the post-processing of the BP limiter, if $ 0 \leq \nu \mu <\frac{1}{3}$,
% $\malH_{1x} u^{n,1 }_i $ computed by \eqref{scheme:1d:full:1} still satisfies $ m \le {u}^{n,1 }_i \le M$ for all $i$.

To insure $ m \leq u^{n,1} \leq M$, it can be divided into the following two cases. 
\begin{itemize}
 
 \item Case 1: On employing the definition of $\mathrm{M}$ matrix in Lemma \ref{lem:M-BP}, we can obtain $ m \leq u^{n,1} \leq M$ by ensuring that $\mathbf H_{1x}$ is $\mathrm{M}$ matrix.
 Clearly, $\mathbf H_{1x}$ defined in \eqref{def:H1x} is a symmetric and strictly diagonally dominant matrix with positive diagonal elements, making it positive definite. Thus, according to Lemma \ref{lem:Mmatrix}, we only need to guarantee $\mathbf H_{1x}$ is $\mathrm L$ matrix, which requires all off-diagonal elements to be nonpositive. Since the diffusion coefficient $\nu$ is positive, this condition converts to 
 $1-3 \nu \mu \leq 0, i.e. \nu \mu \geq \frac{1}{3}$. 

 \item Case 2:  If $\nu \mu < \frac{1}{3}$, although $\mathbf H_{1x}$ is not $\mathbf{M}$-matrix, it's still a convex combination of $\boldsymbol u^{n ,1}$ and satisfies $\malH_{1x} u^{n,1 }_i =\frac{1}{c+2} \left(u_{i-1}^n + c u_{i}^n +u_{i+1}^n \right) \in[m, M] $, where $c= \frac{10+6 \nu \mu}{1-3 \nu \mu}\geq 2$. Thus, we can achieve it through combining \eqref{scheme:1d:full:1} with the BP limiter.
 Then through the post-processing of the BP limiter, if $ 0 < \nu \mu <\frac{1}{3}$,
 $\malH_{1x} u^{n,1 }_i $ computed by \eqref{scheme:1d:full:1} still satisfies $ m \le {u}^{n,1 }_i \le M$ for all $i$.
 \end{itemize}
In summary, for $i=1,\cdots,N_x$, if $\frac{1}{3} \leq \nu \mu \leq \frac{5}{3}$ holds, ${u}^{n,1 }_i $ computed by the original scheme \eqref{scheme:1d:full:1} is bound-preserving; if $\nu \mu < \frac{1}{3}$ holds, we can achieve that ${u}^{n,1 }_i \in [m,M]$ through using the scheme \eqref{scheme:1d:full:1} with the BP limiter. Then we complete the proof.
\end{proof}

By the symmetry of \eqref{scheme:1d:full:4} and \eqref{scheme:1d:full:1}, We have similar conclusions for \eqref{scheme:1d:full:4} under the condition that $u^{n,3} \in[m, M]$.
\begin{corollary} \label{cor:1d:diffusion}
 Under the constraint $\nu \mu \leq \frac{5}{3}$, if $ {u}^{n,3} \in[m, M]$, then $ {u}^{n+1}$ computed by the scheme \eqref{scheme:1d:full:4} with the BP limiter
 satisfies $ m \le {u}^{n+1} \le M$.
\end{corollary}

Next, we consider the second subprocess: the hyperbolic equation \eqref{scheme:1d:full:2}-- \eqref{scheme:1d:full:3} for the entire step from $t^n$ to $t^{n+1}$.
%%-------------------------------
\begin{lemma} \label{lem:1d:convection:1}
  Under the CFL constraint $\lambda \max \limits_{u}\left|f^{\prime}(u)\right| \leq \frac{1}{3}$, if $u^{n ,1} \in[m, M]$, then $u^{n,2 }$ computed by the scheme \eqref{scheme:1d:full:2} with the BP limiter 
 satisfies $ m \leq u^{n,2 } \leq M $.
\end{lemma}
\begin{proof}
 For the scheme \eqref{scheme:1d:full:2}, the following weak monotonicity holds under the condition of
 $\lambda \max \limits_{u}\left|f^{\prime}(u)\right| \leq \frac{1}{3}$,
\begin{equation*}
 \begin{aligned}
  \malB_{x} u_i^{n,2} & = \malB_{x} u_i^{n,1} -\tau D_{\hat{x}} f(u_i^{n,1}) \\
      & = \Big (\frac{1}{6} u_{i-1}^{n,1}
      + \frac{\lambda}{2} f(u_{i-1}^{n ,1}) \Big)
     + \frac{4}{6} u_{i}^{n,1} + \Big(\frac{1}{6} u_{i+1}^{n,1} -\frac{\lambda}{2} f (u_{i+1}^{n,1 }) \Big)\\
      & = \tilde K\big(u_{i-1}^{n ,1}, u_{i}^{n,1}, u_{i+1}^{n ,1} \big)
      = \tilde K(\uparrow, \uparrow, \uparrow).
  \end{aligned}
\end{equation*}
 Therefore, $m \leq u_{i}^{n,1 } \leq M$ implies
 $m = \tilde K(m, m, m) \leq \malB_{x} u_i^{n,2 } \leq$ $\tilde K(M, M, M) = M $.
Similarly, by using the BP limiter with $c=4$, we can ensure the bound-preserving property of $u_i^{n,2 }$.
\end{proof}

 \begin{lemma} \label{lem:1d:convection:2}
 Under the CFL constraint $\lambda \max \limits_{u}\left|f^{\prime}(u)\right| \leq \frac{1}{3}$, if $\left\{u^{n,2}, u^{n,3} \right\} \in[m, M]$, then
 $u^{n,3}$ computed by the scheme \eqref{scheme:1d:full:3} with the BP limiter satisfies $ m \leq u^{n,3}\leq M $.
 \end{lemma}
\begin{proof}
% Firstly, we prove
Substituting \eqref{scheme:1d:full:2} into \eqref{scheme:1d:full:3}, we get
 $$
 \begin{aligned}
 \qquad \malB_{x} u_i^{n,3}
  = & \malB_{x} u_i^{n,1 } -\frac{1}{2} \tau D_{\hat{x}}
  \Big(f\big(u_i^{n,1 }\big)+f\big(u_i^{n,2} \big)\Big) \\
  % = & \frac {1}{2}
  %  \left[ \left(\malB_{x} u_i^{n,1 } -\tau D_{\hat{x}} f(u_i^{n,1 }) \right)
  % +   \left(\malB_{x} u_i^{n ,1} - \tau D_{\hat{x}} f(u_i^{n,2})\right)
  %  \right]\\
   = & \frac {1}{2} \malB_{x} u_i^{n,1} 
    + \frac {1}{2}\left(
     \malB_{x} u^{n,2}- \tau D_{\hat{x}} f(u^{n,2})
     \right)\\
  % = & \frac {1}{2}
  %  \left[ \malB_{x} u^{n,1 }
  %     + \malB_{x} u^{n,2} - \tau D_{\hat{x}} f(u^{n,2})
  %  \right]\\
  = & \frac{1}{2}
   \Big(
    \frac{1}{6} u_{i-1}^{n,1} + \frac{4}{6} u_{i}^{n,1} + \frac{1}{6} u_{i+1}^{n,1}
   \Big) \notag \\
   & +\frac{1}{2}
   \left[
     \left(\frac{1}{6} u_{i-1}^{n,2} + \frac{\lambda}{2} f(u_{i-1}^{n,2}) \right)
     + \frac{4}{6}u_{i}^{n,2} +
     \left(\frac{1}{6} u_{i+1}^{n,2} - \frac{\lambda}{2} f(u_{i+1}^{n,2}) \right)
    \right ].
 \end{aligned}
 $$
% where $\bar{w}_{i}^{n+1} =\mathbf B_x \tilde{u}^{n+1}$.
 So under the
 CFL condition
 $\lambda \max \limits_{u}\left|f^{\prime}(u)\right| \leq \frac{1}{3}$,
 we can insure the right hand of the equation satisfies the weak monotonicity, i.e.,
$  \malB_{x} u_i^{n,3} = \bar K(\uparrow, \uparrow, \uparrow).$
 Therefore, if $u_i^{n,1 }$ and $u_i^{n,2}  \in[m, M]$, by using the BP limiter, we have $ m \leq u_i^{n,3} \leq M $.

\end{proof}

% Then we only need to insure $u^{n,3}$ and $u^{n,2}$ computed by the scheme \eqref{scheme:1d:full:2} and \eqref{scheme:1d:full:3} satisfy $ m \leq \{u^{n,2}, u^{n,3} \} \leq M $.
% We give the following lemma.

Combining Lemmas \ref{lem:1d:diffusion}--\ref{lem:1d:convection:2} and Corollary \ref{cor:1d:diffusion}, we have the following theorem.
\begin{theorem} \label{th:1d:BP:all}
  Under the constraints
  $\nu \mu \leq \frac{5}{3}, \lambda \max \limits_{u}\left|f^{\prime}(u)\right| \leq \frac{1}{3}$,
  if $\boldsymbol u^{n} \in[m, M]$, then $\boldsymbol u^{n+1}$ computed by the scheme \eqref{scheme:1d:full:1}--\eqref{scheme:1d:full:4} with the BP limiter
  satisfies $ m \leq \boldsymbol u^{n+1} \leq M $.
 % after applying the limiter in Algorithm \ref{alg:1d:xie} to $\malA_y \boldsymbol u^{n,1 }, \malB_x \boldsymbol u^{n,2 }, \malB_x \boldsymbol u^{n,3 }$ respectively.
  % In particular, if CFL condition $\lambda \max \limits_{u}\left|f^{\prime}(u)\right| \leq \frac{1}{3}$ satisfied, the scheme is stable.
\end{theorem}
% \begin{proof}
 
%  From the Lemmas \ref{lem:1d:diffusion}--\ref{lem:1d:convection:2} and Corollary \ref{cor:1d:diffusion},  under the constraints $\nu \mu \leq \frac{5}{3}, \lambda \max \limits_{u}\left|f^{\prime}(u)\right| \leq \frac{1}{3}, $  we have $\malB_x u_i^{n,2 }, \malB_x u_i^{n,3 } \in [m, M]$ if $\boldsymbol u^{n,1 } ,\boldsymbol u^{n,2 } \in [m, M]$.
%  Notice that
%  $$
%  \begin{aligned}
%   \malB_x u_i^{n,2 } =\frac{1}{c+2} \left(u_{i-1}^{n,2} +c u_{i}^{n,2} + u_{i+1}^{n,2} \right), \quad \text{where } c=4,\\
%  %
%  \malB_x u_i^{n,3 } =\frac{1}{c+2} \left(u^{n,3}_{i-1} +c u^{n,3}_{i} + u^{n,3}_{i+1} \right), \quad \text{where } c=4,
%  \end{aligned}
%  $$
%  all discussions in Lemma \ref{lem:xie:1d:BPlimiter} are still valid. Then the output values $\boldsymbol u^{n,3}$ and $\boldsymbol {u}^{n,3}$ from Algorithm \ref{alg:1d:xie} are in the range $[m, M]$.
 
 % As for the stability, it is known from Theorem \ref{lem:1d:convection:1} and \ref{lem:1d:convection:2} that when the $\lambda \max \limits_{u}\left|f^{\prime}(u)\right| \leq \frac{1}{3}$ is satisfied, the first and last substeps \eqref{2d:scheme:full2} and \eqref{scheme:1d:full:3} are monotonic. For the treatment of intermediate diffusion terms, C-N method guarantees unconditional stability. Therefore, in summary, the scheme is stable as long as it satisfies the CFL condition $\lambda \max \limits_{u}\left|f^{\prime}(u)\right| \leq \frac{1}{3}$.
 
% \end{proof}
%
\begin{remark} \label{rem:CFL:h_tau:1D}
    In fact, when $\nu$ is large enough so that $h < \frac{\nu}{5\max \limits_{u}\left|f^{\prime}(u)\right|}$, then we only need $\tau < \frac{\nu}{15\max \limits_{u}\left|f^{\prime}(u)\right|}$. 
    On the other hand, when $\nu$ is small enough so that $h \ge\frac{\nu}{5\max \limits_{u}\left|f^{\prime}(u)\right|}$, we only need $\tau \le \frac{h}{3\max \limits_{u}\left|f^{\prime}(u)\right|}$. In this paper, since we mainly consider convection-dominated  problems, i.e., the diffusion coefficient $\nu$ is very small, in such case, we only need to satisfy  $\tau \le \frac{h}{3\max \limits_{u}\left|f^{\prime}(u)\right|}$.
\end{remark}

\begin{remark} \label{rem:sufficient}
The above theorem provides only a sufficient condition for the bound-preserving property of the scheme. In particular, even without the assistance of the limiter Algorithm \ref{alg:1d:xie}, the scheme itself remains bound-preserving in certain cases. See Section \ref{sec:tests} for numerical examples.

\end{remark}

We now summarize the above Strang-HOC scheme \eqref{scheme:1d:full:1}--\eqref{scheme:1d:full:4} with the BP limiter into Algorithm \ref{alg:1d:self}.
%----------------------------------
\begin{algorithm} [ht]
 \caption{Strang-HOC algorithm with the BP limiter}
 \label{alg:1d:self}
 \begin{algorithmic} [1]
  \FOR{$n=0$ to $N_t-1$}
  
  \STATE \textbf{Step 1:} Compute $\boldsymbol u^{n,1}$ by \eqref{scheme:1d:full:1} with initial value $u^n$ on $(t_n,t_{n+ {1}/{2}})$.
  
  \IF {$\nu \mu < \frac{1}{3}$}
  \STATE Applying the limiter in Algorithm \ref{alg:1d:xie} to $\malH_{1x} u^{n,1}$ to obtain $u^{n,1}$ which is bound-preserving.
  \ENDIF
  
  \STATE \textbf{Step 2:} Compute $\malB_x \boldsymbol u^{n,2}$ by \eqref{scheme:1d:full:2} with initial value $u^{n,1}$  on $(t_n,t_{n+1})$.
  
  \STATE \textbf{Step 3:} Applying the limiter in Algorithm \ref{alg:1d:xie} to $\malB_x u^{n,2}$ to obtain $u^{n,2}$ which is bound-preserving.

  \STATE \textbf{Step 4:} Compute $\malB_x u^{n,3}$ by \eqref{scheme:1d:full:3} with initial value $u^{n,1}$ and $u^{n,2}$ on $(t_n,t_{n+1})$.
  
  \STATE \textbf{Step 5:} Applying the limiter in Algorithm \ref{alg:1d:xie} to $\malB_x u^{n,3}$ to obtain $u^{n,3}$ which is bound-preserving.
  
  \STATE \textbf{Step 6:} Compute $\malH_{1x} u^{n+1}$ by \eqref{scheme:1d:full:4} with initial value $u^{n,3}$ to obtain  on $(t_{n+ {1}/{2}},t_{n+1})$.
  
  \IF {$\nu \mu < \frac{1}{3}$}
  \STATE Applying the limiter in Algorithm \ref{alg:1d:xie} to $\malH_{1x} u^{n+1}$ to obtain $u^{n+1}$ which is bound-preserving.
  
  \ENDIF
  \ENDFOR
 \end{algorithmic}
\end{algorithm}

 \section{ Extension to two-dimensional problems}
%%%%%%%%%%%%%%%%%%%%%%%%%%%%%%%%%%%%%%%%%%%%%%%%

 % In the previous section, we have already constructed the Strang-HOC difference scheme for solving the one-dimensional nonlinear convection diffusion equation.
 \qquad In this section, we extend our scheme to two-dimensional equation with periodic boundary as follows:
 \begin{equation}
\begin{aligned} \label{mod:2d}
 & u_t + f(u)_ x + g(u)_y - \nu (u_{xx}+ u_{yy})
=S(x,y,t), \quad 0\leq x,y \leq L,\ 0\leq t \leq T, \\
 & u(\bm x, 0)=u^o(\bm x), \quad 0\leq x,y \leq L,\ 0\leq t \leq T.
 \end{aligned}
 \end{equation}
 where $\bm x = (x,y)$ and $\nu $ is a
 small positive constant compared with $\mathop{\min}\limits_{u} \{|f'(u)|,|g'(u)|\}$.
% Assume $f(u)$, $ g(u) $ and $S(x,t)$ are well-defined smooth functions for any $u \in$ $[ m, M]$,
% where
% ${\mathop{\min} \limits_{x, y} } u^o(x, y)$ and $ M=\mathop{\max} \limits_{x, y} u_{0} (x, y). $
% And its exact solution satisfies
% $$
% \mathop{\min}_{x, y} u^{o} (x, y)= m \leq u(x, y, t) \leq M=\mathop{\max}_{x, y} u^{o} (x, y) \quad \forall t \geq 0.
% $$

 % We now provide the description of fully discrete difference scheme for solving \eqref{mod:2d}.
%We firstly consider the nonlinear convection and diffusion terms by operator splitting.

\subsection{Strang splitting time discretization}

Similar to one dimension case, we define the time step $\tau= {T}/{N_t}$ and $t_n = n\tau(n= 0,1,\dots,N_t)$.
By combining Strang splitting method, \eqref{mod:2d} is reduced to updating two nonlinear systems, specifically the first and third equations of \eqref{2d:Strang:diffusion1} and \eqref{2d:Strang:hyperbolic},
and one diffusion equation \eqref{2d:Strang:diffusion2}) as follows:
\begin{numcases}{}
  u^{*}_{t} - \nu (u_{xx}+u_{yy}) =0,& $t\in [t_n, t_{n+1/2}]; \quad u^{*}(\bm x,{t_n})=u^n(\bm x)$,
  \label{2d:Strang:diffusion1} \\
  u^{**}_{t} + f_x(u)+ g_y(u)=S(x,y,t), & $t\in [t_n, t_{n+1}]; \quad u^{**}(\bm x,t_n)=u^{*}(\bm x,t_{n+{1}/{2}})$,
  \label{2d:Strang:hyperbolic} \\
 u^{***}_{t} - \nu (u_{xx}+u_{yy}) =0, & $t\in [t_{n+1/2}, t_{n+1}]; \quad u^{***}(\bm x,t_{n+{1}/{2}})=u^{**}(\bm x,t_{n+1})$,
  \label{2d:Strang:diffusion2}
 \end{numcases}
 for $n = 0,1,\cdots, N_t-1$ with $u^{n+1}(\bm x) = u^{***}(\bm x,t_{n+1})$.

To further reduce computational cost for the diffusion terms, we apply a Douglas-type ADI scheme \cite{hundsdorfer2003numerical,douglas1962alternating} for time discretization, which is second-order accurate and easily generalizable to high-dimensional problems. By combining the ADI method, high-dimensional problems can be transformed into the calculation of several one-dimensional problems.
In the first half-step $\left(t_{n+ {1}/{2}}, t_{n+1} \right)$, we have
\begin{numcases}{}
\frac{U^{n ,1} -U^{n} }{\frac{\tau}{2} } - \nu \frac{U_{xx}^{n,1 } +U_{xx}^{n} }{2} -\nu U_{yy}^{n} =0, \notag\\
 \frac{U^{n ,2} -U^{n ,1} }{\frac{\tau}{2} } -\nu \frac{U_{yy}^{n ,2} -U_{yy}^{n} }{2} =0. \notag
\end{numcases}
The scheme for the remaining half-step $\left(t_{n+ {1}/{2}}, t_{n+1} \right)$ can be similarly given.

 For convection items, we apply the improved Euler scheme  as time discretization on entire interval $\left(t_n, t_{n+1} \right) $, then we have
 \begin{align*} 
 \left\{  
  \begin{aligned}
 &\frac{U^{n,3} -U^{n,2 } }{\tau} + f\left(U^{n ,2} \right)_x+g\left(U^{n ,2} \right)_y = S(x,y,t^n), \\
 & \frac{{U}^{n,4} -U^{n ,2 } }{\tau} +\frac{1}{2} \left(f\left(U^{n,2 } \right)_x+f\left(U^{n ,3} \right)_x\right) +
 \frac{1}{2} \left(g\left(U^{n ,2} \right)_y+g\left(U^{n ,3} \right)_y\right) \\
 & \quad = \frac{1}{2} \left(S(x,y, t^{n}) + S(x,y,t^{n+1}) \right).
  \end{aligned}
 \right.
 \end{align*} 

Then we have the semi-discrete scheme as follows :
 \begin{numcases}{}
 \frac{U^{n ,1} -U^{n} }{\frac{\tau}{2} } - \nu \frac{U_{xx}^{n,1 } +U_{xx}^{n} }{2} -\nu U_{yy}^{n} =0, \\
 \frac{U^{n ,2} -U^{n ,1} }{\frac{\tau}{2} } -\nu \frac{ U_{yy}^{n ,2} -U_{yy}^{n} }{2} =0, \\
 \frac{U^{n,3} -U^{n,2 } }{\tau} + f\left(U^{n ,2} \right)_x+g\left(U^{n ,2} \right)_y = S(x,y,t^n), \\
 \frac{{U}^{n,4} -U^{n ,2 } }{\tau} +\frac{1}{2} \big(f\left(U^{n,2 } \right)_x+f\left(U^{n ,3} \right)_x\big) +
 \frac{1}{2} \big(g\left(U^{n ,2} \right)_y+g\left(U^{n ,3} \right)_y\big) \notag \\
\qquad \quad = \frac{1}{2} \left(S(x,y,t^n) + S(x,y,t^{n+1}) \right), \\
 \frac{U^{n ,5} -U^{n,4} }{\frac{\tau}{2} } - \nu \frac{U_{xx}^{n,5 } +U_{xx}^{n,4} }{2} -\nu U_{yy}^{n,4} =0, \\
 \frac{U^{n+1} -U^{n ,5} }{\frac{\tau}{2} } -\nu \frac{ U_{yy}^{n+1} -U_{yy}^{n,5} }{2} =0,
 \end{numcases}
where $U^{n,s}$ for $s=1,\cdots,5$ are intermediate variables, $U^{n+1}$ denotes the semi-discrete numerical solution at $t=t_{n+1}$.

%%%---------------------------------------
\subsection{HOC spatial discretization}
%%%---------------------------------------

 Similar to the one-dimensional case, we introduce the following notations:
 $$\Omega_h=\left\{(x_i,y_j) \mid 0 \leq i \leq N_x, 0 \leq j \leq N_y \right\}, ~~h_x = L/N_x, ~~h_y = L/N_y, $$
$$\mathcal{Z}_h = \left\{u = \{u_{i,j}\}| u_{i+N_x,j} = u_{i,j}, u_{i,j+N_y} = u_{i,j} \right\},$$
where $N_x$ and $N_y$ are given positive integers.
For any grid function $v \in \mathcal{Z}_h $, we define the difference operators 
 $\mathcal D_{\hat{z}} v$, $\delta_z^2 v$, and compact operators $\malA_{z}, \malB_{z}$ and $\malH_{1z}$ , for $z=x$ and $y$, accordingly. In addition, 
we rewrite $\mathbf H_{1z}=\mathbf A_z - \frac{\tau}{4}\nu \delta_z^2$ in the corresponding matrix form as $\malH_{1z}$ for $z=x,y$, similar to the one-dimensional case. Denote $u:=(u_{i,j})_{N_x\times N_y}$ and $\boldsymbol{u^{n,s}}:=(u_{i,j}^{n,s})_{N_x\times N_y}$ for $s=1,\cdots,5$.

Applying the fourth order compact finite difference scheme for spatial discretization, we have the Strang-ADI-HOC  scheme as follows:
%\begin{small}
\begin{numcases}{}
 \malH_{1 x} \malA_{y}  u_{i,j}^{n,1}
 = \left(\malA_{x} \malA_{y} +\frac{\nu \tau}{4} \malA_{y}
 \mathcal \delta_x^2 + \frac{\nu \tau}{2} \malA_{x} \delta_y^2 \right)u_{i,j}^{n},
 \label{2d:scheme:full1} \\
 \malH_{1 y} u_{i,j}^{n,2}
 = \malA_{y} u_{i,j}^{n,1} -\frac{\nu \tau }{4} \mathcal \delta_y^2 u_{i,j}^{n},
 \label{2d:scheme:full2} \\
 \malB_{x} \malB_{y} u_{i,j}^{n,3} = \malB_{x} \malB_{y} u_{i,j}^{n,2} -\tau \malB_{y} D_{\hat{x}} f_{i,j}^{n,2}-\tau \malB_{x} D_{\hat{y}} g_{i,j}^{n,2} + {\tau} \malB_{x} \malB_{y} S^{n}_{i,j},
 \label{2d:scheme:full3} \\
 \malB_{x} \malB_{y} u_{i,j}^{n,4}
 = \malB_{x} \malB_{y} u_{i,j}^{n,2} -\frac{\tau}{2} \malB_{y} \mathcal
 D_{\hat{x}} \left(f_{i,j}^{n,2} + f_{i,j}^{n,3} \right)-\frac{\tau}{2} \malB_{x} D_{\hat{y}} \left(g_{i,j}^{n,2} + g_{i,j}^{n,3} \right) \notag \\
 \qquad \qquad \quad + \frac{\tau}{2} \malB_{x} \malB_{y} \left(S^{n}_{i,j} + S^{n+1}_{i,j} \right),
 % \text {on} \left(t_n, t_{n+1} \right)
 \label{2d:scheme:full4}
 \\
 \malH_{1 x} \malA_{y} u_{i,j}^{n,5}
 =
 \left(\malA_{x} \malA_{y} +\frac{\nu \tau}{4} \malA_{y} \mathcal \delta_x^2 +
 \frac{\nu \tau}{2} \malA_{x}  \delta_y^2
 \right)u_{i,j}^{n,4},
 % \quad \text {on} \left(t_n, t_{n+1 / 2} \right)
 \label{2d:scheme:full5}
 \\
 \malH_{1 y} u_{i,j}^{n+1}
 = \malA_{y} u_{i,j}^{n,5} -\frac{\nu \tau }{4}  \delta_y^2 u_{i,j}^{n,4}, %\quad \text {on} \left(t_n, t_{n+1 / 2} \right)
 \label{2d:scheme:full6}
\end{numcases}
%\end{small}
where $f^{n,s}_{i,j}:=f(u^{n,s}_{i,j})$, $g^{n,s}_{i,j}:=g(u^{n,s}_{i,j})$ for $s=1,\cdots,5$ and $S^{k}_{i,j}:=S(x_i,y_j,t^{k})$.  
 
\begin{remark}\rm
 When solving \eqref{2d:scheme:full1}, we observe that we only need to solve $\malA_{x} u_{i,j}^{n,1},\malA_{x} u_{i,j}^{n,5}$, instead of $u_{i,j}^{n,1}$ and $u_{i,j}^{n,5}$, by directly incorporating \eqref{2d:scheme:full2}. Consequently, the resulting algebraic systems from \eqref{2d:scheme:full1}--\eqref{2d:scheme:full2} involve only cyclic tridiagonal constant-coefficient matrices, which need to be generated only once. Furthermore, the matrices in \eqref{2d:scheme:full5}--\eqref{2d:scheme:full6} are identical to those in \eqref{2d:scheme:full1}--\eqref{2d:scheme:full2}.

\end{remark}

\subsection{Discrete mass conservation}

We define the discrete $L^2$ inner product as:
$$(w,\ v)= \sum_{i=1}^{N_x} \sum_{j=1}^{N_y} h_x h_y w_i v_j,\quad v,w \in \mathcal{Z}_h.$$
Below, we use $u^n=\{u_{i,j}^n\}\in \mathcal{Z}_h$ to denote the fully-discrete numerical solution at $t=t^n$, and $u^{n,s}=\{u_{i,j}^{n,s}\}$ for $s=1,\cdots,5$ to represent the fully-discrete intermediate solutions. Similar to one-dimensional case, we have the following results.
\begin{lemma}\label{lem:int:2d}
For any grid function $v\in \mathcal{Z}_h $, for $z = x,y$, we have
\begin{align*}
& (\mathcal{A}_z v, \ 1)=(v,\ 1),\quad (\mathcal{B}_z v,\ 1)=(v,\ 1), \\
& (\delta^2_z v,\ 1)=0,\quad (D_{\hat{z}} v,\ 1)=0.
\end{align*}
\end{lemma}

\begin{theorem}
Let $u^n \in \mathcal{Z}_h$ be the solution of the Strang-ADI-HOC difference scheme \eqref{2d:scheme:full1}--\eqref{2d:scheme:full4}. Then there holds
 \begin{align*}
  (u^{n+1},1) = (u^{o},1) + \frac{\tau}{2} \sum_{k=0}^{n} (S^{k} + S^{k+1} ,\ 1) .
 \end{align*}
\end{theorem}
\begin{proof}

  Taking the inner product on both sides of \eqref{2d:scheme:full1}--\eqref{2d:scheme:full2} with 1 and using Lemma \ref{lem:int:2d}, we obtain
  \begin{align}\label{2d:MassCons:D1}
    (u^{n,1},\ 1) = (u^{n},\ 1), \quad (u^{n,2},\ 1) = (u^{n,1},\ 1),
  \end{align}
 which means $$(u^{n,2},\ 1) = (u^{n},\ 1).$$
 Similarly, \eqref{2d:scheme:full5}--\eqref{2d:scheme:full6} yields
 \begin{align}\label{2d:MassCons:D2}
   (u^{n+1},\ 1) =({u}^{n,4} , \ 1).
 \end{align}
Taking the inner product on both sides of \eqref{2d:scheme:full3} with 1 and using Lemma \ref{lem:int:2d}, we have
   \begin{align}\label{2d:MassCons:D3}
    ({u}^{n,3} , 1) = (u^{n,1}, 1) + \frac{\tau}{2} \left(S^{n} + S^{n+1} , 1\right).
  \end{align}
  Collecting \eqref{2d:MassCons:D1}--\eqref{2d:MassCons:D3} together we obtain
   \begin{align}\label{2d:MassCons:D4}
    (u^{n+1},\ 1) =({u}^{n} ,\ 1)+ \frac{\tau}{2} \left(S^{n} + S^{n+1} , 1\right),
   \end{align}
  which completes the proof.
\end{proof}

%%%---------------------------------------
\subsection{Bound-preserving analysis}
%%%---------------------------------------

%We abuse the notation by using $H_{1 x} u_{i,j}$ to denote the $(i, j)$ entry of $H_{1 x}{u}$.
%For convenience, we abbreviate $f(u^{n,i}) $ as $f^{n,i}$ and $S(u^{n,i}) $ as $S^{n,i}$, respectively. In addition, we denote $\mu_1= \frac {\tau }{h_x^2}, \mu_2= \frac {\tau }{h_ y^2}, \lambda = \frac {\tau }{h_ x}, \lambda_2 = \frac {\tau }{h_ y}$.
Without loss of generality, we assume $S(x,y,t)=0$. Assume the exact solution satisfies
$$
 \mathop{\min}_{x,y} u^{o} (x,y)= m \leq u(x,y, t) \leq M=\mathop{\max}_{x,y} u^{o} (x), \quad \forall t \geq 0.
$$
For convenience, we introduce
\begin{align} \label{def:W}
 W^{n,s} =
  \left(\begin{array}{lll}
  u_{i-1, j+1}^{n,s} & u_{i, j+1}^{n,s} & u_{i+1, j+1}^{n,s} \\
  u_{i-1, j}^{n,s} & u_{i, j}^{n,s} & u_{i+1, j}^{n,s} \\
  u_{i-1, j-1}^{n,s} & u_{i, j-1}^{n,s} & u_{i+1, j-1}^{n,s}
 \end{array} \right), \quad s=0,\cdots,5,
\end{align}
where we denote $W^{n,0}:=W^n$. In addition, we define $\lambda_1 =\frac{\tau}{h_x}, \lambda_2=\frac{\tau}{h_y}$ and $\mu_1 = \frac{\tau}{h_x^2}, \mu_2 =\frac{\tau}{h_y^2}$. 
Then we consider the bound-preserving property of the above scheme.
Firstly, we consider the scheme \eqref{2d:scheme:full1}--\eqref{2d:scheme:full2}.

\begin{lemma} \label{lem:2d:diffusion}
 Under the constraint ${\nu \mu_1=\nu \mu_2 \leq \frac{5}{3}}$, if $u^{n} \in[m, M]$, then $  u^{n,2}$ computed by the scheme \eqref{2d:scheme:full2} with the BP limiter,
satisfies $ m \leq  u^{n,2}  \leq M $.
\end{lemma}
\begin{proof}
 Combining \eqref{2d:scheme:full1}--\eqref{2d:scheme:full2} and the definition of $W^{n}$ in \eqref{def:W}, we have
\begin{align*}
 \bar u_{i,j}^{n,2} 
 &= \left(\malA_{x} \malA_{y} +\frac{\nu \tau}{4} \malA_{y} \mathcal \delta_x^2 + \frac{\nu \tau}{2} \malA_{x}  \delta_y^2 \right)u^{n}_{i,j}
 - \frac{\nu \tau }{4} \malH_{1x}  \delta_y^2 u^{n}_{i,j} \\
 &= \left(\frac{1}{144} 
 \left(\begin{array}{ccc}
   1 & 10 & 1 \\
  10 & 100 & 10 \\
  1 & 10 & 1
 \end{array} \right)
 +\frac{\nu \mu_{1} }{48} \left(\begin{array}{ccc}
  1 & -2 & 1 \\
  10 & -20 & 10 \\
  1 & -2 & 1
 \end{array} \right)
 +\frac{\nu \mu_{2} }{24} \left(\begin{array}{ccc}
  1 & 10 & 1 \\
  -2 & -20 & -2 \\
  1 & 10 & 1
 \end{array} \right) \right) : W^n \\
 & \quad -\frac{\nu \mu_{2} }{48} 
 \left(\begin{array}{ccc}
  1- {3\nu \mu_{1} } &  10+ {6\nu \mu_{1} } &  1- 3{\nu \mu_{1} } \\
  -2+ 6{\nu \mu_{1} } &  -20-12\nu \mu_{1} &  -2+ 6{\nu \mu_{1} } \\
  1- {3\nu \mu_{1} } &  10+ {6\nu \mu_{1} } &  1- 3{\nu \mu_{1} }
 \end{array} \right) : W^{n},
 \end{align*}
 where : denotes the sum of all entrywise products in two matrices of the same size and $\bar u_{i,j}^{n,2} := \malH_{1x} \malH_{1y} u_{i,j}^{n,2}$. Obviously, the right-hand side above is a monotonically increasing function with respect to $u_{l m}$ for $i-1 \leq l \leq i+1, j-1 \leq m \leq j+1$ if $\nu \mu_1= \nu \mu_2 \leq \frac{5}{3}$ holds. The monotonicity implies the bound-preserving result of $\bar u_{i,j}^{n,2}$.

 Given $\bar u^{n,2}$, we can recover the point values $u^{n,2}$ by first obtaining $\hat u^{n,2} = \mathbf H_{1 x}^{-1} \bar u^{n,2}$, then $u^{n,2} = \mathbf H_{1 y}^{-1} \hat u^{n,2}$. Similar to the one-dimensional case, this process can be similarly divided into the following two situations.
\begin{itemize}
 
 \item Case 1: if $\frac{1}{3} \leq \nu \mu_1, \leq \frac{5}{3}$ holds, following the arguments in the proof of Lemma \ref{lem:1d:diffusion}, we can verify  $\mathbf H_{1 x}$ is $\mathbf{M}$-matrix, which doesn't affect the bound-preserving property by Lemma \ref{lem:M-BP}. The same can be obtained for the $y$ direction. Thus, we can derive that $u^{n,2}$ computed by the original scheme \eqref{2d:scheme:full1}--\eqref{2d:scheme:full2} is bound-preserving under the condition that $\frac{1}{3} \leq \nu \mu_1=\nu \mu_2\leq \frac{5}{3}$.

 \item Case 2: And if $\nu \mu_1 < \frac{1}{3}$, although $\mathbf H_{1x}$ isn't $\mathbf{M}$-matrix, it's still a convex combination of $\bar u^{n,2}$, which doesn't change the upper and lower bound with the help of BP limiter.
 Similar to one-dimensional case, we can use the limiter in Algorithm \ref{alg:1d:xie}  dimension by dimension several times to enforce ${u}_{i,j}^{n,2} \in \left[m, M\right]$ : 

(a) Given $\bar{u}_{i,j}^{n,2}$, first compute $\hat{u}_{i,j}^{n,2} =\malA_{x}^{-1} \bar{u}_{i,j}^{n,2}$ which are not necessarily in the range $[m, M]$. 
Notice that 
$$
\bar u_{i,j}^{n,2} =\frac{1}{c+2} \left(\hat{u}_{i-1, j}^{n,3} +c \hat{u}_{i,j}^{n,2}+\hat{u}_{i+1, j}^{n,2} \right), \quad c=10,
$$
thus all discussions in Section \ref{sec:1d} are still valid. 
Then apply the limiter in Algorithm \ref{alg:1d:xie}  to $\hat{u}_{i,j}^{n,2}$ for each fixed $j$ and
denote the output of the limiter as $\hat{v}_{i,j}^{n,2} (i=1,\cdots,N_x) $, thus we have $\hat{v}_{i,j}^{n,2} \in[m, M]$. 

(b) Compute ${u}_{i,j}^{n,2} = \malA_{y}^{-1} \hat {v}_{i,j}^{n,2}$. Then we have
$$\hat {v}_{i,j}^{n,2} =\frac{1}{c+2} \left({u}_{i, j-1}^{n,2} +c {u}_{i,j}^{n,2} +{u}_{i, j+1}^{n,2} \right),\quad c=10. $$
 Apply the limiter in Algorithm \ref{alg:1d:xie} to ${u}_{i,j}^{n,2} (j=1,\cdots,N_y) $ for each fixed $ i $. Then the output values are in the range $[m, M]$.

 \end{itemize}
 
\end{proof}

By the symmetry of \eqref{2d:scheme:full1}--\eqref{2d:scheme:full2} and \eqref{2d:scheme:full5}--\eqref{2d:scheme:full6}, we have similar conclusions for \eqref{2d:scheme:full5}--\eqref{2d:scheme:full6}. Then we have the corollary as follows.

\begin{corollary}
%Under the constraints
Under the constraint $\nu \mu_1 \leq \frac{5}{3}, \nu \mu_2 \leq \frac{5}{3}$, if $u^{n,4}\in[m, M]$, then $u^{n+1}$ computed by the scheme \eqref{2d:scheme:full6} with the BP limiter,
satisfies $ m \leq u^{n+1} \leq M $.
\end{corollary}

Next, we consider the scheme \eqref{2d:scheme:full3}--\eqref{2d:scheme:full4} and give the following results similar to the one-dimensional case.
%%%---------------------------------------------- 
 \begin{lemma} \label{lem:2d:convection:1}
 Under the constraints $\lambda_1 \max \limits_u\left|f^{\prime}(u)\right|+ \lambda_2 \max \limits_u\left|g^{\prime}(u)\right| \leq \frac{1}{3}$, if $u^{n,2} \in[m, M]$, then $ u^{n,3}$ computed by the scheme \eqref{2d:scheme:full3} with the BP limiter satisfies $ m \leq {u}^{n,3}\leq M $.
 \end{lemma}
 \begin{proof}
 Applying the definition of $W^{n,2}$, \eqref{2d:scheme:full3} can be written as
  \begin{align*}
  \bar{u}_{i,j}^{n,3}
   & =  \malB_{x} \malB_{y} u_{i,j}^{n,2} -\tau \malB_{y} \mathcal
     D_{\hat{x}} f_{i,j}^{n,2}-\tau \malB_{x} D_{\hat{y}} g_{i,j}^{n,2} \\
   &=
    \frac{1}{36} \left(\begin{array}{ccc}
    1 & 4& 1 \\
    4 & 16 & 4\\
    1 & 4 & 1 
    \end{array} \right) : W^{n,2}
   + \frac{\lambda}{12} \left(\begin{array}{ccc}
    1 & 0 & -1 \\
    4 & 0 & -4 \\
    1 & 0 & -1
    \end{array} \right)f : W^{n,2}
   + \frac{\lambda_2}{12} \left(\begin{array}{ccc}
    -1 & -4 & -1 \\
    0 & 0 & 0 \\
    1 & 4 & 1
    \end{array} \right) g : W^{n,2},
 \end{align*}
 where $\bar{u}_{i,j}^{n,3} :=\malB_{x} \malB_{y}u_{i,j}^{n,3}$.
Thus, if $\lambda_1 \max \limits_u\left|f^{\prime}(u)\right|+ \lambda_2 \max \limits_u\left|g^{\prime}(u)\right| \leq \frac{1}{3}$ holds, similar to the discussions in Lemma \ref{lem:2d:diffusion}, the weak monotonicity of \eqref{2d:scheme:full3} can be guaranteed.

 Similarly, given $\bar{u}_{i,j}^{n,3}$, we can recover the point values $u_{i,j}^{n,3}$ by first obtaining $\hat{u}_{i,j}^{n,3} = \malB_{x}^{-1}\bar{u}_{i,j}^{n,3}$, then $u_{i,j}^{n,3} = \malB_{y}^{-1} \hat{u}_{i,j}^{n,3}$. Following the notations and arguments in the proof of Lemma \ref{lem:2d:diffusion}, we can ensure that ${u}_{i,j}^{n,3} \in\left[m, M\right]$ by using the limiter in Algorithm \ref{alg:1d:xie} dimension by dimension several times with $c=4$. 
 
\end{proof}

%--------------------------------------
\begin{lemma} \label{lem:2d:convection:2}
 Under the CFL constraints $\lambda_1 \max \limits_{u}\left|f^{\prime}(u)\right|+ \lambda_2 \max \limits_u\left|g^{\prime}(u)\right| \leq \frac{1}{3}$, if $\left\{u^{n,2}, u^{n,3} \right\} \in[m, M]$, then $u^{n,4}$ computed by the scheme \eqref{2d:scheme:full4}
 with the BP limiter satisfies $ m \leq u^{n,4} \leq M $.
 
\end{lemma}  
\begin{proof} 
 Combining \eqref{2d:scheme:full3} and \eqref{2d:scheme:full4}, we get
\begin{align*}
 \begin{aligned}
  \malB_{x} \malB_{y} u_{i,j}^{n,4}
  =& \malB_{x} \malB_{y} u_{i,j}^{n,2} -\frac{\tau}{2} \malB_{y} \mathcal
 D_{\hat{x}} \big(f_{i,j}^{n,2}+ f_{i,j}^{n,3}\big)-\frac{\tau}{2} \malB_{x} D_{\hat{y}} \left(g_{i,j}^{n,2} + g_{i,j}^{n,3}\right)
  \\
  =& \frac {1}{2} \big(\malB_{x} \malB_{y} u_{i,j}^{n,2} -{\tau} \malB_{y} D_{\hat{x}} f_{i,j}^{n,2} -{\tau} \malB_{x} D_{\hat{y}} \underline{} g_{i,j}^{n,2} \big) 
  \\
  \quad & +\frac {1}{2} \big(\malB_{x} \malB_{y} u_{i,j}^{n,2} -{\tau} \malB_{y} D_{\hat{x}} f_{i,j}^{n,3}-{\tau} \malB_{x} D_{\hat{y}} g_{i,j}^{n,3} \big)\\
  =& \frac {1}{2} \malB_{x} \malB_{y} u_{i,j}^{n,2} + \frac {1}{2} \big(\malB_{x} \malB_{y}u_{i,j}^{n,3} -{\tau} \malB_{y} D_{\hat{x}} f_{i,j}^{n,3} -{\tau} \malB_{x} D_{\hat{y}} g_{i,j}^{n,3} \big).
 \end{aligned}
\end{align*}
 Following the notations and arguments in the proof of Lemma \ref{lem:2d:convection:1}, we can show
  $  {u}_{i,j}^{n,4} \in [m,M]. $

\end{proof}
% Similarly, given ${z}_{i,j}^{n,4}$, we can recover point values $u_{i,j}^{n,4}$ by obtaining first $\bar{u}_{i,j}^{n,4} = \malB_{x}^{-1}{z}_{i,j}^{n,4}$, then $u_{i,j}^{n,4} = \malB_{y}^{-1} \bar{u}_{i,j}^{n,4}$. Thus similar to the discussions in the previous subsection, given point values $z_{i,j}^{n,4}$ satisfying $z_{i,j}^{n,4} \in[m, M]$ for any $i$ and $j$, we can use the limiter in Algorithm \ref{alg:1d:xie}  dimension by dimension several times to enforce ${u}_{i,j}^{n,4} \in\left[m, M\right]$ :\\
% 1. Given ${z}_{i,j}^{n,4}$, compute $\bar{u}_{i,j}^{n,4} = \malB_{x}^{-1}{z}_{i,j}^{n,4}$ which are not necessarily in the range $[m, M]$. Then apply the limiter in Algorithm \ref{alg:1d:xie}  to $\bar{u}_{i,j}^{n,4}$ for each fixed $j$. Since we have
% $$
% {z}_{i,j}^{n,4} =\frac{1}{c+2} \left(\bar{u}_{i-1, j}^{n,4} +c \bar{u}_{i,j}^{n,4} +\bar{u}_{i+1, j}^{n,4} \right), c=4,
% $$
% all discussions in section \ref{alg:1d:xie}  are still valid. Let $\bar{v}_{i,j}^{n,4}$ denote the output of the limiter; thus we have $\bar{v}_{i,j}^{n,4} \in[m, M]$. \\
% 2. Compute ${u}_{i,j}^{n,4} = \malB_{y}^{-1} \bar{v}_{i,j}^{n,4}$. Then we have
% $$\bar {v}_{i,j}^{n,4} =\frac{1}{c+2} \left({u}_{i, j-1}^{n,4} +c {u}_{i,j}^{n,4} +{u}_{i, j+1}^{n,4} \right), c=4. $$
% Apply the limiter in Algorithm \ref{alg:1d:xie}  to ${u}_{i,j}^{n,4}$ for each fixed $ i $. Then the output values are in the range $[m, M]$.  

 Combining with the above results, we have the following theorem without proof.
\begin{theorem} \label{th:2d:all}
 Under the constraint
 \begin{align}
  & \nu\mu_1= \nu \mu_2 \leq \frac{5}{3},  \\
  &\lambda_1 \max \limits_u  |f^{\prime}(u)| +\lambda_2  \max \limits_u  |g^{\prime}(u)| \leq \frac{1}{3}, 
\end{align}
 if $u^n \in[m, M]$, then $u^{n+1}$ computed by the scheme \eqref{2d:scheme:full1}--\eqref{2d:scheme:full6} with the BP limiter
 satisfies $ m \leq u^{n+1} \leq M $.
 % , after applying the limiter in Algorithm 1 to $\malA_{y} u^{n,5}_{i,j}, \malB_{x} \malB_{y} u^{n,3}_{i,j}$ and $\malB_{x} \malB_{y} u^{n,4}_{i,j}$
 % for all $ (i=0, \ldots, N_x, j=0, \ldots, N_y)$.
 % In particular, if the CFL constraints $\lambda_1 \max \limits_{u}\left|f^{\prime}(u)\right|+ \lambda_2 \max \limits_u\left|g^{\prime}(u)\right| \leq \frac{1}{3}$ satisfied, the scheme is stable.
\end{theorem}

Similarly, we give Algorithm \ref{alg:3} for the two-dimensional problem.
\begin{algorithm} [ht]
 \caption{Strang-ADI-HOC algorithm with the BP limiter}
 \label{alg:3}
 \begin{algorithmic} [1]
  \FOR{$n=0$ to $N_t-1$}
  \STATE \textbf{Step1:} Compute $\malA_{y} u^{n,1}$ along $x$-direction by \eqref{2d:scheme:full1} with initial value $u^n$ on $ (t_n,t_{n+1/2})$. 
  
   \IF {$\{\nu \mu_1,\nu \mu_2\} < \frac{1}{3}$}
  \STATE Applying the limiter in Algorithm \ref{alg:1d:xie} to $\malH_{1x} \malA_{y} u^{n,1}$ to obtain $\malA_{y} u^{n,1}$ which is bound-preserving.
  \ENDIF
  
  \STATE \textbf{Step2:} Compute $\malH_{1y} u^{n,2}$ along $y$-direction by \eqref{2d:scheme:full2} with initial value $u^n$ and $\malA_{y}u^{n,1}$ on $(t_n,t_{n+ 1/2})$.

     \IF {$\{\nu \mu_1,\nu \mu_2\} < \frac{1}{3}$}
  \STATE Applying the limiter in Algorithm \ref{alg:1d:xie} to $\malH_{1y} u^{n,2}$ to obtain $u^{n,2}$ which is bound-preserving.
  \ENDIF
 
  \STATE \textbf{Step3:} Compute $\malB_{x} \malB_{y} u^{n,3}$ by \eqref{2d:scheme:full3} with initial value $u^{n,2}$ on $(t_n,t_{n+1})$.
 
  \STATE \textbf{Step4:} Applying the limiter in Algorithm \ref{alg:1d:xie} to $\malB_{x} \malB_{y} u^{n,3}$ dimension by dimension to obtain $u^{n,3}$ which is bound-preserving.
 
  \STATE \textbf{Step4:} Compute $\malB_x \malB_x u^{n,4}$ by \eqref{2d:scheme:full4} with initial value $u^{n,2}$ and $u^{n,3}$ on $(t_n,t_{n+1})$.
 
  \STATE \textbf{Step5:} Applying the limiter in Algorithm \ref{alg:1d:xie} to $\malB_x \malB_x u^{n,4}$ dimension by dimension to obtain $u^{n,4}$ which is bound-preserving.
  
  \STATE \textbf{Step6:} Compute $\malA_{y} u^{n,5}$ by \eqref{2d:scheme:full5} with initial value $u^{n,4}$ on $(t_{n+ \frac{1}{2}},t_{n+1})$.

    \IF {$\{\nu \mu_1,\nu \mu_2\} < \frac{1}{3}$}
  \STATE Applying the limiter in Algorithm \ref{alg:1d:xie} to $\malH_{1x} \malA_{y} u^{n,1}$ to obtain $\malA_{y} u^{n,1}$ which is bound-preserving.
  \ENDIF
  
  \STATE \textbf{Step7:} Compute $\malH_{1y} u^{n+1}$ by \eqref{2d:scheme:full6} with initial value $u^{n,4}$ and $\malA_{y} u^{n,5}$ to obtain  on $(t_{n+ \frac{1}{2}},t_{n+1})$.

   \IF {$\{\nu \mu_1,\nu \mu_2\} < \frac{1}{3}$}
  \STATE Applying the limiter in Algorithm \ref{alg:1d:xie} to $\malH_{1y} u^{n,1}$ to obtain $ u^{n,1}$ which is bound-preserving.
  \ENDIF
  
  \ENDFOR
 \end{algorithmic}
\end{algorithm}

\begin{remark} \label{rem:CFL:h_tau:2D}
    Similarly to Remark \ref{rem:CFL:h_tau:2D}, We can discuss the selection of step size based on the value of $\nu$. For simplicity, here we assume $\nu_1=\nu_2=\nu$ and $\lambda_1\max \limits_{u}\left|f^{\prime}(u)\right|=\lambda_2\max \limits_{u}\left|g^{\prime}(u)\right|$. When $\nu$ is large enough  so that  $h < \frac{\nu}{10\max \limits_{u}\left|f^{\prime}(u)\right|}$,   we only need $\tau < \frac{\nu}{60\max \limits_{u}\left|f^{\prime}(u)\right|}$. 
    On the other hand, when $\nu$ is large enough so that $h>\frac{\nu}{10\max \limits_{u}\left|f^{\prime}(u)\right|}$,   we only need $\tau < \frac{h}{6\max \limits_{u}\left|f^{\prime}(u)\right|}$.
\end{remark}

 \section{Numerical tests} \label{sec:tests}
 In this section, we perform several numerical examples
 including linear and nonlinear problems with different dimension to test accuracy, efficiency, and effectiveness in preserving bound or/and mass of the proposed schemes.
% All the numerical experiments are performed using Matlab R2022b on a desktop with the configuration:
% 12th Gen Intel(R) Core(TM) i5-12400F @ 2.50GHz and 16 GB RAM.
 Under the periodic boundary condition, the error of mass is measured by $$Mass_{err}:= h_x \sum_{i=0}^{N_x} u_{i}^n - h_x \sum_{i=0}^{N_x} u_{i}^o$$ for one dimension and $$Mass_{err}:= h_x h_y \sum_{i=0}^{N_x} \sum_{j=0}^{N_y} u_{i,j}^n - h_x h_y \sum_{i=0}^{N_x} \sum_{j=0}^{N_y} u_{i,j}^o$$ for two dimension. Similarly, the mass error under the Dirichlet boundary condition can be defined similarly
In the following tests, we always assume that $N = N_x = N_y$ and thus $h = h_x = h_y$. 
 And in this section, assume $u$ represents the exact solution and $U$ represents the numerical solution. The upper and lower bound errors are respectively defined as
 $ M_{err}:= M-\mathop{\max} \{U\}$ and $m_{err}:= \mathop{\min} \{U\}-m$. Based on the discussion of Remark \ref{rem:CFL:h_tau:1D} and \ref{rem:CFL:h_tau:2D},  $\tau$ and $h$ satisfy the proposed bound-preserving condition in Theorem \ref{th:1d:BP:all} because of the very small $\nu$.

 %%%%%%%%%%%%%%%%%%%%%%%%%%%%%%%%%%%%%%%%%%%%%%%%%%%%
 \subsection{ One-dimensional problems} \label{sec:test:1d}
 %%%%%%%%%%%%%%%%%%%%%%%%%%%%%%%%%%%%%%%%%%%%%%%%%%%%
% In this part, we present the numerical results on several 1D problems. We tested accuracy and bound-preserving property of \eqref{2d:scheme:full1}-\eqref{scheme:1d:full:4}.
%
% The order of our method will be measured using discrete $L_{2}$ and $L_\infty$ error respectively, which errors and convergence rate are defined as follows:
% \begin{align*}
% &  \Vert \text{Error} \Vert_{\infty} =\max_{0 \leq i \leq N} \left|U_{i,m}-u\left(x_i, t_M\right)\right|, \quad  \Vert \text{Error} \Vert_2=\sqrt{h \sum_{i=0}^N\left[U_{i,m}-u\left(x_i, t_M\right)\right]^2},\\
%  & \text {Rate} =\frac{\log \left[ \Vert \text{Error} \Vert_{i} \left(h_1\right) / \Vert \text{Error} \Vert_{i} \left(h_2\right)\right]}{\log \left(h_1 / h_2\right)}, i =2 \text{or} \infty.
% \end{align*}

 \paragraph{Example 1} \label{eg1} Firstly, we consider the following linear convection diffusion problem \cite{acosta2010mollification}
 $$
 u_t + u_x = \nu u_{xx}, \quad 0 \leqslant x \leqslant 5, 0 \leqslant t \leqslant 2
 $$
 with exact solution
 $$
 u(x, t)= \frac{1}{\sqrt{1+t}} \exp (-\frac{{(x- (1+t))^2}}{{4\nu (1+t)} }),
 $$
 where the initial and boundary conditions can be derived from this exact solution.

%% =========================================
 \begin{figure} [!htbp]
 \centering
 \subfigure[{Without the BP limiter} ]
 {
 %\label{Fig. sub. 1}
 \includegraphics[width=0.45\textwidth]{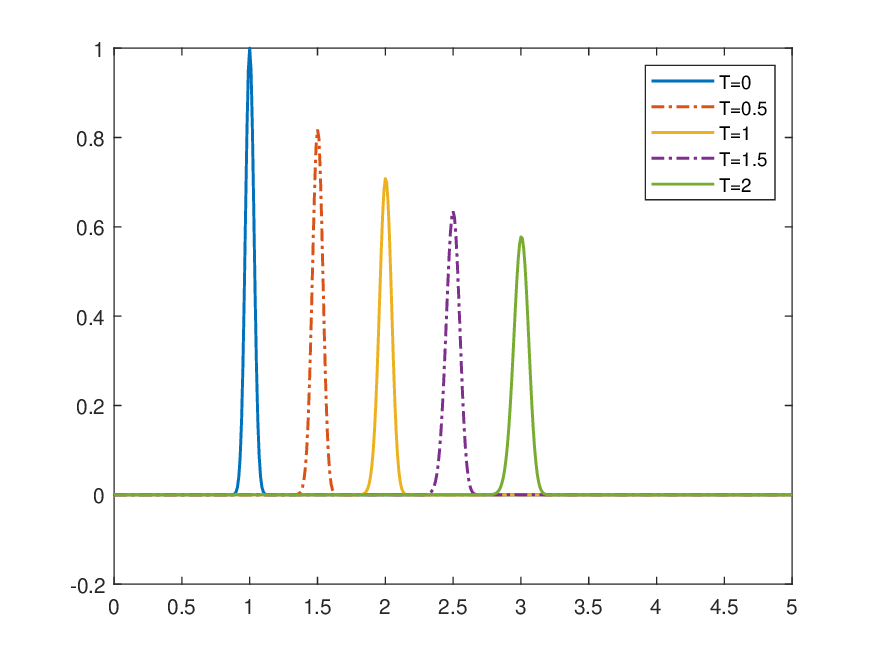}
}
 \subfigure[{With the BP limiter} ]
 {
% \label{Fig. sub. 2}
 \includegraphics[width=0.45\textwidth]{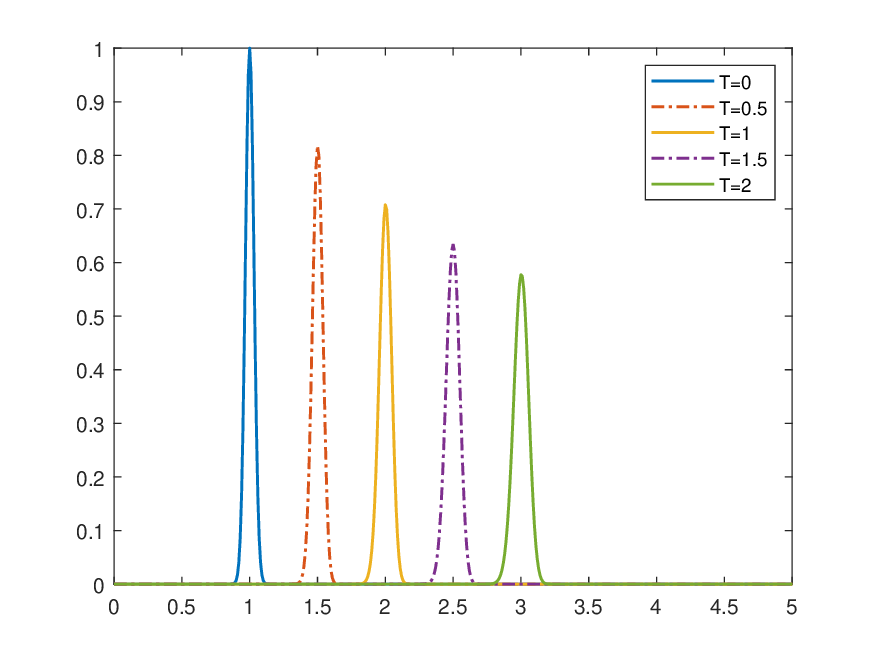}
}
 \subfigure[{Without the BP limiter} ]
 {
% \label{Fig. sub. 1}
 \includegraphics[width=0.45\textwidth]{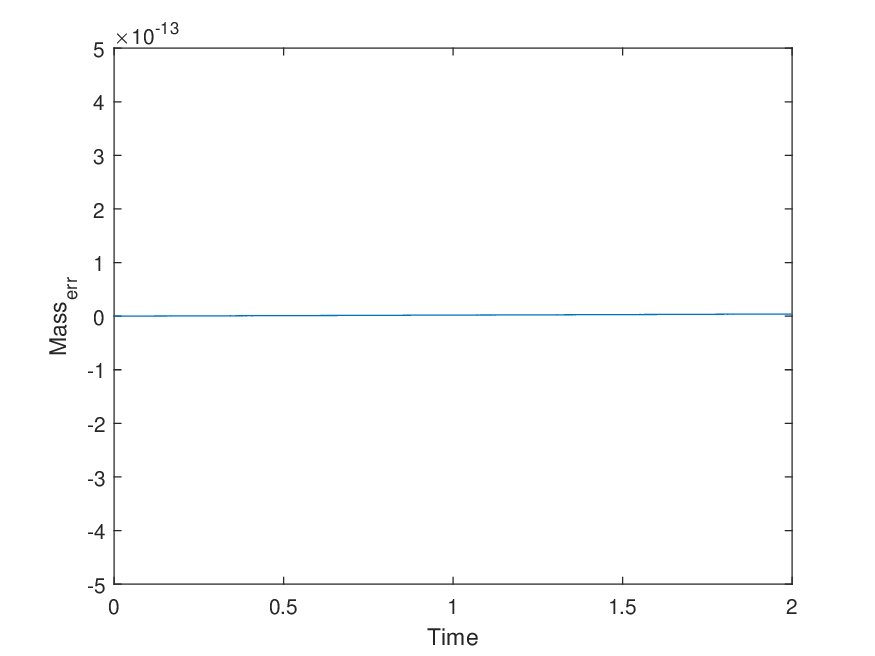}
}
 \subfigure[{With the BP limiter} ]
 {
 %\label{Fig. sub. 2}
 \includegraphics[width=0.45\textwidth]{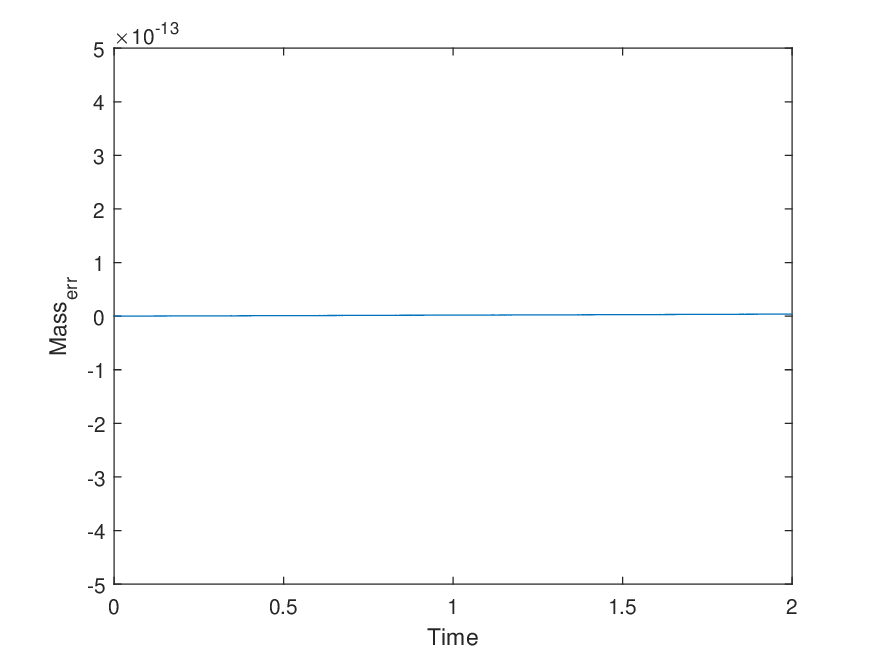}
}
 \caption{Numerical solutions and the evoution of  $Mass_{err}$ without/with the BP limiter at different time.}
 \label{fig:eg1:diff_u:mass_err}
 \end{figure}
%%%%% ============================================

 In Fig. \ref{fig:eg1:diff_u:mass_err}, we show the numerical solutions and the evoution of  $Mass_{err}$ without/with the BP limiter for $\nu =5 \times 10^{-4}$ at different time $T=0,0. 5,1$ and 2. We set $N= 500, \tau= \frac{h}{3 \max \limits_u\left|f^{\prime} \right|}$, and $\max \limits_u\left|f^{\prime} \right|=1.$  It clearly shows the Strang-HOC difference scheme \eqref{scheme:1d:full:1}--\eqref{scheme:1d:full:4} can only ensure mass conservation but not bound-preserving property, while the scheme with BP limiter can ensure both properties at the same time, which clearly illustrates the role of the limiter.
 Furthermore, numerical errors and convergence rates are presented in {Table \ref{tab:eg1:order}} with $\tau = h^2, \nu = 5\times 10^{-4}$ to test both time and space accuracy simultaneously.
 It shows that the schemes without/with BP limiter both achieve fourth order in space and second order in time.
 The corresponding upper  bound error $M_{err}$ and lower bound error  $m_{err}$ are also shown in Table \ref{tab:eg1:bound}.
Based on the definition of $M_{err}$($m_{err}$), if a negative value appears, it means that the numerical solution exceeds the upper (lower) bound.
 It can be found that although the upper bound of the two schemes are both maintained, the lower bound can only be maintained the scheme using the limiter.
 Thus, combining Table \ref{tab:eg1:order} and Table
 \ref{tab:eg1:bound}, it clearly shows that the use of BP limiter works well and does not affect numerical accuracy, which coincide with the theoretical analysis in Remark \ref{rem:order}.
% In addition, we find that the constraints we have prosed only the sufficient condition. In fact, in many cases, we only need to guarantee the CFL condition from the convection terms, which gives a more appropriate spatial step size compared with the constraints in Theorem \ref{th:1d:BP:all}. This also confirms the conclusion in Remark \ref{rem:sufficient}.

%%%%% ============================================
 \begin{table} [!htbp]
 \vspace{-0.1em}
 \centering
 \caption{
 Accuracy test with $\nu =5 \times 10^{-4}, \tau= {h}^2 $.
}
 \label{tab:eg1:order}
 \setlength{\tabcolsep}{1.2mm}{
 \begin{tabular}{|c|cccc|cccc|}
 \hline & \multicolumn{4}{c|}{\text {The scheme without limiter} }
  & \multicolumn{4}{c|}{\text {With limiter } } \\
 \hline
 $N$ & $L^{\infty}$ error & Order & $ L^{2}$ error & Order
   & $L^{\infty}$error & Order & $ L^{2}$ error & Order \\
 \hline
 300 & 6. 1928 E -3 & ---
   & 1. 7708 E -3 & ---
   & 6. 1924 E -3 & ---
   & 1. 7704 E -3 & --- \\
 \hline
 400 & 1. 8850 E -3 & 4. 13
   & 5. 4629 E -4 & 4. 09
   & 1. 8850 E -3 & 4. 13
   & 5. 4629 E -4 & 4. 09 \\
 \hline
 500 & 7. 5069 E -4 & 4. 13
   & 2. 2092 E -4 & 4. 06
   & 7. 5069 E -4 & 4. 13
   & 2. 2092 E -4 & 4. 06 \\
 \hline
 600 & 3. 6287 E -4 & 3. 99
   & 1. 0580 E -4 & 4. 04
   & 3. 6287 E -4 & 3. 99
   & 1. 0580 E -4 & 4. 04 \\
 \hline
 700 & 1. 9357 E -4 & 4. 08
   & 5. 6870 E -5 & 4. 03
   & 1. 9357 E -4 & 4. 08
   & 5. 6870 E -5 & 4. 03 \\
 \hline
 \end{tabular}
}
 \vspace{-0.1em}
 \end{table}
%%%%% =========================================

%%%%% =========================================
 \begin{table} [!htbp]
 \vspace{-0.8em}
  \centering
  \caption{Verification of the bound-preserving property with $\nu =5 \times 10^{-4}, \tau= {h}^2 $.
  }
  \label{tab:eg1:bound}
 \setlength{\tabcolsep}{1.2 mm}{
  \begin{tabular}{|c|cc|cc|}
   \hline & \multicolumn{2}{c|}{\text {Without the BP limiter} } &
      \multicolumn{2}{c|}{\text {With the BP  limiter } } \\
   \hline
    $N$
    & $ M_{err}$
    & $ m_{err}$
    & $ M_{err}$
    & $ m_{err}$
 \\
   \hline 300 & 4. 2343 E -1 & -2.1743 E -5
       & 4. 2343 E -1 & 0 \\
   \hline 400 & 4. 2271 E -1 & -1.1938 E -7 &
        4. 2271 E -1 & 0 \\
   \hline 500 & 4. 2265 E -1 & -8.0145 E -11 &
       4. 2265 E -1 & 0 \\
   \hline 600 & 4. 2265 E -1 & -5.8332 E -15 &
       4. 2265 E -1 & 0 \\
    \hline 700 & 4. 2265 E -1 & -3.1338 E -19 &
       4. 2265 E -1 & 0 \\
   \hline
  \end{tabular}
}
 \vspace{-0.8em}
 \end{table}
%%%%% =========================================
 \paragraph{Example 2}
 Considering the viscous Burgers equation
 \begin{align*}
  u_t + f(u)_x = \nu u_{xx}
 \end{align*}
 where $f(u) = \frac{1}{2} u^2 $
 and give the following three forms of solutions or initial condition: 
 \begin{itemize}
  \item Case {1}\cite{ali1992}:  $u(x, t)=\frac{x / t}{1 + \sqrt{t} \exp \left((x^{2} -0. 25t) / 4 \nu t\right)}, \quad (x, t) \in [0,3]\times [1,4]$.

  \item Case {2} \cite{acosta2010mollification}: $u(x,t) = 1.5-0.5 \tanh\left(0.5(x-1.5t)/2\nu\right), \quad (x, t) \in [-1,3] \times [0,0.5]$.

  \item Case {3} \cite{xie2008numerical}: $u(x,0)=\sin(2\pi x), \quad (x, t) \in [0,1] \times [0,2]$.
 \end{itemize}
 The initial and Dirichlet boundary conditions are derived from the exact solution in Case  {1} and {2}. For Case {3}, periodic boundary condition are applied.

The presence of large gradient in the solution to the Burgers equation is a well-known challenge of numerical simulation. Even if the initial data are sufficient smooth, large gradient will still occur when the characteristic curves of the Burgers equation intersect. 
Thus, a robust and accurate numerical algorithm should be used to capture large gradient, and the numerical solution is anticipated to exhibit the correct physical behavior. 
% In addition, from the expression of exact solution in Case {1} and {2}, we see that the evolution of solution always keep non-negative over the spatial domain, which requires that the numerical schemeshould also preserve this property. 

\begin{figure} [!htbp]
  \centering
  \subfigure[ Without any limiter ]
  {
   \includegraphics[width=0.45\textwidth]{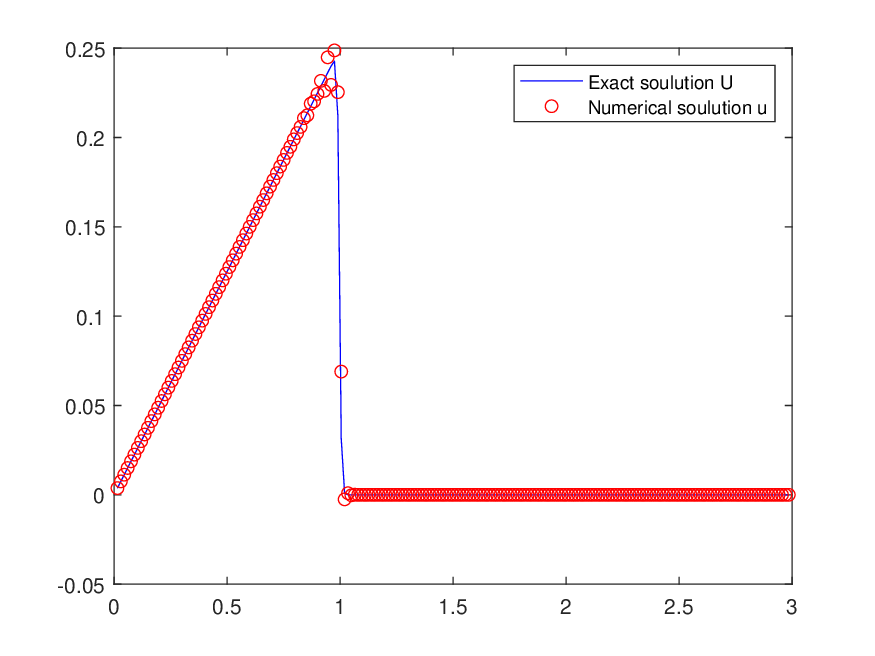}
  }
  \subfigure[ Only with TVB limiter ]
  {
   \includegraphics[width=0.45\textwidth]{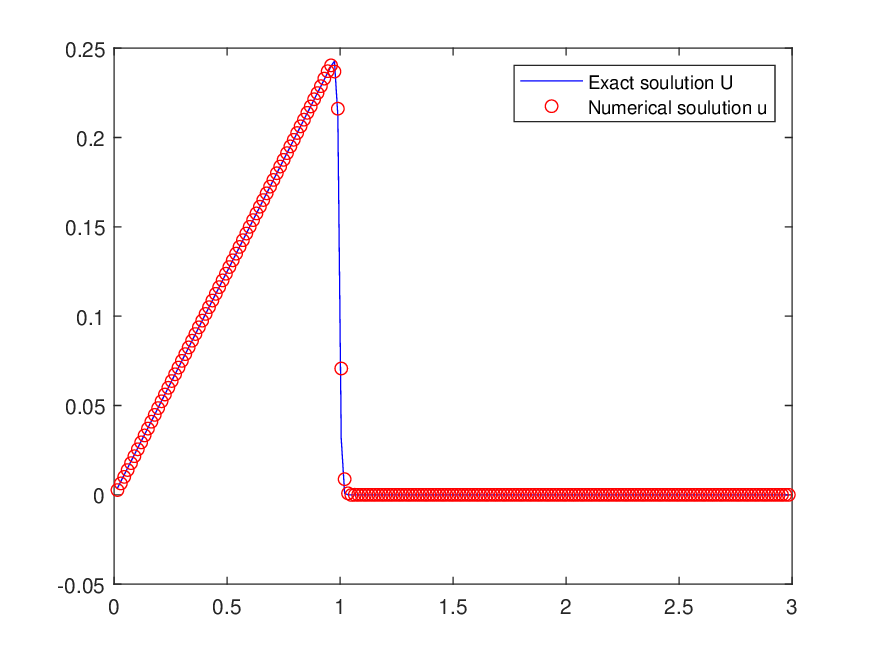}
  }
  \subfigure[ Only with BP limiter ]
  {
   \includegraphics[width=0.45\textwidth]{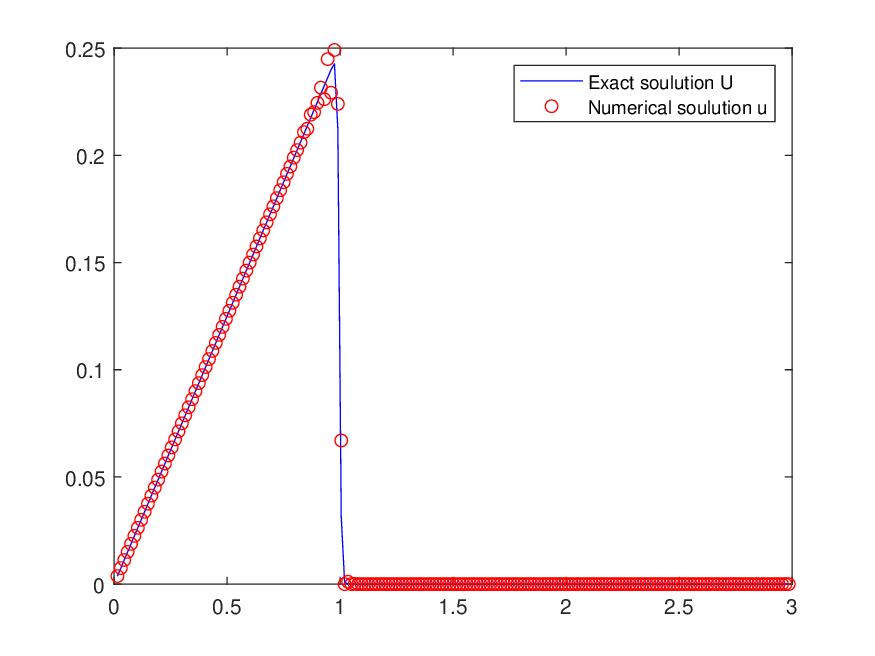}
  }
  \subfigure[ With both limiters ]
  {
   \includegraphics[width=0.45\textwidth]{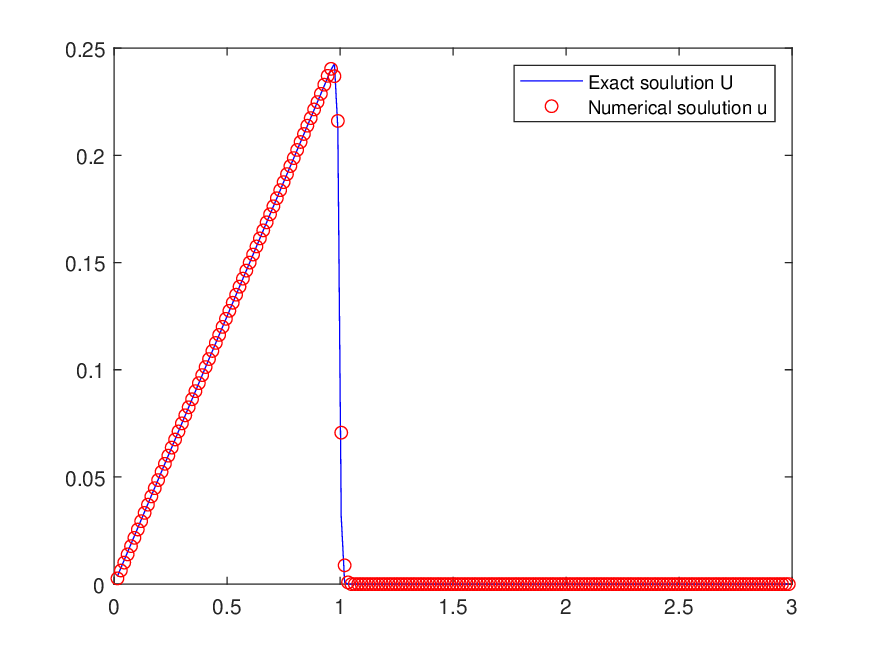}
  }
  \caption{Comparisons of numerical solutions obtained without/with TVB or/and BP limiter. }
  \label{fig:2g2:case1:u}
 \end{figure}

\begin{figure} [!htbp]
  \centering
  \subfigure[Without any limiter ]
  {
   \includegraphics[width=0.45\textwidth]{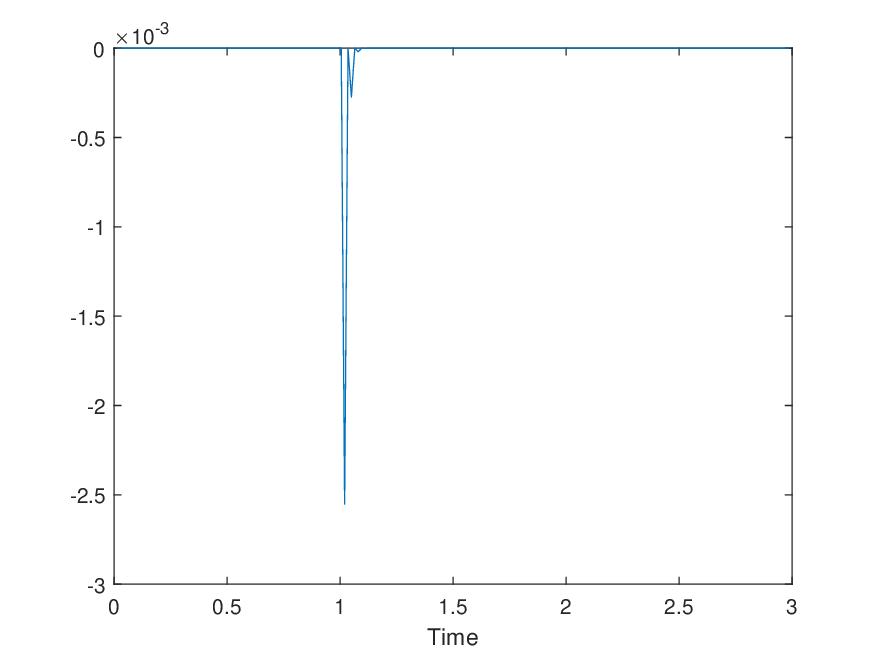}
  }
  \subfigure[ Only with TVB limiter ]
  {
   \includegraphics[width=0.45\textwidth]{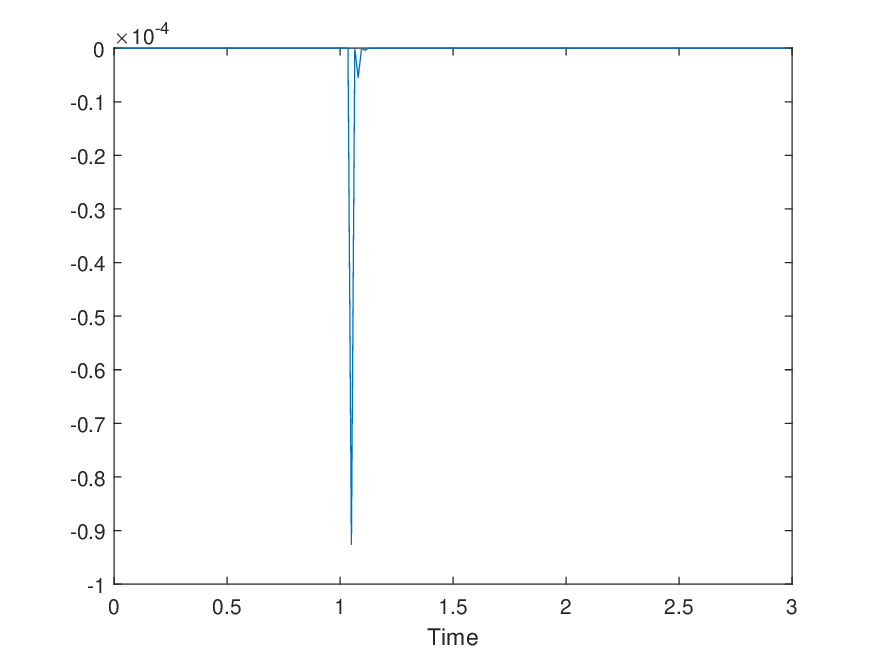}
  }
  \subfigure[ Only with BP limiter ]
  {
   \includegraphics[width=0.45\textwidth]{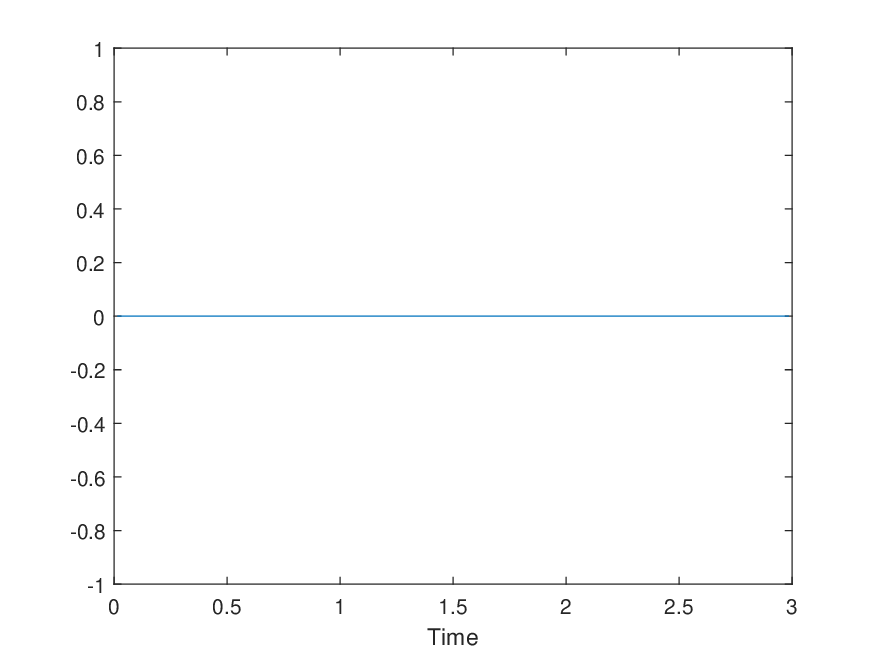}
  }
  \subfigure[ With both limiters ]
  {
   \includegraphics[width=0.45\textwidth]{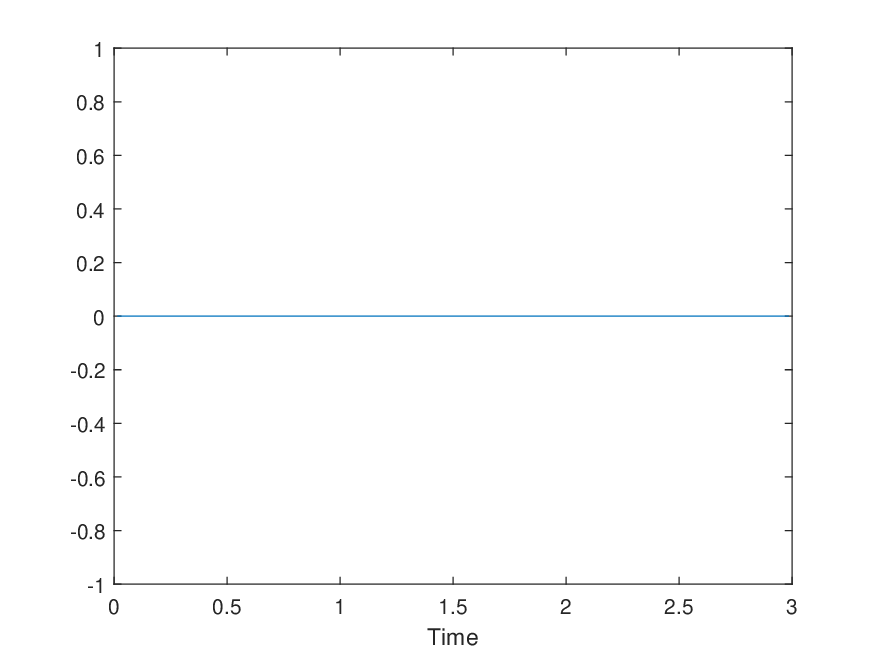}
  }
  \caption{Comparisons of $min(U,0)$ obtained without/with TVB or/and BP limiter. }
  \label{fig:2g2:case1:min}
 \end{figure}

 Firstly, we test the efficiency, bound preservation, and accuracy of the proposed scheme \eqref{scheme:1d:full:1}--\eqref{scheme:1d:full:4} in Case {1}.
 To better visualize the large gradients, we plot the numerical and exact solutions in Fig. \ref{fig:2g2:case1:u} and compare numerical solutions obtained without/with BP limiter or/and TVB limiter at $T =4$ with $\nu =5 \times 10^{-4},N = 200, \tau= \frac{h}{3 \max \left|f^{\prime}(u)\right|}$, where $\max |f^{\prime}| \approx 0.25$. In addition, from the expression of exact solution in Case {1}, we see that the evolution of solution always keep non-negative over the spatial domain, which requires that the numerical scheme should also preserve this property. 
 In Fig. \ref{fig:2g2:case1:min}, we plot the minimum values $min(U,0)$ of the different solutions under the same condition as Fig. \ref{fig:2g2:case1:u} to clearly compare the bound-preserving effects. Here, values less than $10^{-20}$ are considered negligible and treated as 0.
 It is easy to find that the TVB limiter eliminates oscillations but does not remove the overshoot/undershoot. However, when both limiters are used together, a non-oscillatory, bound-preserving numerical solution is achieved.
 In order to test the temporal and spatial accuracy simultaneously, we take $\tau=h^2$ with $\nu=5 \times 10^{-4}$. 
 The numerical results in Table \ref{tab:eg2:order} and \ref{tab:eg2:bound} show that the proposed scheme with the BP limiter is bound-preserving while maintaining its original accuracy. Specifically, the scheme achieves fourth-order accuracy in space and second-order accuracy in time in both the $L_2$ norm and $L_{\infty}$ norm. 

%% ======================================
%% =========================================
% 
% 
%% =========================================
%% ==========================================
% %============================================
 \begin{table} [!htbp]
 \vspace{-0.1em}
 \centering
 \caption{
 Accuracy test with $\nu =5 \times 10^{-4}, \tau= {h}^2 $.
}
 \label{tab:eg2:order}
 \setlength{\tabcolsep}{1.2mm}{
 \begin{tabular}{|c|cccc|cccc|}
 \hline & \multicolumn{4}{c|}{\text {Without the BP  limiter} }
  & \multicolumn{4}{c|}{\text {With the BP  limiter } } \\
 \hline
 $N$ & $L^{\infty}$ error & Order & $ L^{2}$ error & Order
   & $L^{\infty}$error & Order & $ L^{2}$ error & Order \\
 \hline
 400 & 4.3780 E -3 & ---
   & 4.3187 E -4 & ---
   & 4.3782 E -4 & ---
   & 4.3189 E -4 & --- \\
 \hline
 500 & 1.2769 E -3 & 5.52
   & 1.3548 E -4 & 5.20
   & 1.2769 E -3 & 5.52
   & 1.3548 E -4 & 5.20 \\
 \hline
 600 & 5.4157 E -4 & 4.70
   & 5.0094 E -5 & 5.46
   & 5.4157 E -4 & 4. 70
   & 5.0094 E -5 & 5.46 \\
 \hline
 700 & 2.8648 E -4 & 4. 13
   & 2.4117 E -5 & 4.74
   & 2.8648 E -4 & 4. 13
   & 2.4117 E -5 & 4.74 \\
 \hline
 \end{tabular}
}
 \vspace{-0.1em}
 \end{table}

%%%% ------------------------------------------------------------
 \begin{table} [!htbp]
 \vspace{-0.8em}
  \centering
  \caption{Verification of the bound-preserving effect with $\nu = 5 \times 10^{-4}, \tau= {h}^2 $.
  }
  \label{tab:eg2:bound}
 \setlength{\tabcolsep}{1.2 mm}{
  \begin{tabular}{|c|cc|cc|}
   \hline & \multicolumn{2}{c|}{\text {The scheme without the BP limiter} } &
      \multicolumn{2}{c|}{\text {The scheme with the BP limiter } } \\
   \hline
    $N$
    & $ M_{err}$
    & $ m_{err}$
    & $ M_{err}$
    & $ m_{err}$
 \\
   \hline 400 & 2.4415 E -1 & -1.4056 E -8
       & 2.4415 E -1 & 0 \\
   \hline 500 & 2.4291 E -1 & -3.4915 E -15 &
        2.4291 E -1 & 0 \\
   \hline 600 & 2.4383 E -1 & 0 &
       2.4383 E -1 & 0 \\
   \hline 700 & 2.4426 E -1 & 0 &
       2.4426 E -1 & 0 \\
   \hline
  \end{tabular}
}
 \vspace{-0.8em}
 \end{table}
%%%%%% ------------------------------------------------------------

\begin{figure} [!htbp]
  \centering
%  \subfigure[ $N = 200, \nu = 1/64$ without the limiter. ]
%  {
%   \label{fig:1d:step:64_non}
%   \includegraphics[width=0.45\textwidth]{figure/eg2_Case2_64_u_non.eps}
%  }
%  \subfigure[ $N = 200, \nu = 1/64 $ with the BP limiter. ]
%  {
%  \label{fig:1d:step:64}
%   \includegraphics[width=0.45\textwidth]{figure/eg2_Case2_64_u.eps}
%  }
  \subfigure[ Without the  BP limiter. ]
  {
  %\label{fig:1d:step:wan_non}
   \includegraphics[width=0.45\textwidth]{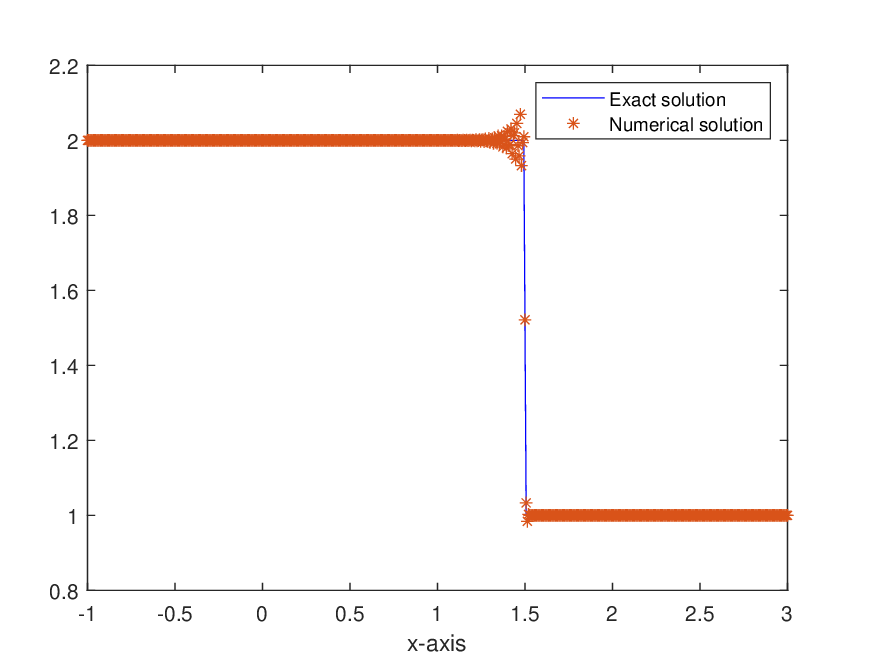}
  }
  \subfigure [With the BP limiter. ]
  {
  %\label{fig:1d:step:wan}
   \includegraphics[width=0.45\textwidth]{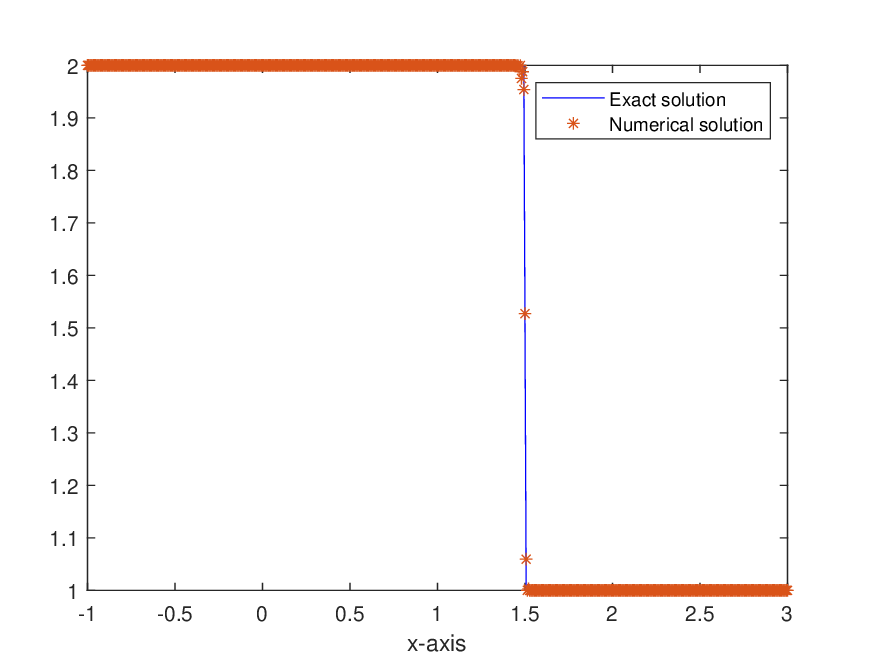}
  }
  \caption{Comparisons of numerical solutions obtained without/with BP limiter under $N =600, \tau= \frac{h}{3 \max \limits_{u}\left|f^{\prime}(u)\right|}$ and $\nu= 5 \times 10^{-4}$. }
  \label{fig:1d:step:nonhomo}
 \end{figure}

 Next, we consider Case {2} with non-homogeneous Dirichlet boundary. From the expression of exact solution, we see that the evolution of solution always remains within $[1,2]$ over the spatial domain. This requires that the numerical scheme should also preserve this property.   Fig. \ref{fig:1d:step:nonhomo} shows numerical solutions without/with the BP limiter at $T =1$ and set $\nu = 5 \times 10^{-4},N=600,\tau= \frac{h}{3 \max \left|f^{\prime}(u)\right|}$, where $\max \left|f^{\prime}(u)\right|=2$. 
 With the help of BP limiter, the numerical solution agrees well with the exact solution and preserves the bound. Indeed, excellent wave characteristics can be clearly observed with the applied approaches. In contrast, the numerical solution without the limiter produced significant oscillations at large gradients and failed to preserve the bound.

 %%%%%% ==========================================
%%%%%% ==========================================

\begin{figure} [!htbp]
  \centering
  \includegraphics[width=6cm,height=4.5cm]{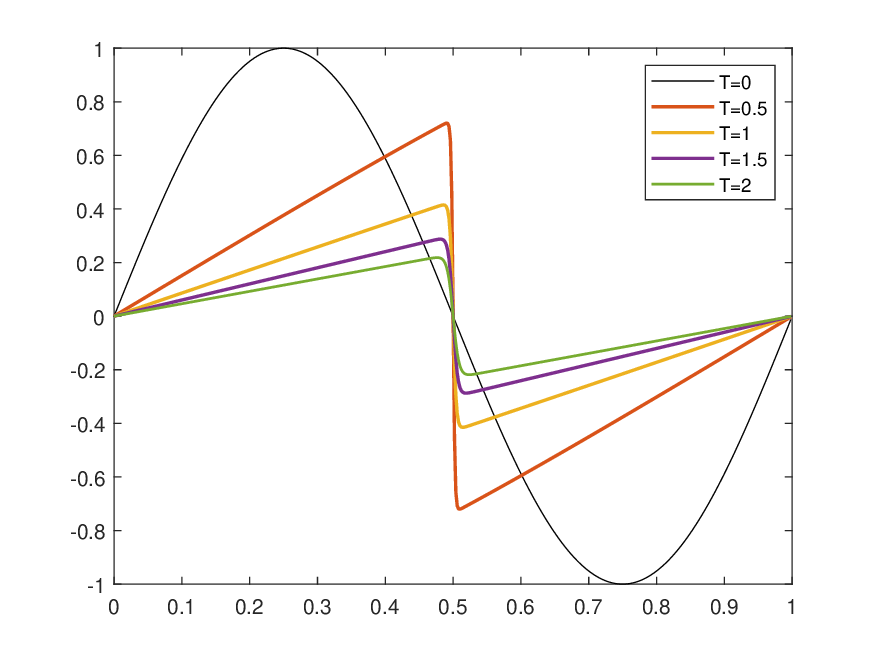}
  \caption{Numerical solutions at different time with $\nu = 10^{-3 },N =500,\tau=\frac{h}{2 \max{|f'|} }$. }
  \label{fig:1d:sin:u}
 \end{figure}

\begin{figure} [!htbp]
  \centering
  \includegraphics[width=6cm,height=4.5cm]{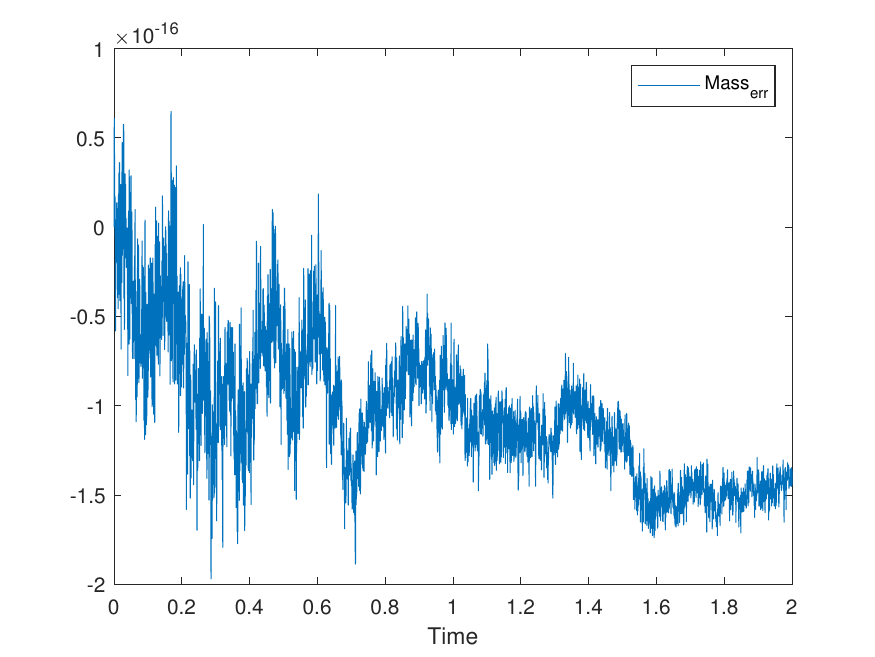}
  \caption{$Mass_{err}$ with $\nu = 10^{-3 },N =500,\tau=\frac{h}{2 \max{|f'|} }$. }
  \label{fig:1d:sin:mass}
 \end{figure}

Finally, we consider Case {3} with periodic boundary condition, which is a typical example for simulating shock formation. The initial solution and numerical solutions without using any limiter at different time are shown in Fig. \ref{fig:1d:sin:u}, where we set $\nu = 10^{-3},N =500,\tau= \frac{h}{2 \max |f'|}> \frac{h}{3 \max |f'|}$ and $\max |f'|=1$.
 It's clear that the proposed scheme \eqref{scheme:1d:full:1}--\eqref{scheme:1d:full:4} effectively captures solutions with sharp changes. This indicates that in certain situations, the BP limiter is not necessary, and the original Strang-HOC difference scheme can still perform effectively. 
To verify the mass conservation property of our original scheme without the limiter, we show the $Mass_{err}$ of the solution in Fig. \ref{fig:1d:sin:mass}.
The result indicates that the mass errors obtained with our proposed scheme approach the accuracy level of $10^{-16}$, demonstrating the conservative nature of our Strang-HOC difference scheme. This is consistent with the theoretical result in Theorem \ref{thm:mass:1d}. 
%%%%%% ==========================================
%%%%%% ==========================================
 
%%%%%% ==========================================
%%%%%% ==========================================

\begin{figure} [!htbp]
  \centering
  \includegraphics[width=6cm,height=4.5cm]{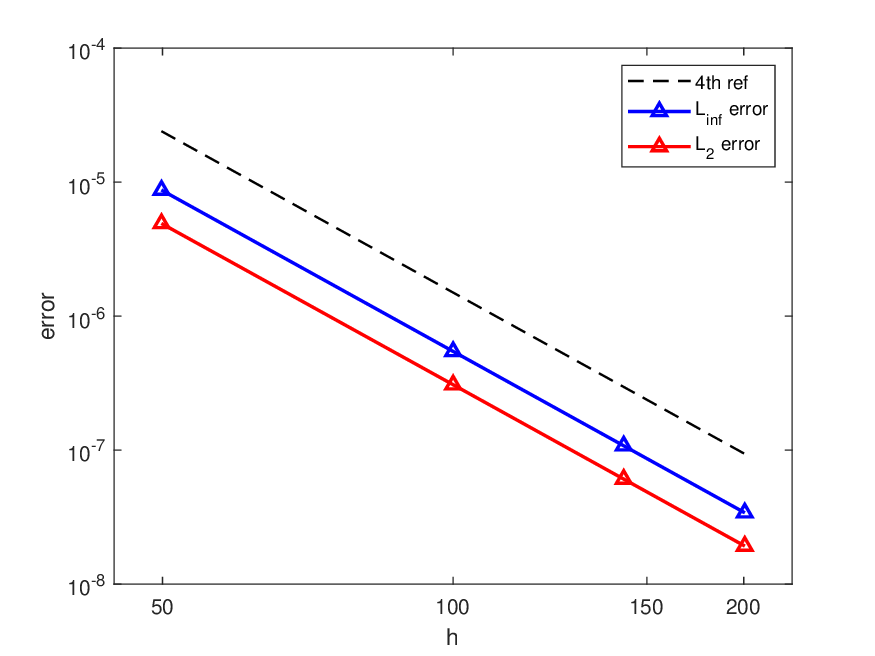}
  \caption{ Accuracy test with $\tau= h^2 $ at $T=1$. }
  \label{fig:1d:coupled_burger_rate}
 \end{figure}

 \paragraph{Example 3}
 To further gain insight into the performance of the proposed scheme, we also test it on the coupled viscous Burgers' equations:
 \begin{equation*}
  \left\{\begin{aligned}
 u_t-u_{x x}-2 u u_x+(u v)_x=0, \\ v_t-v_{x x}-2 v v_x+(u v)_x=0.
\end{aligned}\right. 
 \end{equation*}
Both the initial condition and boundary conditions are derived from exact solution \cite{kaya2001explicit} 
$$u(x, t)=v(x, t)=\exp (-t) \sin (x).$$
We compute the numerical solution over the domain $x \in[-\pi, \pi]$. We use the time step $\tau =h^2$ for the algorithm at $T=1$ to evaluate convergence in both spatial and temporal scales. It shows that the Strang-HOC difference scheme performs very well for the system and indeed holds fourth-order accurate in space and second-order accurate in time for both $L_2$ and $L_{\infty}$ norm, as shown in Fig. \ref{fig:1d:coupled_burger_rate}. 
% The results are also graphically depicted for $u(x, t)$ in Fig. \ref{fig:1d:coupled_burger_u} at the finial time $T=1$ with $\tau=\frac{h}{\max{|f'|} }>\frac{h}{3\max{|f'|} }$, where $\max{|f'|} \approx0.8$. This indicates that our scheme performs well even with very small step sizes for this example.
Due to the symmetric initial and boundary conditions, similar results are obtained for $v(x, t)$. 

%%%%%%% -----------------------------------------
% \begin{table} [!htbp]
% %\vspace{-0.1em}
% \centering
% \caption{
% Accuracy test with $\tau= {h}^2 $.
% }
% \label{tab:1d:coupled_burger:order}
% \setlength{\tabcolsep}{1.2mm}{
% \begin{tabular}{|c|cccc|}
% %\hline & \multicolumn{4}{c|}{\text {The scheme without limiter} } \\
% \hline
%  $N $ & $L^{2}$ error & Order & $ L^{\infty}$ error & Order \\
% \hline
%  50 & 8.7289 E -6 & ---
%   & 4.9150 E -6 & --- \\
% \hline
%  100 & 5.4575 E -7 & 4.00
%   & 3.0790 E -7 & 4.00 \\
% \hline
%  150 & 1.0780 E -7 & 4.00
%   & 6.0809 E -8 & 4.00  \\
% \hline
%  200 & 3.4110 E -8 & 4.00
%   & 1.9245 E -8 & 4.00 \\
% \hline
% \end{tabular}
% }
% %\vspace{-0.1em}
% \end{table}
%%%%%% ------------------------------------------------------------------

 % \begin{figure} [!htbp]
 %  \centering
 %  \includegraphics[width=6cm,height=4.5cm]{figure/1d_coupled_burgers_u.eps}
 %  \caption{Solution profiles with $N =50,\tau=\frac{h}{\max{|f'|} }$ at $T=1$. }
 %  \label{fig:1d:coupled_burger_u}
 % \end{figure}

 %%%%%%%%%%%%%%%%%%%%%%%%%%%%%%%%%%%%%%%%%%%%%%%%%
 \subsection {Two-dimensional problems }

\paragraph{Example 4} 
 
 Consider the two-dimensional Burgers equation
\begin{equation}
 u_{t} + f(u)_{x} + g(u)_{y} = \nu \left(u_{x x} +u_{y y} \right) + S(x,y,t), 
\end{equation}
with $f(u) = g(u) = \frac{1}{2}u^2$ 
and give the exact solution as follows:
$$u(x, y, t)= 2 \exp(-t) \sin (\pi x) \sin(\pi y), \quad
 (x, y, t) \in[0, 2] \times[0,2] \times[0, 2]. $$
 The source term and initial condition in both cases are given by the exact solutions.
% \begin{itemize}
%  \item Case {1}: $u(x, y, t)= 2 \exp(-t) \sin (\pi x) \sin(\pi y), \quad
%  (x, y, t) \in[0, 2] \times[0,2] \times[0, 1] $.
 
%  \item Case {2}: $
%  u(x, y, 0)=\frac{\sigma^{2} }{\sigma^{2} +2\nu t} \exp\left(-\frac{({x} -t - x_0)^{2} +({y} -t - y_0)^{2} }{2 \sigma^{2} +4 \nu t} \right), \quad
%  (x, y, t) \in[-1, 2] \times[-1,2] \times[0, 0.5] \notag
% $
% with $\sigma = 0.07, x_0=y_0=0.5$.  

% \end{itemize}
% The source term and initial condition in both cases are given by the exact solutions.
% In Case {1}, the source term and initial condition with periodic boundary condition. And in Case {2}, we only set the initial condition and the source term are zero to test the property of the bound-preserving better with homogeneous boundary conditions.
 
\begin{figure} [!htbp]
 \centering 
 \subfigure
 {
  %\label{Fig. sub. 1}
  \includegraphics[width=0.45\textwidth]{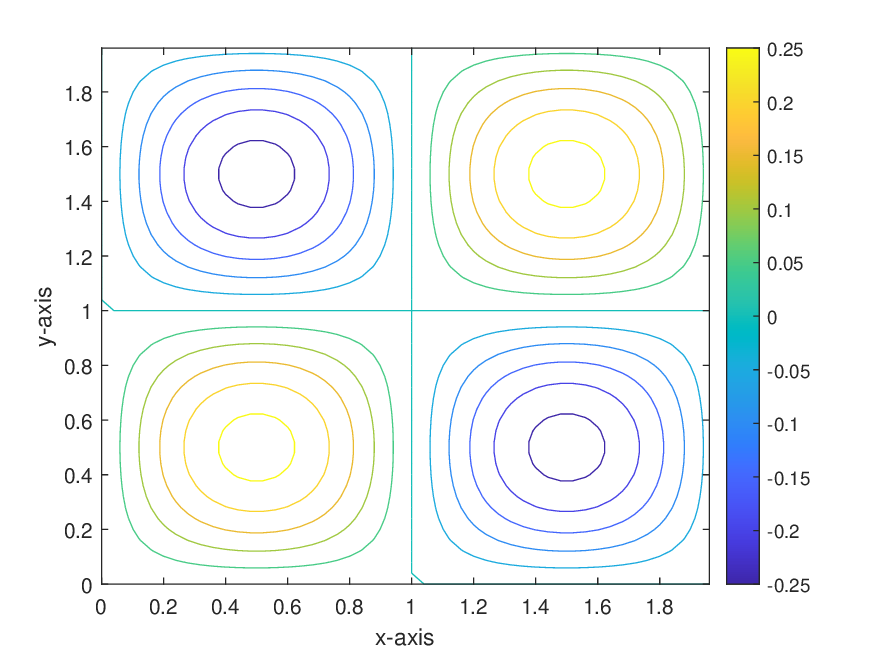}
 }
 \subfigure
 {
  %\label{Fig. sub. 2}
  \includegraphics[width=0.45\textwidth]{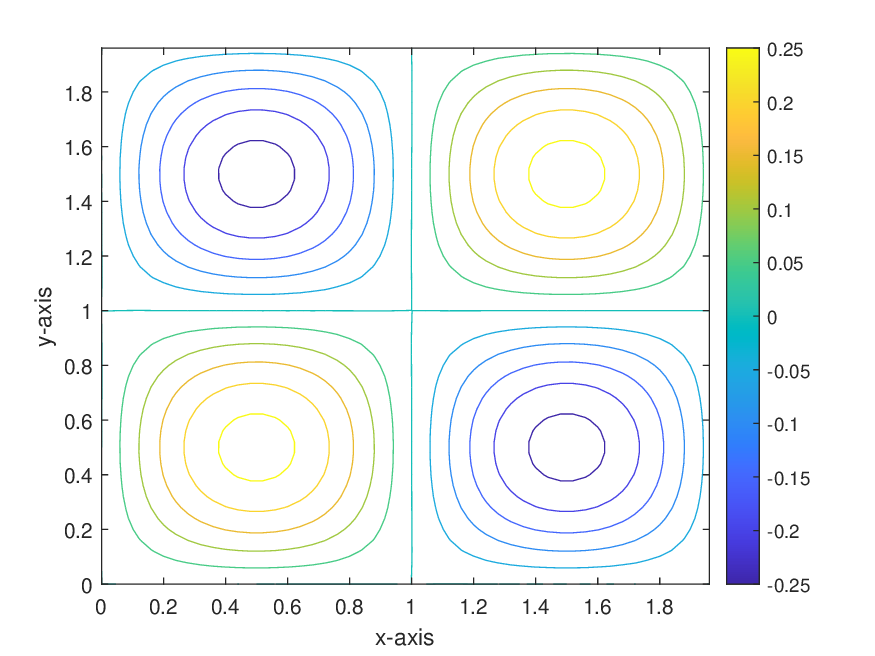}
 }
  % \subfigure
 % {
 %  %\label{Fig. sub. 2}
 %  \includegraphics[width=0.3\textwidth]{figure/2d_burger_case2_split_u.eps}
 % }
 \caption{The comparison of the exact solution (left) and obtained numerical solutions by Strang-ADI-HOC scheme (right).}
 \label{2d:burger:exp:u}
\end{figure}

\begin{figure} [!htbp]
 \centering
 \includegraphics[width=6cm,height=4.5cm]{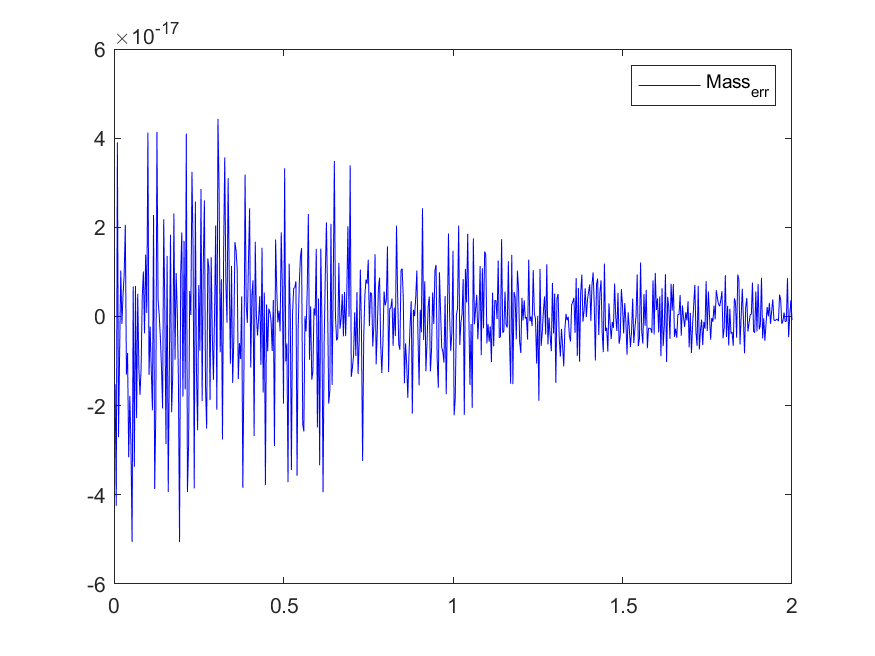}
 \caption{$Mass_{err}$ with $\nu = 5\times10^{-5 },N =100,\tau=\frac{h}{6 \max{|f'|} }$. }
 \label{2d:burger:perio:mass}
\end{figure}

\begin{figure} [!htbp]
 \centering
 \includegraphics[width=6cm,height=4.5cm]{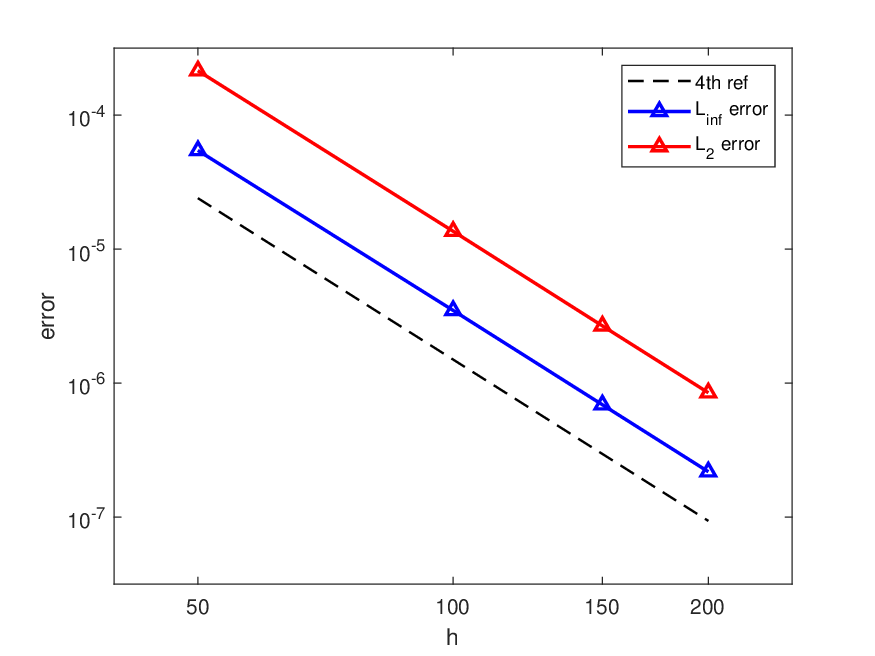}
 \caption{Convergence rates with $\nu = 5\times10^{-5 }, \tau=h^2$. }
 \label{2d:burger:perio:order}
\end{figure}

 In this example, we verify the convergence rates and mass conservation property of the  Strang-ADI-HOC difference scheme \eqref{2d:scheme:full1}--\eqref{2d:scheme:full6}.
To begin with, we present the solution graphs  with $\nu =  10^{-6}$, $N = 50$, and $\tau = \frac{h}{6 \max{|f'|}}$ in Fig. \ref{2d:burger:exp:u}, where $\max{|f'|}=2$. It is observed that our numerical scheme accurately captures the essential characteristics of the exact solution in this example, showing satisfactory agreement between the numerical and exact solutions. 
Fig. \ref{2d:burger:perio:mass} shows that the mass errors of the numerical solutions obtained through our proposed scheme approach the machine accuracy level of $10^{-17}$, demonstrating the conservative nature of the Strang-ADI-HOC difference scheme. Furthermore, numerical errors and convergence rates in discrete $L_2$ and $L_{\infty}$ norms are displayed in Table \ref{2d:burger:perio:order} with $\tau = h^2$. The results confirm that the convergence rates in $L_2$ and $L_{\infty}$ norms are indeed fourth order in space and second order in time, which aligns well with the expected results.
 
% %==============================================
% \begin{figure} [!htbp]
%  \centering
%  \subfigure%[{Numerical solution} ]
%  {
%   \includegraphics[width=0.6\textwidth]{figure/2D_burger_perio_u.eps}
%  }
%  % \subfigure[ Pointwise error distributions ]
%  % {
%  %  \includegraphics[width=0.45\textwidth]{figure/2D_burger_perio_err.eps}
%  % }
%  \caption{ Numerical solution and pointwise error distributions of the Strang-ADI-HOC difference scheme .}
%  \label{2d:burger:perio:u}
% \end{figure}

%==============================================

%==========================================
%==========================================

%==========================================
\paragraph{Example 6} 
In this test, we consider the linear variation coefficient convection diffusion problem %\cite{qin2023}
\begin{align*}
  u_t + c_1 u_ x + c_2 u_y - \nu (u_{xx}+ u_{yy}) =0, \quad 0\leq x,y \leq 1.
\end{align*}
 and give different forms of solutions and the velocity field as follows. 
\begin{itemize}
  \item Case {1} \cite{FU2017jsc}: 
The first example is the Gaussian pulse with a velocity field $\left(c_1,c_2\right)=(-4 y, 4 x)$ and the initial concentration distribution and boundary conditions obtained by the analytically exact solution is given as :
\begin{align*}
u(x, y, t)=\frac{ \sigma^2}{ \sigma^2+2 \nu t} \exp \Big(-\frac{\left(\bar{x}(t)-x_0\right)^2+\left(\bar{y}(t)-y_0\right)^2}{2 \sigma^2+4 \nu t}\Big),
\end{align*}
where $\bar{x}(t)=x \cos (4 t)+y \sin (4 t), \bar{y}(t)=-x \sin (4 t)+y \cos (4 t)$.
Here, $\left(x_0, y_0\right)=(-0.35,0)$ represents the initial center of the Gaussian pulse, and $\sigma^2=5 \times 10^{-4}$ is the standard deviation. 

  \item Case {2} \cite{FU2019sisc}: 
In this example, the transport of another Gaussian hump in a vortex velocity field is considered. The velocity field is divergence-free, $(c_1,c_2) = (\psi_y, \psi_x)$, where $\psi (x, y) = \frac{1}{2}\exp(\sin(\pi x))\exp(\sin(\pi y))$ in the domain. 
The initial distribution of the Gaussian hump is given as 
\begin{align*}
u(x, y, t)= \exp \Big(-\frac{\left({x} -x_0\right)^2+\left({y} -y_0\right)^2}{2 \sigma^2}\Big),
\end{align*}
where the center of the initial Gaussian hump is specified as $x_0 = 0.75$, $y_0 = 0.5$ and $\sigma^2=1.6 \times 10^{-3}$.

  \item Case {3} \cite{qin2023}: 
  In this case, the initial value is given as
\begin{align*}
u_{0} (x, y)= \begin{cases} 1, & \text {for } 0.1 \leq x \leq 0. 3, 0.1 \leq y \leq 0. 3, \\ 0, & \text {otherwise. } \end{cases}
\end{align*}
with $\sigma = 0.07, x_0=y_0=0.5$ and the velocity field of $(c_1,c_2) = (\sin(\pi y), \sin(\pi x))$. For simplicity, the periodic boundary condition is employed.
  
\end{itemize}

In Case {1}, we test the accuracy and property of bound-preserving and mass-preserving of the proposed Strang-ADI-HOC scheme with BP limiter. To illustrate the second-order accuracy in time and the fourth-order accuracy in space at the same time, we set $\tau = h^2$ with $\nu = 5 \times 10^{-3}$ and $T=0.25$.  
 From Table \ref{tab:roate:order:BP}, we can observe that the convergence rate of the Strang-ADI-HOC difference scheme with and without the BP limiter tends to fourth order accuracy. In addition, negative values appear  if the BP technique is not applied, and the limiter remedied the negative values in a conservative way. This once again clearly illustrates that the BP technique does not destroy the accuracy when it works, which is consistent with our theoretical result in Remark \ref{rem:order}.

\begin{figure} [!htbp]
 \centering 
 \subfigure[] %[(a)]
 {
  % \label{fig:burger:ini:0:sub1}
  \includegraphics[width=0.3\textwidth]{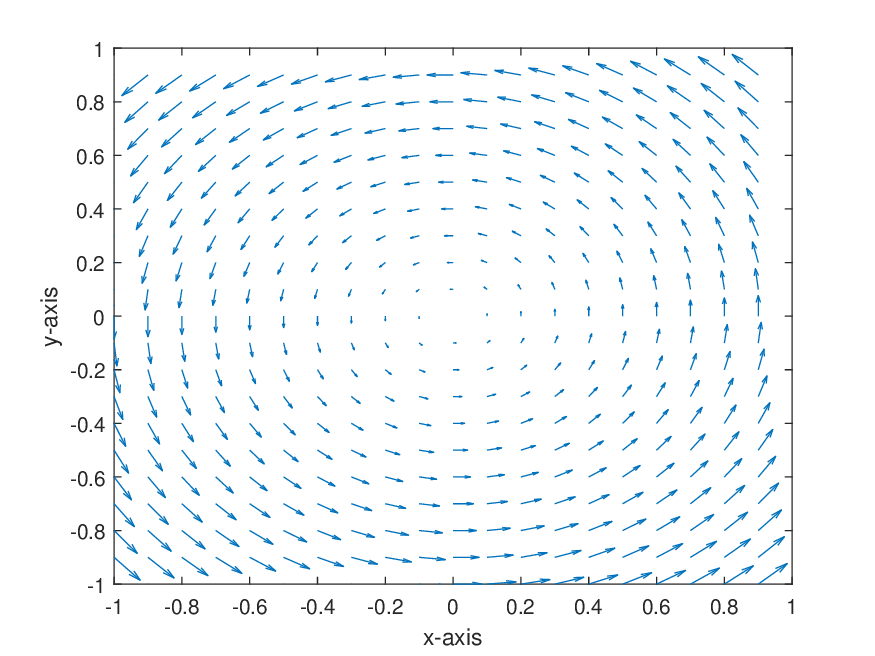}
 }
 \subfigure[]%[(b)]
 {
  %\label{Fig. sub. 2}
  \includegraphics[width=0.3\textwidth]{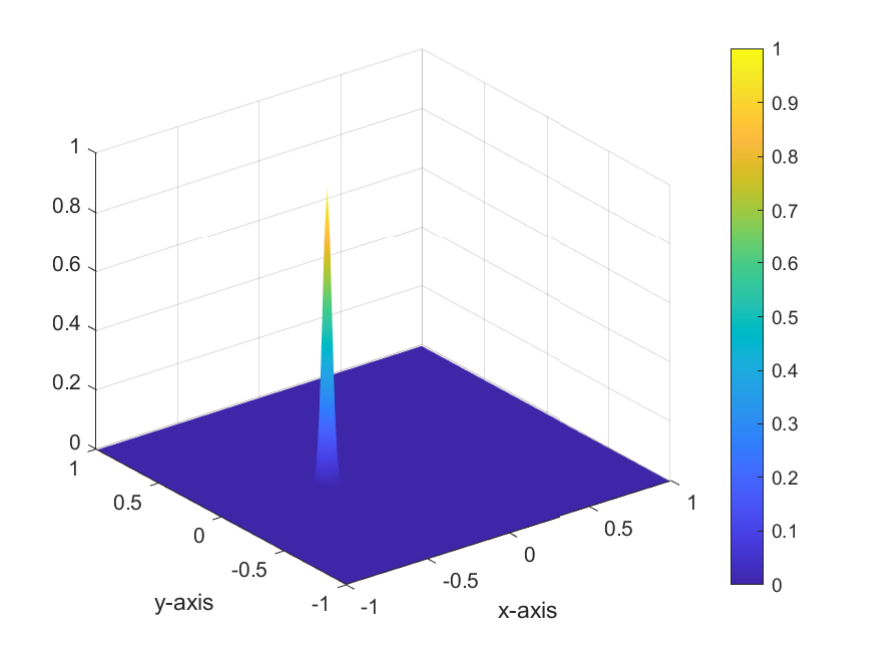}
 }
  \subfigure[]%[(c)]
 {
  %\label{Fig. sub. 2}
  \includegraphics[width=0.3\textwidth]{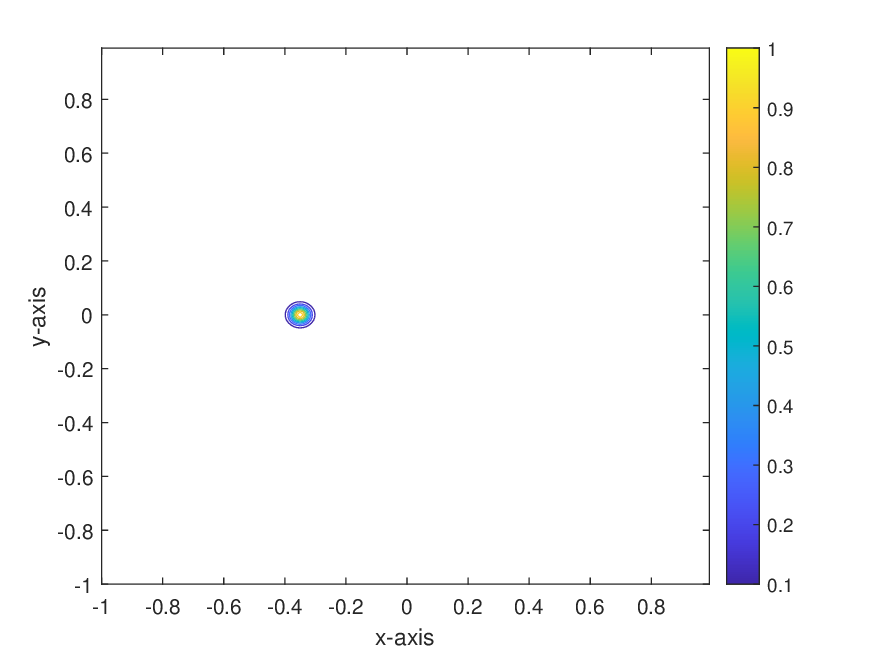}
 }
 \caption{ (a) The velocity field of $(c_1,c_2) = (-4y, 4x)$. (b) The three-dimensional surface plot of the initial value distribution. (c) The corresponding top view of the initial value distribution.}
 \label{fig:burger:4yx_UO}
\end{figure}

\begin{figure} [!htbp]
 \centering 
 \subfigure[Without the BP limiter]
 {
  % \label{fig:burger:ini:sub1}
  \includegraphics[width=0.3\textwidth]{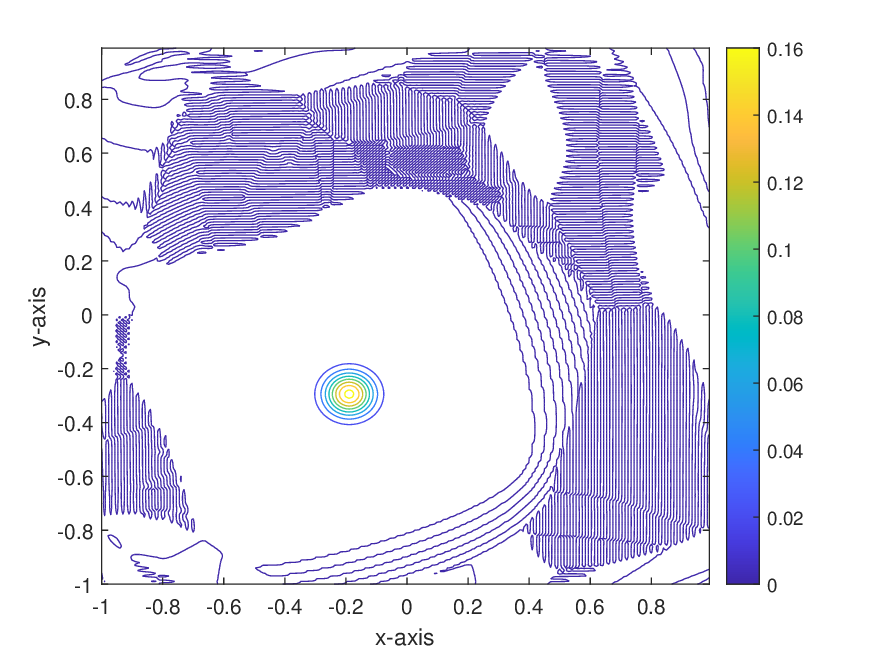}
 }
 \subfigure[With the BP limiter]
 {
  % \label{fig:burger:ini:sub2}
  \includegraphics[width=0.3\textwidth]{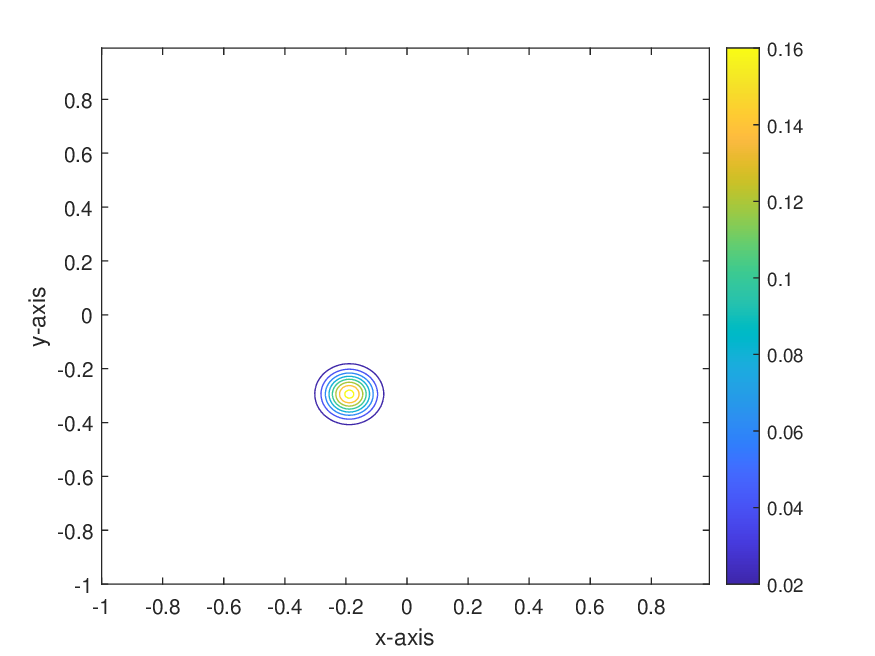}
 }
  \subfigure[Exact solution]
 {
   % \label{fig:burger:ini:sub3}
  \includegraphics[width=0.3\textwidth]{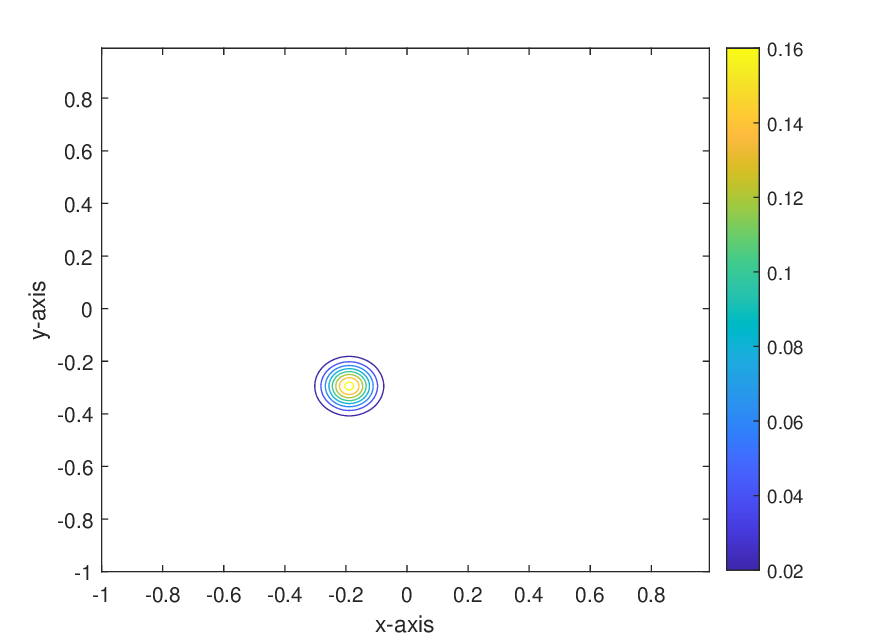}
 }
 \caption{ Comparison  of contour plots obtained by the solutions without/with the BP limiter and the exact solutions.}
 \label{fig:burger:4yx:comprare_ref:contour}
\end{figure}

\begin{figure} [!htbp]
 \centering 
 \subfigure[Without the BP limiter]
 {
  % \label{fig:burger:ini:sub1}
  \includegraphics[width=0.45\textwidth]{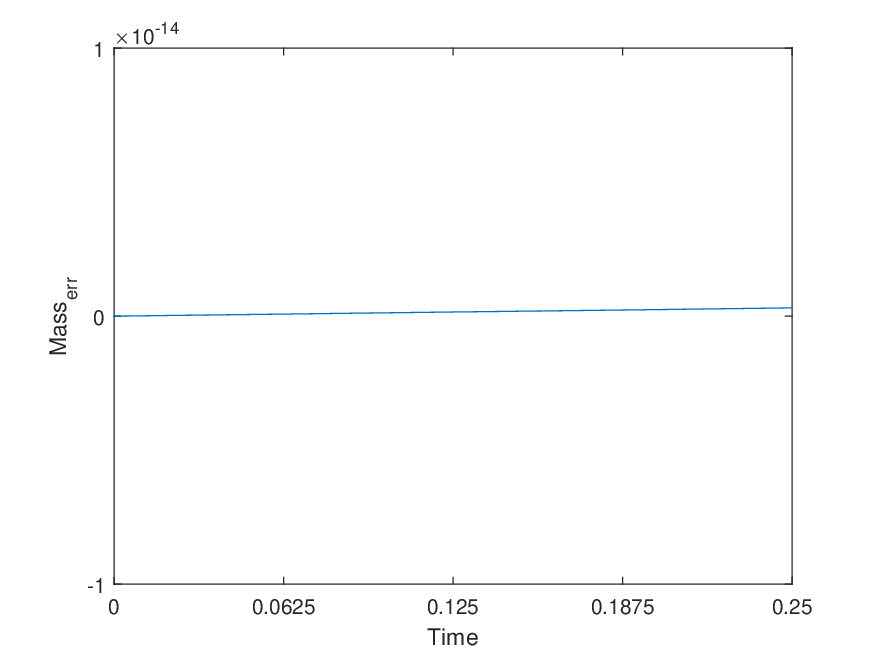}
 }
 \subfigure[With the BP limiter]
 {
  % \label{fig:burger:ini:sub2}
  \includegraphics[width=0.45\textwidth]{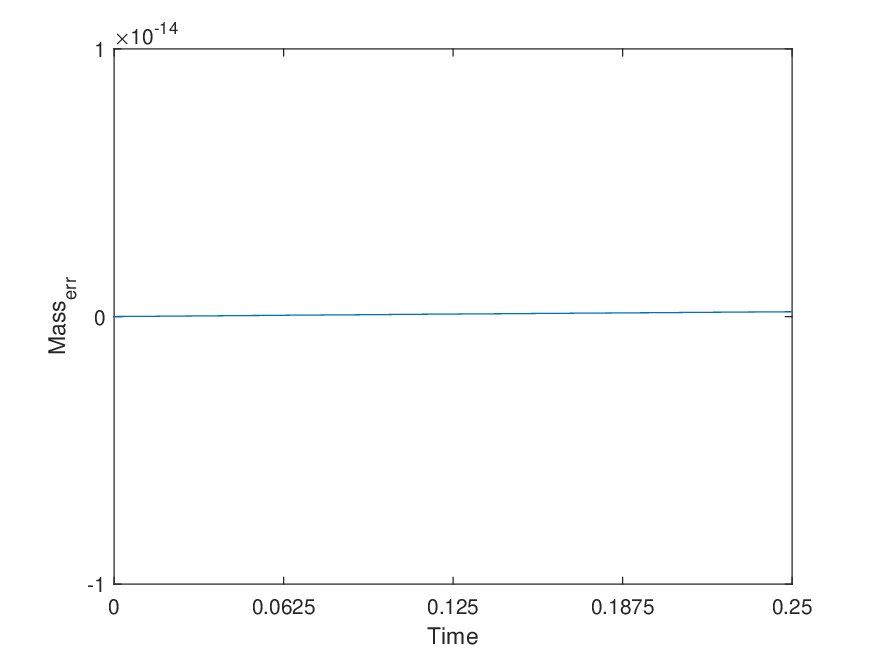}
 }
 %  \subfigure[Exact solution]
 % {
 %  % \label{fig:burger:ini:sub3}
 %  \includegraphics[width=0.3\textwidth]{figure/2d_burger_4yx_mass_U.eps}
 % }
 \caption{ Comparison  of mass obtained by the solutions without/with the BP limiter.}
 \label{fig:burger:4yx:mass}
\end{figure}

 To better observe the effect of the BP technique, we present a series of plots. First, the divergence-free vortex shear and contour plots at the initial time $T=0$ are illustrated in Fig. \ref{fig:burger:4yx_UO}. In Fig. \ref{fig:burger:4yx:comprare_ref:contour}, we compare the surface and contour plots without/with the BP limiter and the exact solution with $\nu = 5 \times 10^{-3}, \tau = \frac{h}{6 {\max{|f'|} } }$ at $T=0.25$, where $\max{|f'|}\approx 4$. 
It is clear that the numerical solution calculated without the BP limiter exhibits negative values and oscillations, which pollute the numerical results. In contrast, the numerical solution calculated with the BP limiter shows good consistency with the exact solutions.
Finally, the time evolution of the $Mass_{err}$ without/with the BP limiter is plotted in Fig. \ref{fig:burger:4yx:mass}. It is observed that mass conservation is well-preserved numerically in both cases, with or without the BP limiter, which is consistent with our numerical results.

 %============================================
\begin{table} [!htbp]
%\vspace{-0.1em}
\centering
\caption{
Accuracy test with $\nu =5 \times 10^{-3}, \tau= {h}^2 $.
}
\label{tab:roate:order:BP}
 \scalebox{0.95}{%此处放置命令
\setlength{\tabcolsep}{1.0mm}{
\begin{tabular}{|c|ccccc|ccccc|}
\hline & \multicolumn{5}{c|}{\text {The scheme without the BP limiter} }
  & \multicolumn{5}{c|}{\text {The scheme with the BP limiter } } \\
\hline
 $N$ & $L^{2}$ error & Order & $ L^{\infty}$ error & Order & $m_{err}$ 
  & $L^{2}$error & Order & $ L^{\infty}$ error & Order & $m_{err}$\\
\hline
 50 & 3.9493 E -4 & ---
   & 3.8901 E -3 & ---
   & -8. 9772 E -5   
   & 8.0875 E -4 & ---
   & 6.3116 E -3 & --- 
   & 0\\
\hline
 100 & 1.5794 E -5 & 4. 64
   & 1.6433 E -4 & 4. 57
   & -3.0292 E -12
   & 1.5477 E -5 & 5. 70
   & 1.5550 E -4 & 5. 34 
   & 0\\
\hline
 200 & 9.4946 E -7 & 4. 06
   & 9.7788 E -6 & 4. 07
   & -3.9502 E -30 
   & 9.4946 E -7 & 4. 03
   & 9.7788 E -6 & 4. 00 
   & 0\\
\hline
 300 & 1.8625 E -7 & 4. 02
   & 1.9315 E -6 & 4. 00
   & -2.1723 E -43
   & 1.8625 E -7 & 4. 02
   & 1.9315 E -6 & 4. 00 
   & 0\\
\hline
\end{tabular}
}
}%这里要包回来
%\vspace{-0.1em}
\end{table}
%=========================================
%
%
%==========================================

%================================================
%===============================================

\begin{figure} [!htbp]
 \centering 
 \subfigure[] %[(a)]
 {
 % \label{fig:burger:ini:0:sub1}
  \includegraphics[width=0.3\textwidth]{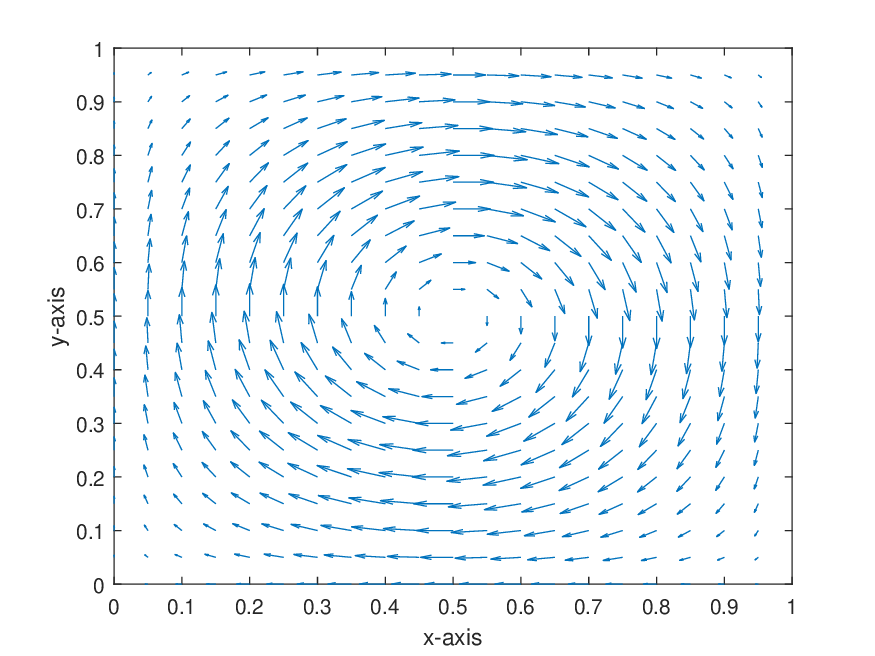}
 }
 \subfigure[]%[(b)]
 {
  %\label{Fig. sub. 2}
  \includegraphics[width=0.3\textwidth]{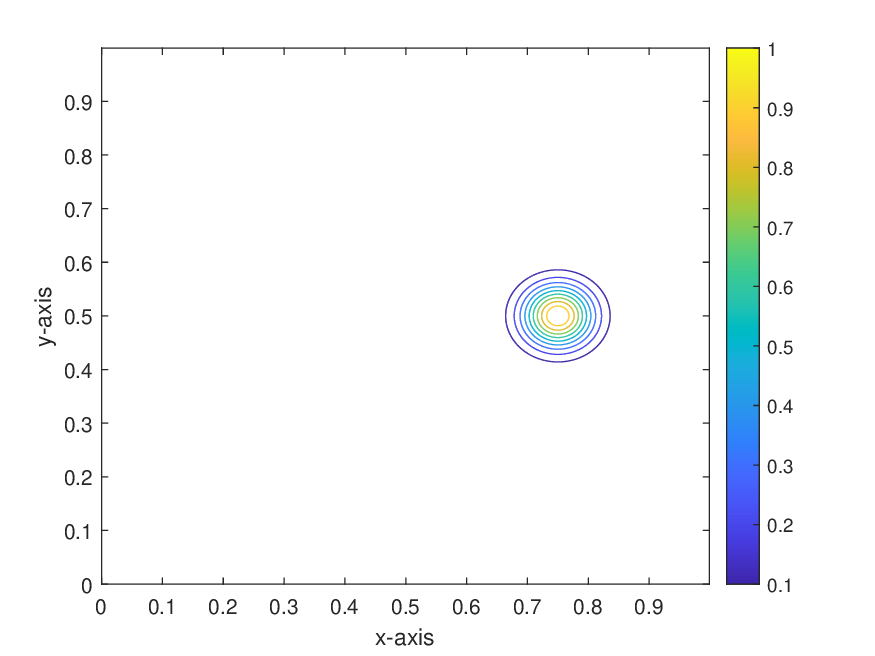}
 }
 \subfigure[]%[(c)]
 {
  %\label{Fig. sub. 2}
  \includegraphics[width=0.3\textwidth]{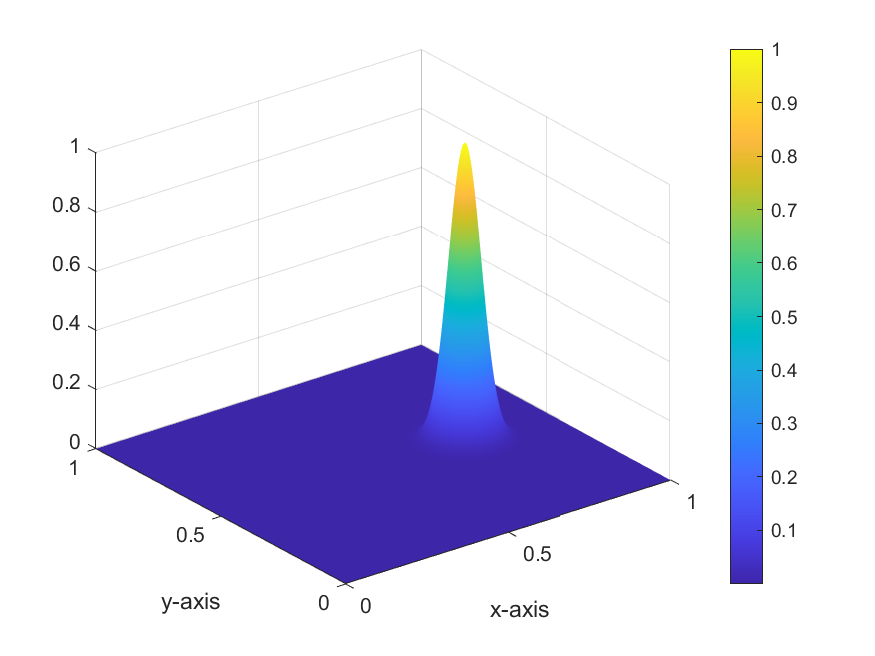}
 }
 \caption{ (a) The velocity field of $(c_1,c_2) = (\psi_y, \psi_x)$, where $\psi (x, y) = \frac{1}{2}\exp(\sin(\pi x))\exp(\sin(\pi y))$. (b) The contour plot of the initial value distribution. (c) The three-dimensional surface plot of the initial value distribution.}
 \label{fig:roate_phi:initial}
\end{figure}

\begin{figure} [!htbp]
 \centering 
 \subfigure[$T=0.25$]
 {
 % \label{fig:roate_case3:u:1}
  \includegraphics[width=0.3\textwidth]{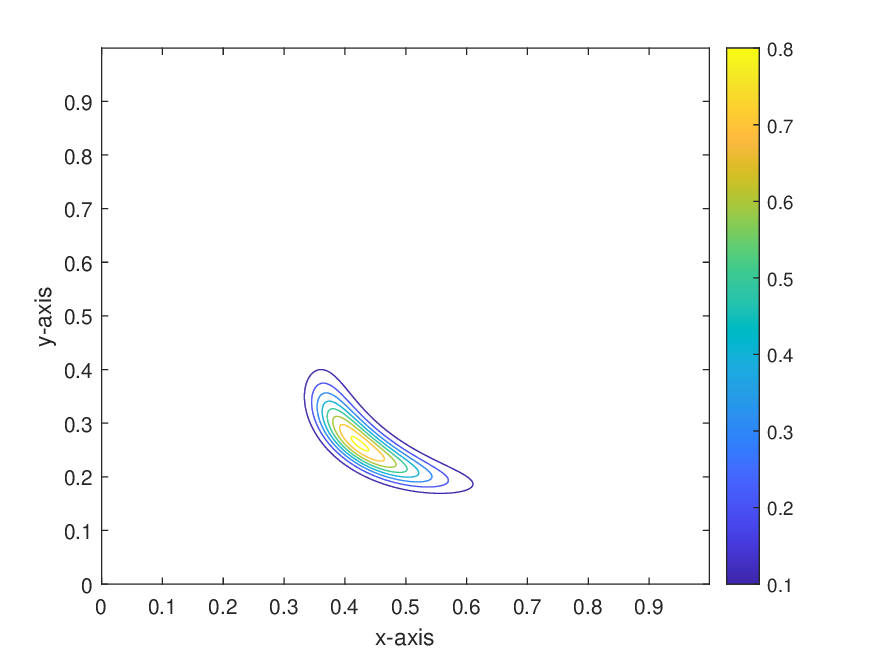}
 }
 \subfigure[$T=0.5$]
 {
 % \label{fig:roate_case3:u:4}
  \includegraphics[width=0.3\textwidth]{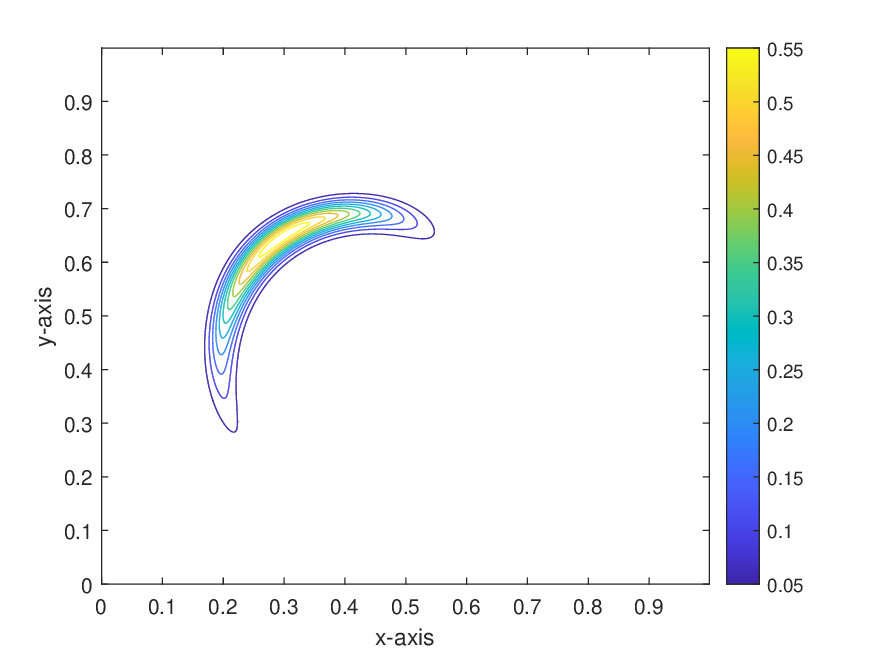}
 }
  \subfigure[$T=0.75$]
 {
 % \label{fig:roate_case3:u:7}
  \includegraphics[width=0.3\textwidth]{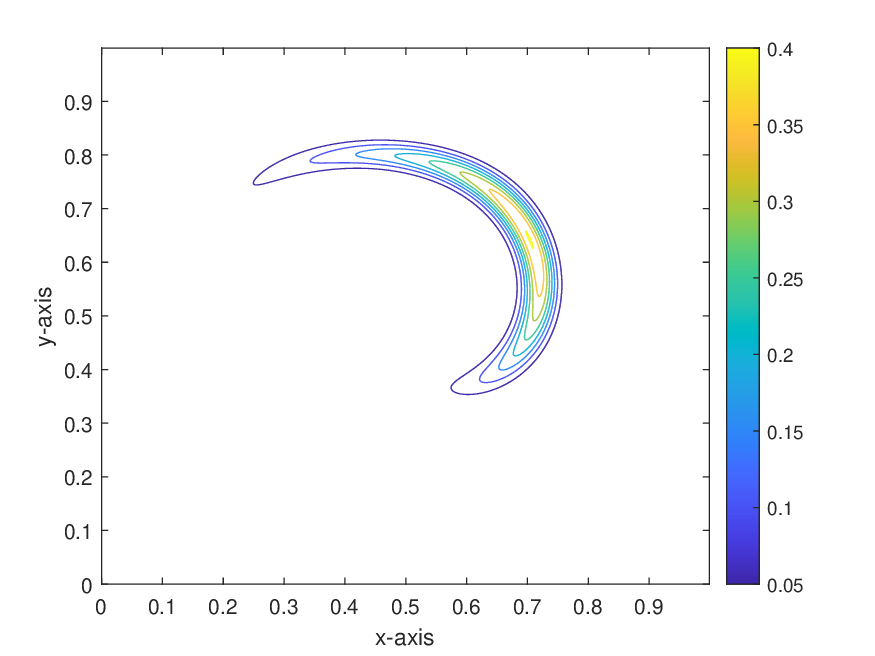}
 }
 \subfigure[$T=0.25$]
 {
 % \label{fig:roate_case3:u:2}
  \includegraphics[width=0.3\textwidth]{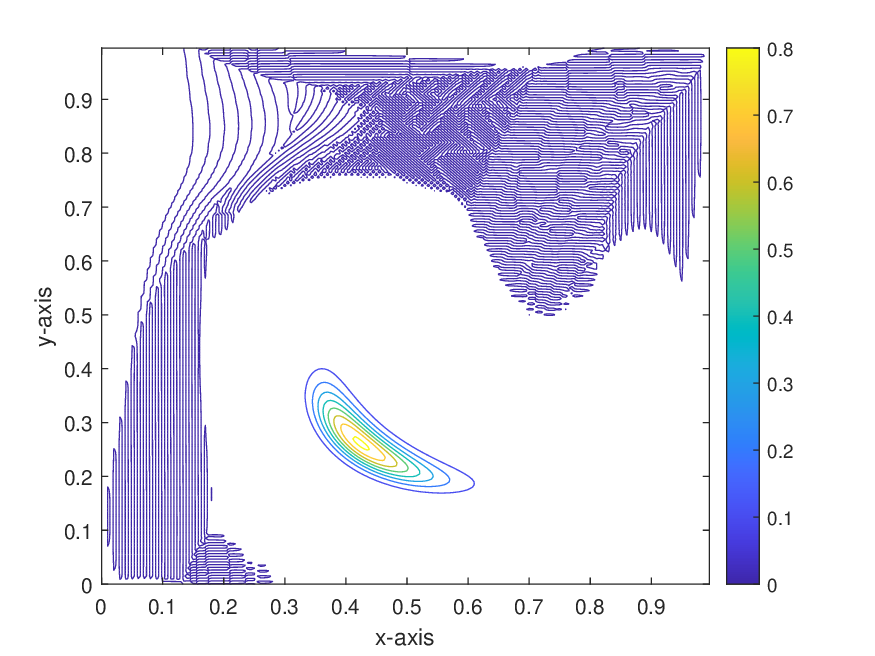}
 }
 \subfigure[$T=0.5$]
 {
 % \label{fig:roate_case3:u:5}
  \includegraphics[width=0.3\textwidth]{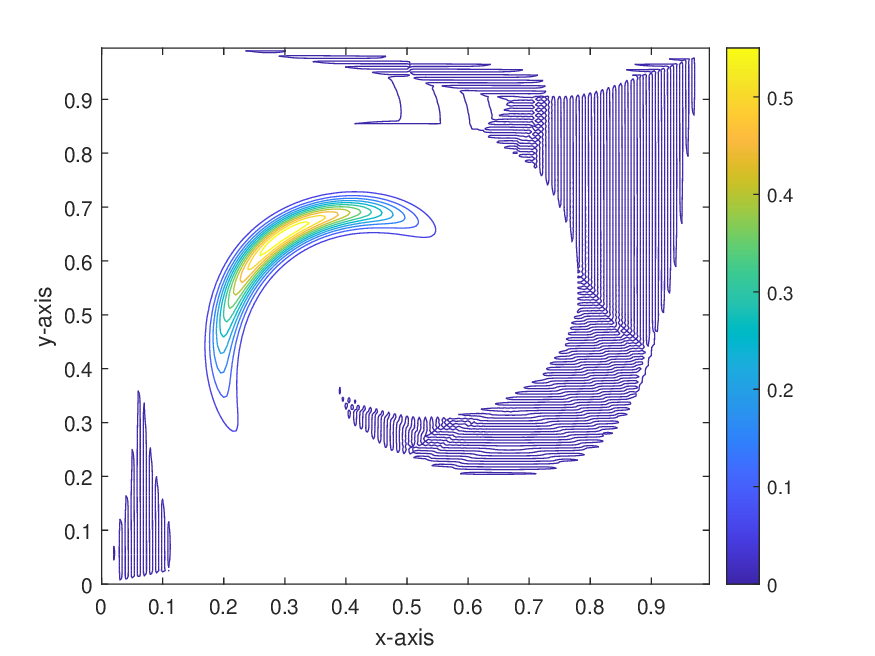}
 }
  \subfigure[$T=0.75$]
 {
 % \label{fig:roate_case3:u:8}
  \includegraphics[width=0.3\textwidth]{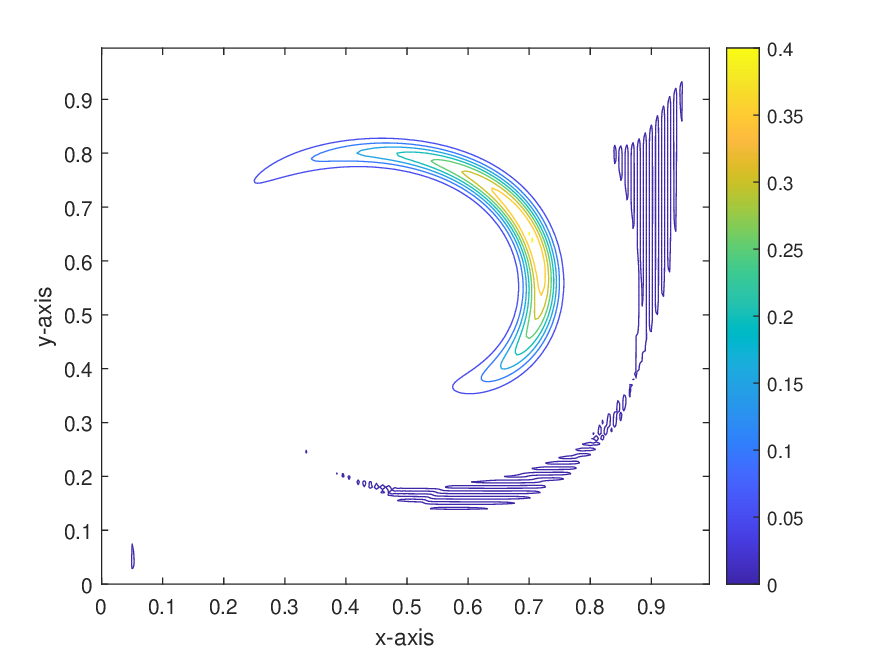} 
 }
 \subfigure[$T=0.25$]
 {
 % \label{fig:roate_case3:u:3}
  \includegraphics[width=0.3\textwidth]{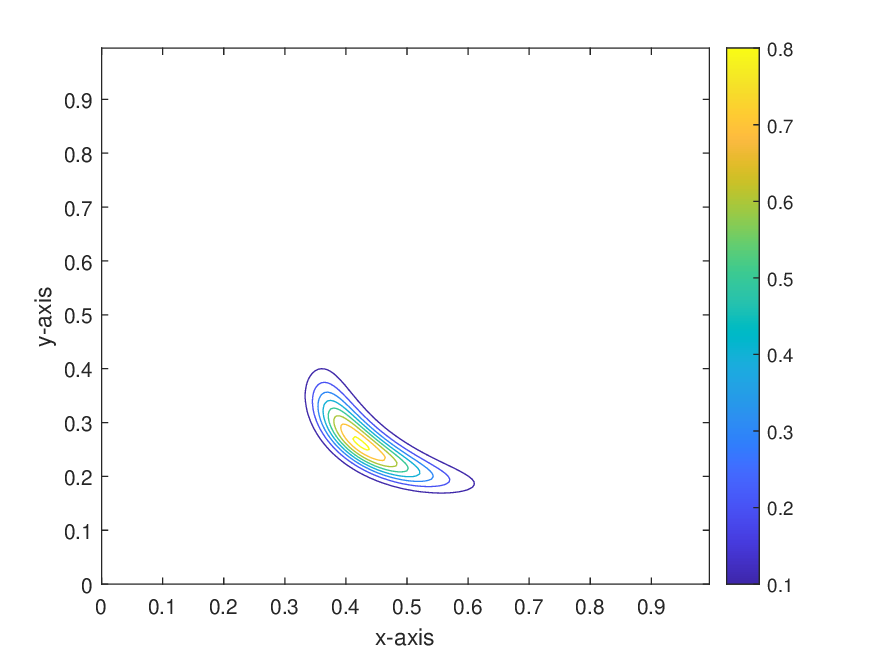}
 }
 \subfigure[$T=0.5$]
 {
 % \label{fig:roate_case3:u:6}
  \includegraphics[width=0.3\textwidth]{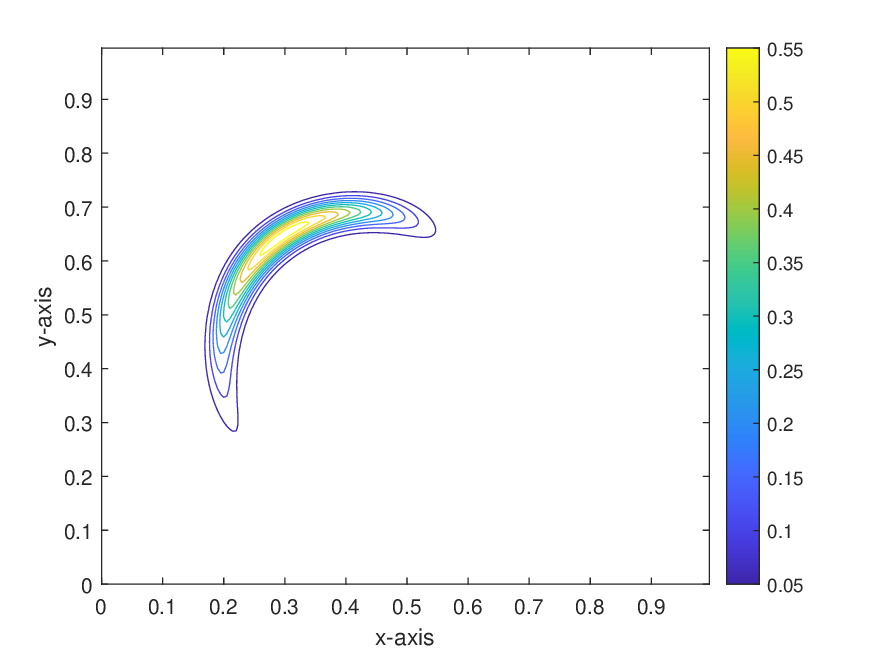}
 }
 \subfigure[$T=0.75$]
 {
 % \label{fig:roate_case3:u:9}
  \includegraphics[width=0.3\textwidth]{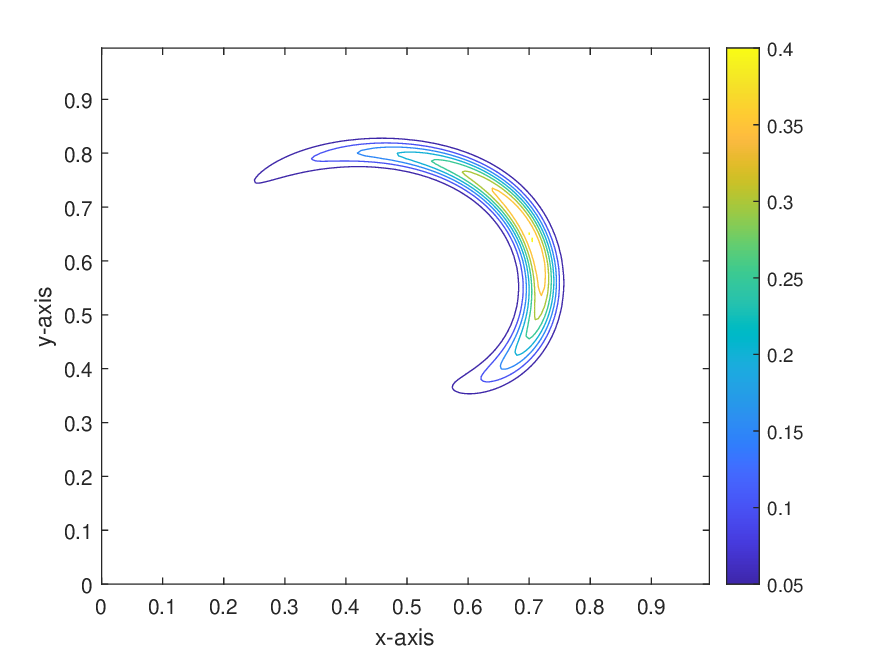}
 }
 \caption{Comparison  of contour plots between the reference solution (a-c) and the obtained numerical solution by HOC-ADI-Splitting scheme without  (d-f) and with  the BP limiter (h-g).}
 \label{fig:roate:phi:u:contour} 
\end{figure}

\begin{figure} [!htbp]
 \centering 
 \subfigure[$T=0.25$]
 {
  % \label{fig:roate_case3:u:1}
  \includegraphics[width=0.3\textwidth]{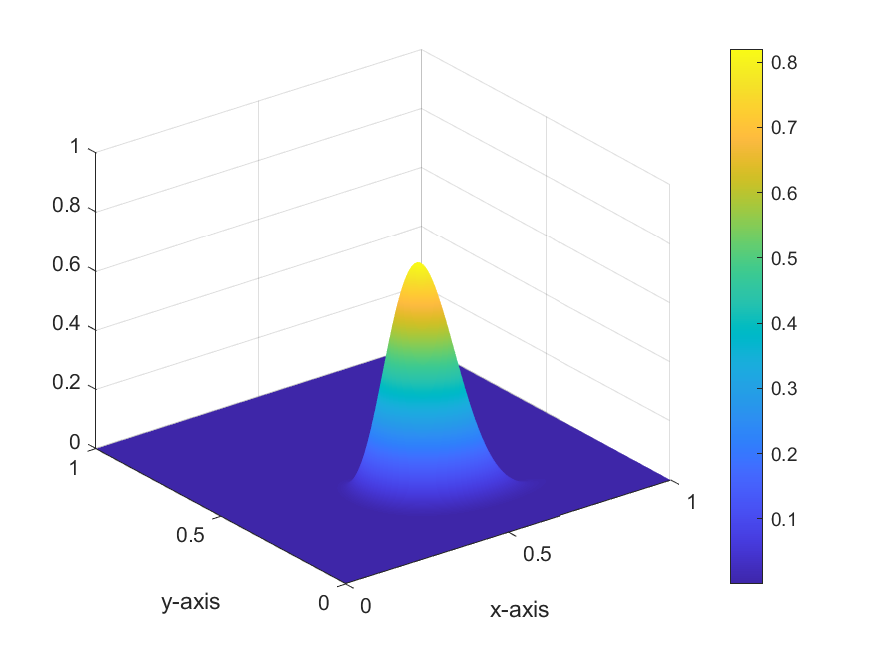}
 }
 \subfigure[$T=0.5$]
 {
  % \label{fig:roate_case3:u:4}
  \includegraphics[width=0.3\textwidth]{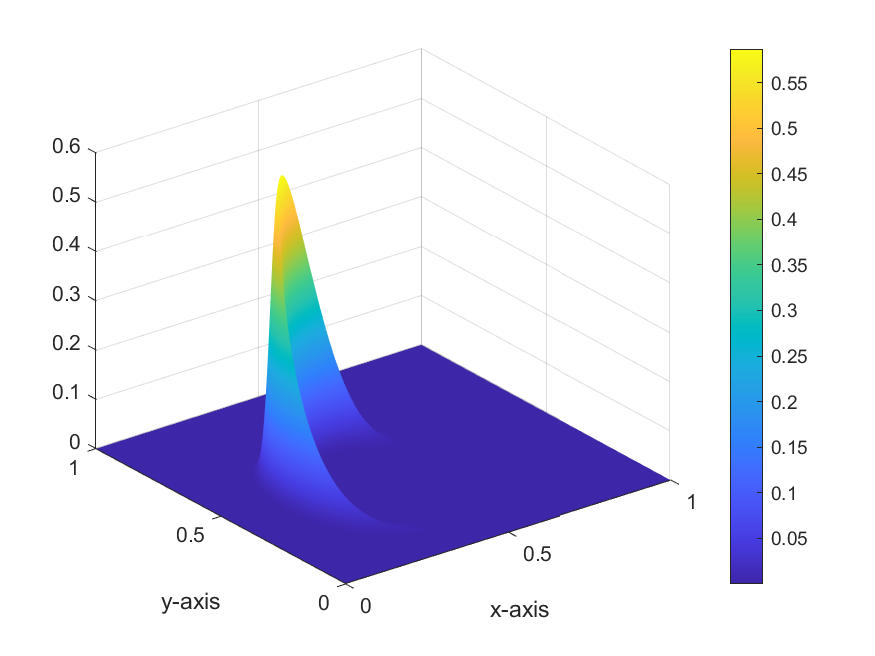}
 }
 \subfigure[$T=0.75$]
 {
  % \label{fig:roate_case3:u:7}
  \includegraphics[width=0.3\textwidth]{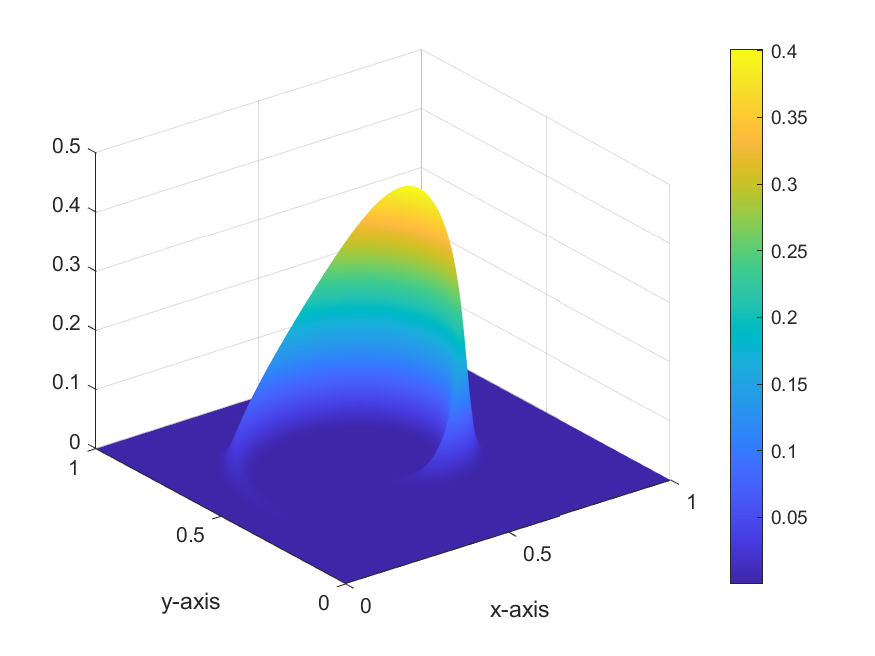}
 }
 \subfigure[$T=0.25$]
 {
  % \label{fig:roate_case3:u:2}
  \includegraphics[width=0.3\textwidth]{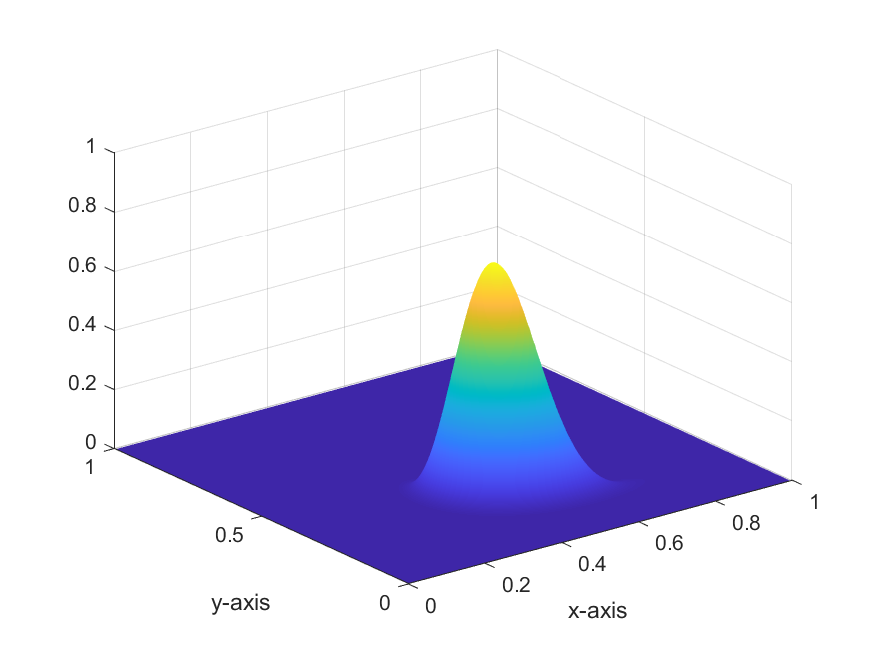}
 }
  \subfigure[$T=0.5$]
 {
  % \label{fig:roate_case3:u:5}
  \includegraphics[width=0.3\textwidth]{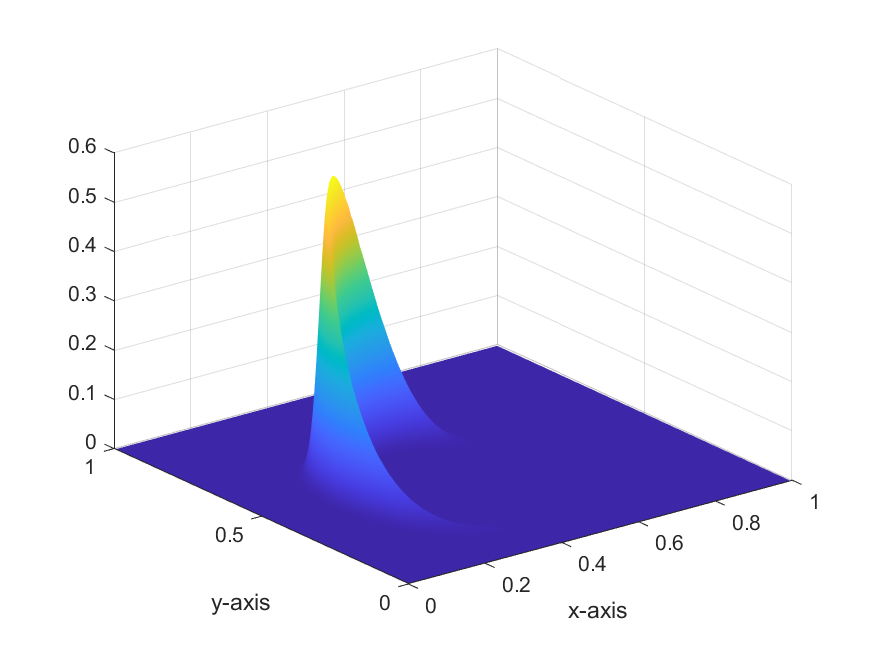}
 }
 \subfigure[$T=0.75$]
 {
  % \label{fig:roate_case3:u:8}
  \includegraphics[width=0.3\textwidth]{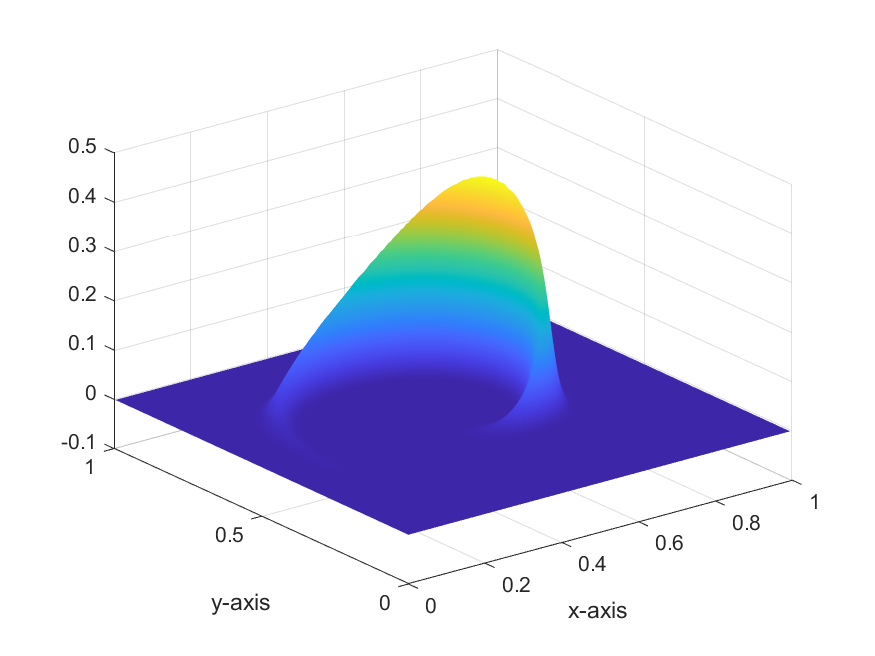}
 }
 \subfigure[$T=0.25$]
 {
  % \label{fig:roate_case3:u:3}
  \includegraphics[width=0.3\textwidth]{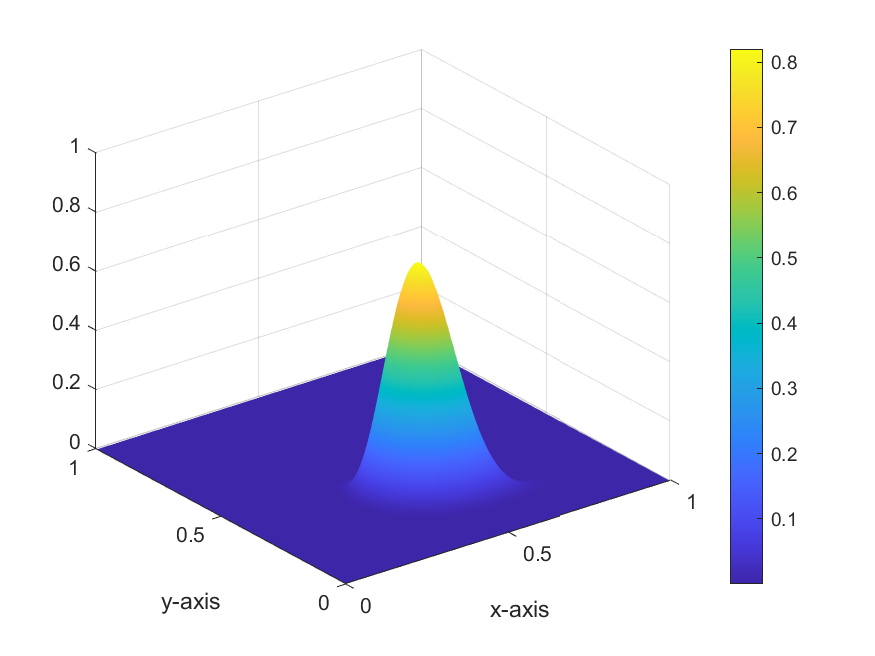}
 }
 \subfigure[$T=0.5$]
 {
  % \label{fig:roate_case3:u:6}
  \includegraphics[width=0.3\textwidth]{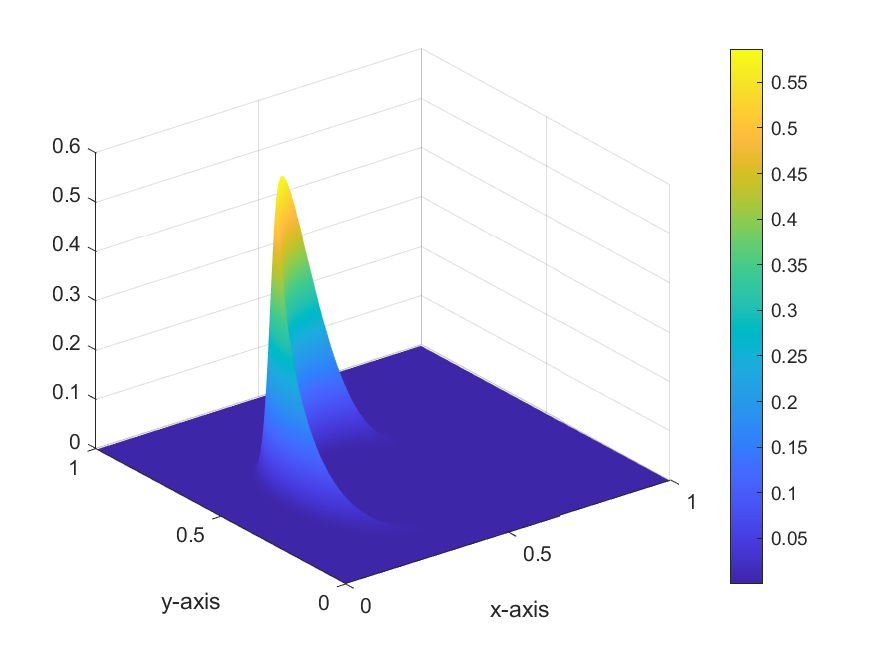}
 }
 \subfigure[$T=0.75$]
 {
  % \label{fig:roate_case3:u:9}
  \includegraphics[width=0.3\textwidth]{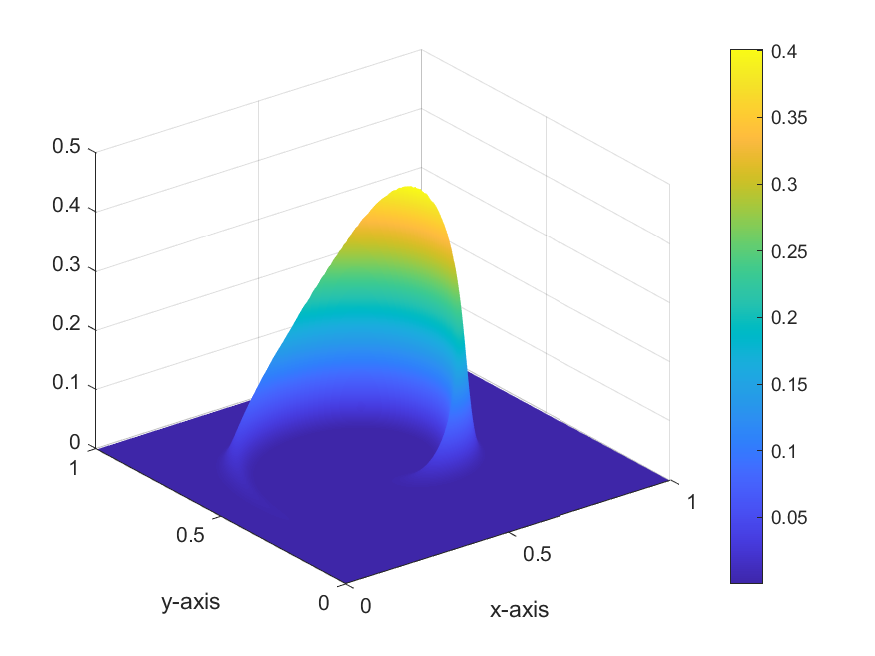}
 }
 \caption{Comparison  of three-dimensional views between the obtained numerical solution by the reference solution (a-c), HOC-ADI-Splitting scheme without  (d-f) and with the BP limiter (g-i).}
 \label{fig:roate:phi:u:surf}
\end{figure}

In Case {2}, the transport of a Gaussian hump in a vortex velocity field is considered. The velocity field and the initial value distribution are presented in Fig. \ref{fig:roate_phi:initial}. 
 Fig. \ref{fig:roate:phi:u:contour} and \ref{fig:roate:phi:u:surf} visually compare the contour plots  of the moving square wave computed by the HOC-ADI-Splitting scheme without/with the BP limiter at different time. We set $\nu = 5 \times 10^{-4}, N =200$ and $\tau =\frac{h}{4 {\max{|f'|} }} \leq \frac{h}{6 {\max{|f'|} } }, \max{|f'|}=5 $ at different time. The distribution solved using a fine mesh of $N=1000$
 and $\tau =\frac{h}{6 {\max{|f'|} } }$
is used as reference for the exact solution.
This comparison highlight the growing performance gap: as time progresses, the numerical solution with the BP limiter gets a perfect match with the exact solution, in which the scalar quantity gradually being stretched by the swirling flow over time. In contrast, the solution without the limiter shows significant numerical dissipation of concentration, leading to a less accurate representation of the flow.
 This example also demonstrates that the condition stated in  Theorem \ref{th:2d:all} is sufficient but not necessary for preserving the bound. In practice, the net ratio condition is more lenient and can be effectively applied in actual calculations.

%==========================================

%================================================

%================================================

\begin{figure} [!htbp]
 \centering 
 \subfigure[] %[(a)]
 {
   \label{fig:burger:ini:0:sub1}
  \includegraphics[width=0.3\textwidth]{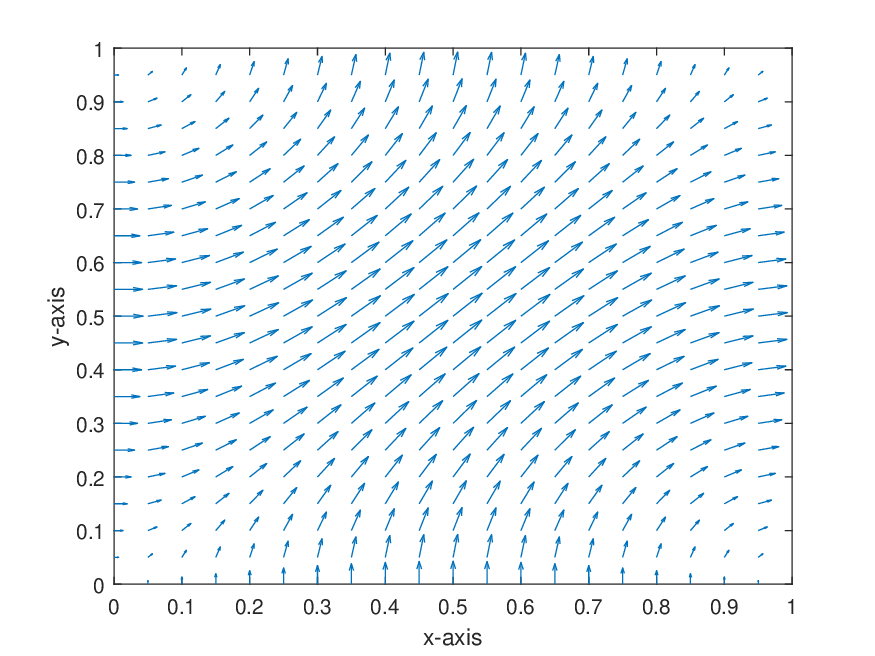}
 }
 \subfigure[]%[(b)]
 {
  %\label{Fig. sub. 2}
  \includegraphics[width=0.3\textwidth]{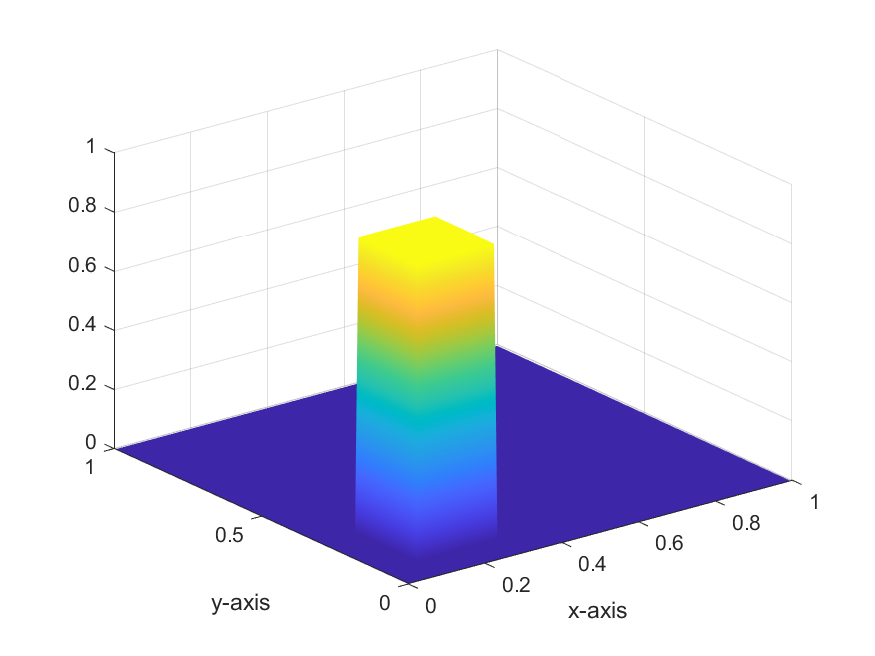}
 }
  \subfigure[]%[(c)]
 {
  %\label{Fig. sub. 2}
  \includegraphics[width=0.3\textwidth]{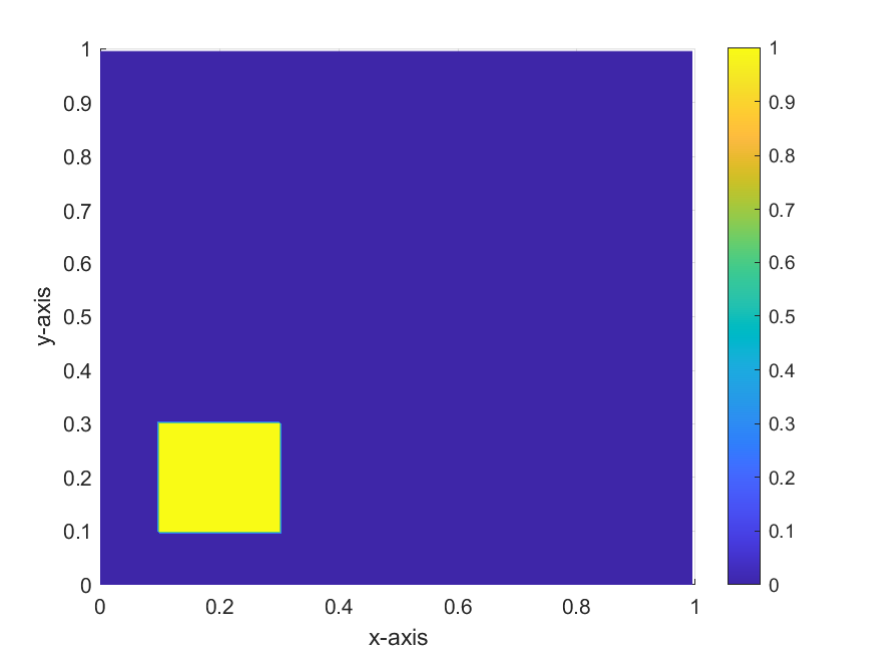}
 }
 \caption{ (a) The velocity field of $(c_1,c_2) = (\sin(\pi y), \sin(\pi x))$. (b) The three-dimensional surface plot of the initial value distribution. (c) The corresponding top view of the initial value distribution.}
 \label{fig:burger:ini:0}
\end{figure}

\begin{figure} [!htbp]
 \centering 
 \subfigure[$T=0.25$]
 {
  % \label{fig:burger:ini:sub1}
  \includegraphics[width=0.4\textwidth]{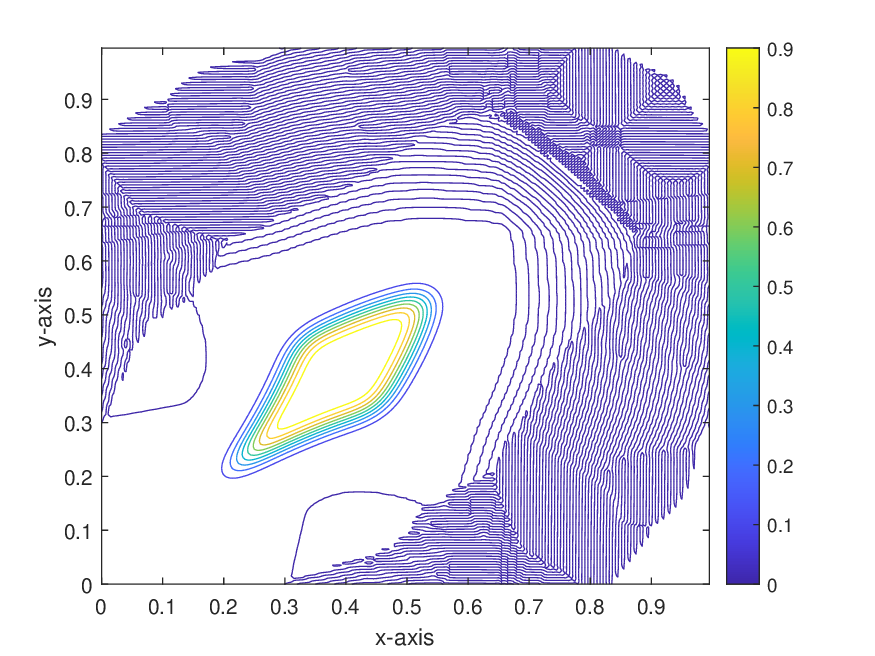}
 }
  \subfigure[$T=0.5$]
 {
  % \label{fig:burger:ini:sub3}
  \includegraphics[width=0.4\textwidth]{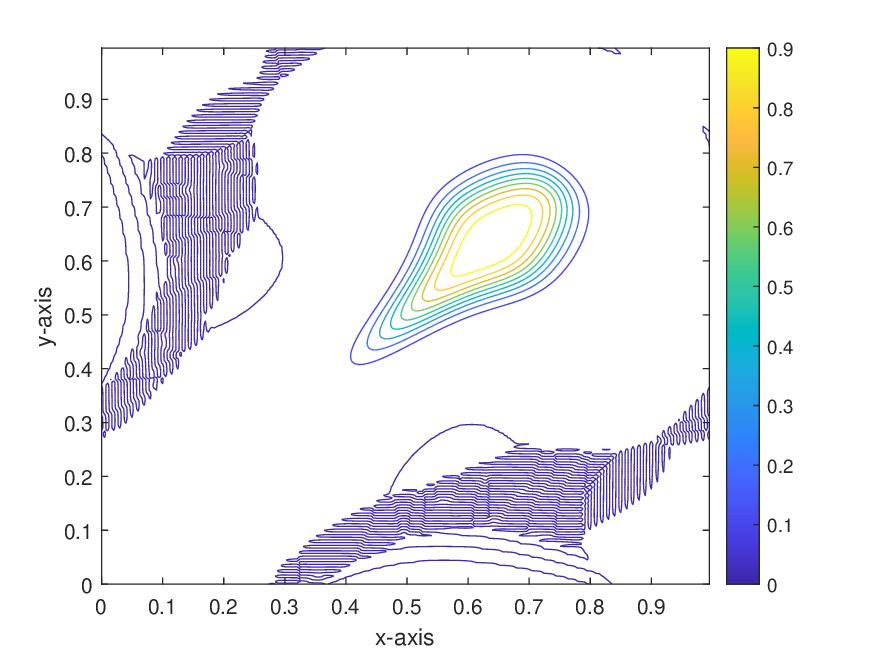}
 }
  \subfigure[$T=0.25$]
 {
  % \label{fig:burger:ini:sub2}
  \includegraphics[width=0.4\textwidth]{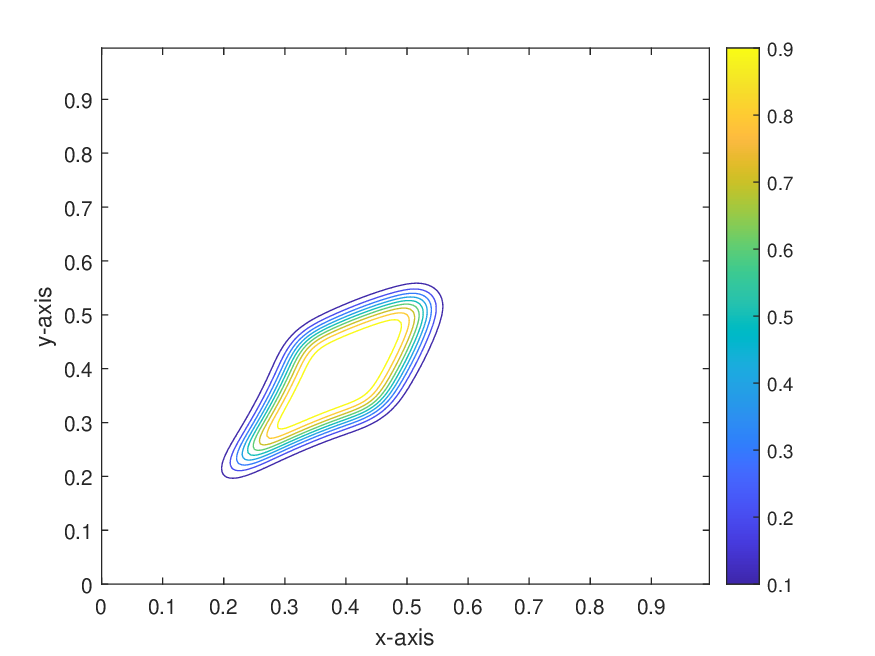}
 }
  \subfigure[$T=0.5$]
 {
  % \label{fig:burger:ini:sub4}
  \includegraphics[width=0.4\textwidth]{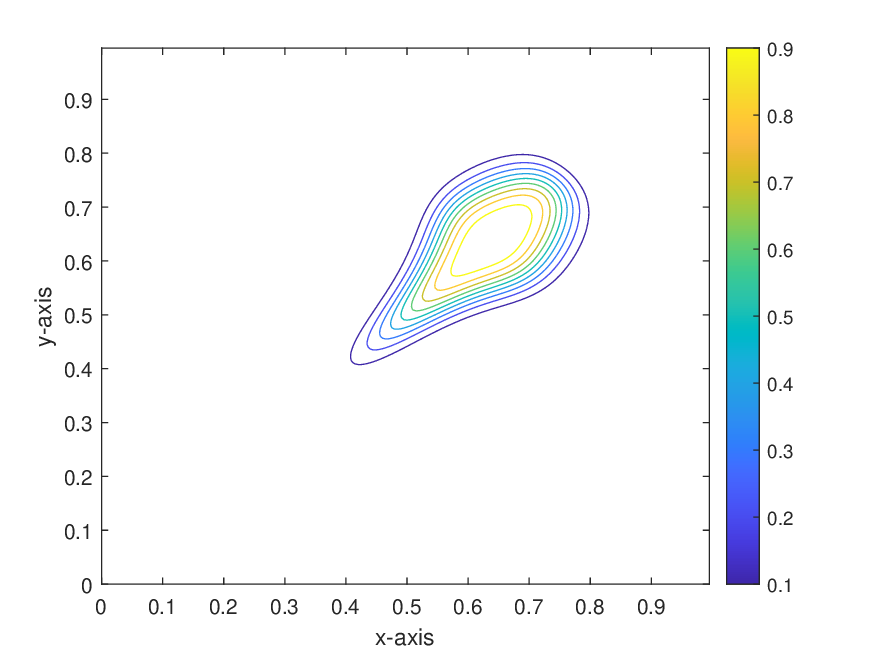}
 }
 \caption{Comparison  of contour plots obtain by HOC-ADI-Splitting scheme without the BP limiter (a-b) and with BP limiter (c-d) using $\nu=10^{-3}$.}
 \label{fig:case3:initical:u:nu1e3}
\end{figure}

\begin{figure} [!htbp]
 \centering 
 \subfigure[$T=0.25$]
 {
  % \label{fig:burger:ini:sub1}
  \includegraphics[width=0.4\textwidth]{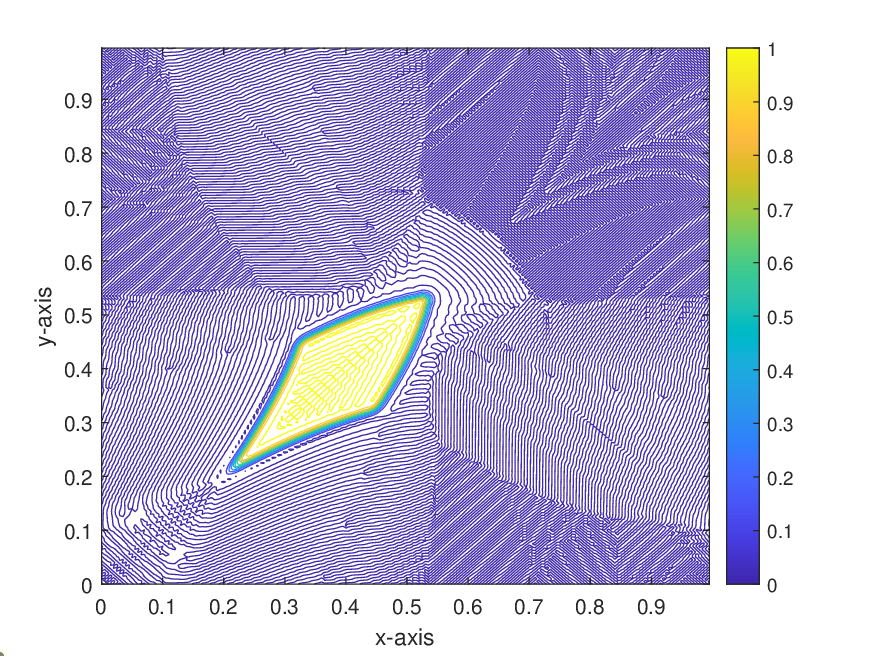}
 }
  \subfigure[$T=0.5$]
 {
  % \label{fig:burger:ini:sub3}
  \includegraphics[width=0.4\textwidth]{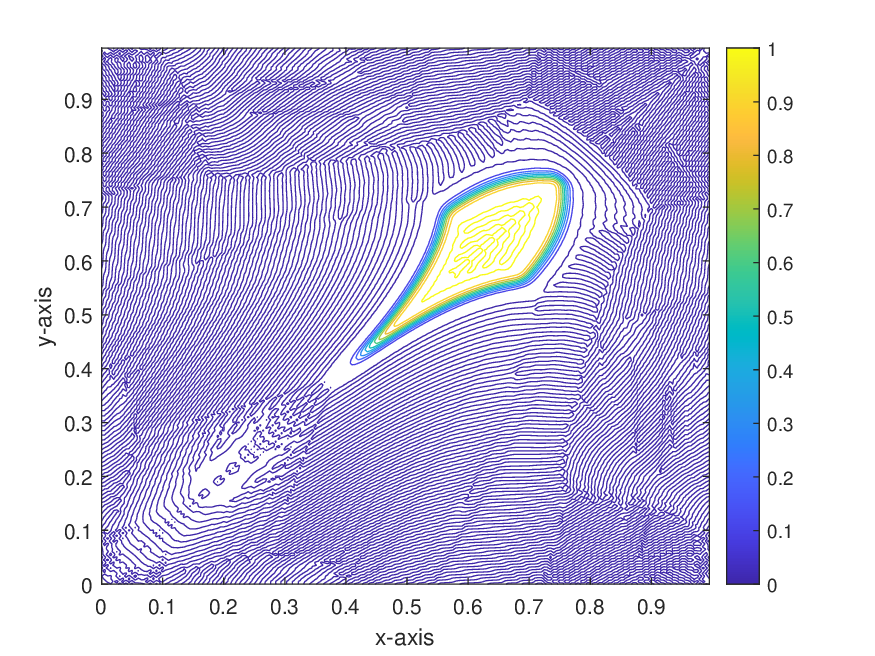}
 }
 \subfigure[$T=0.25$]
 {
  % \label{fig:burger:ini:sub2}
  \includegraphics[width=0.4\textwidth]{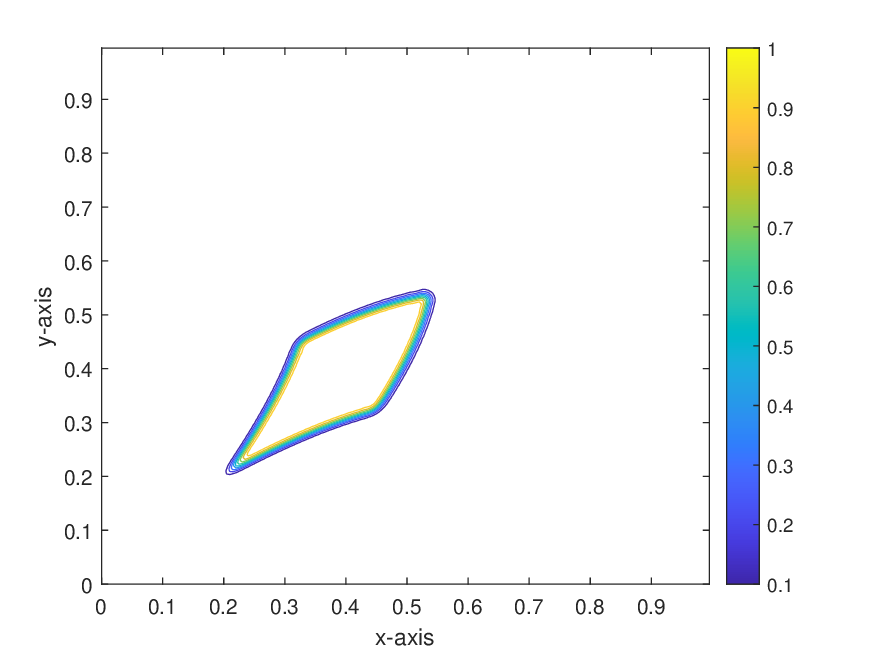}
 }
  \subfigure[$T=0.55$]
 {
  % \label{fig:burger:ini:sub4}
  \includegraphics[width=0.4\textwidth]{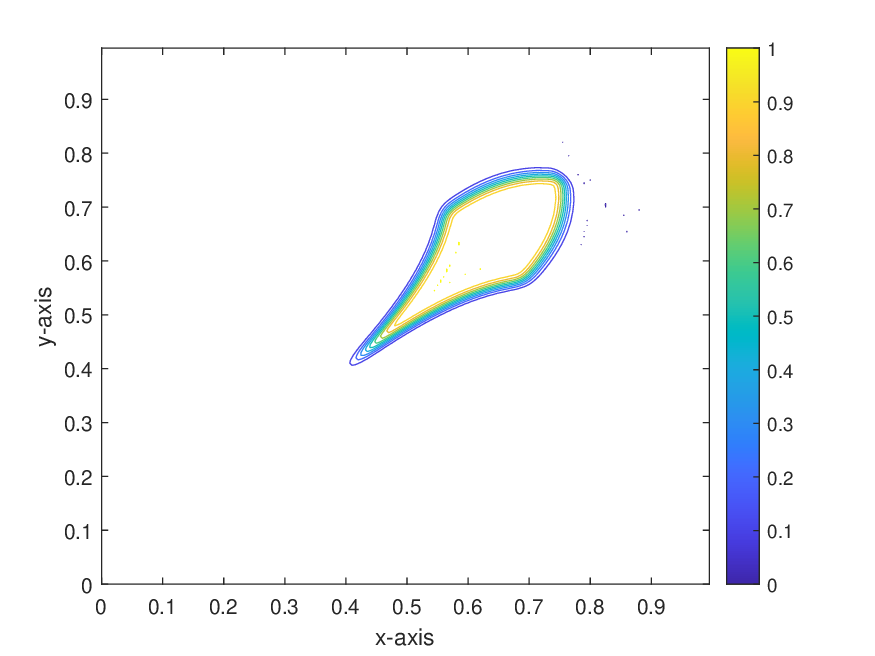}
 }
 \caption{Comparison  of contour plots between  the obtained numerical solution by HOC-ADI-Splitting scheme without the BP limiter (a-b) and with BP limiter (c-d) using $\nu=10^{-4}$.}
 \label{fig:case3:initical:u:nu1e4}
\end{figure}

In Case {3}, we further test the efficiency of our proposed scheme  with different $\nu$. The space steps are chosen as $N=200$, and the time steps $\tau =\frac{h}{2 {\max{|f'|} } }>\frac{h}{6 {\max{|f'|} } }$ are tested, where ${\max{|f'|} }=1$.
The velocity field and the initial value distribution are presented
in Fig. \ref{fig:burger:ini:0}.  In Fig. \ref{fig:case3:initical:u:nu1e3} and \ref{fig:case3:initical:u:nu1e4}, we draw the contour plots of the numerical solution calculated by the proposed Strang-ADI-HOC scheme without/with the BP limiter at $T= 0.25, 0.35$ and 0.5 by setting different $\nu$. As depicted in Fig. \ref{fig:case3:initical:u:nu1e3}, a larger diffusion coefficient of $\nu=10^{-3}$ leads to the rapid dissipation of the square wave over a short time period. In contrast, with a smaller diffusion coefficient  of $\nu=10^{-4}$, the steep gradients at the edges of the square wave persist throughout the dynamic evolution, as illustrated in Fig. \ref{fig:case3:initical:u:nu1e4}. Thus,  the results with the limiter accurately captures steep fronts without introducing numerical oscillations, even when convection dominates the dynamics.

%==========================================

%================================================

%================================================

\begin{figure} [!htbp]
 \centering 
 \subfigure[Without the BP limiter]
 {
   \label{fig:ploy:nonlim:surf}
  \includegraphics[width=0.45\textwidth]{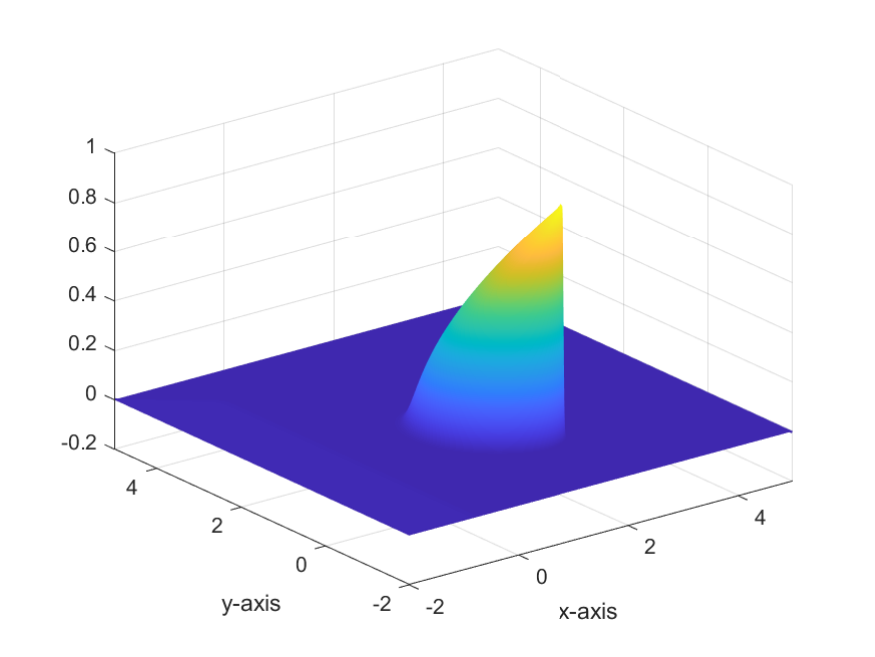}
 }
   \subfigure[Without the BP limiter]
 {
   \label{fig:ploy:nonlim:contour}
  \includegraphics[width=0.45\textwidth]{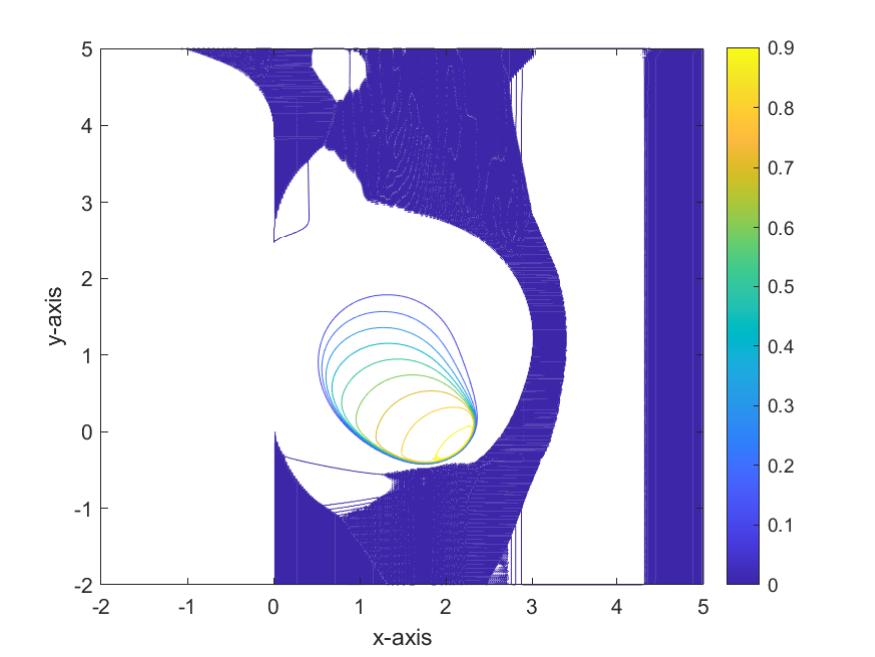}
 }
 \subfigure[With the BP limiter]
 {
  %\label{Fig. sub. 2}
  \includegraphics[width=0.45\textwidth]{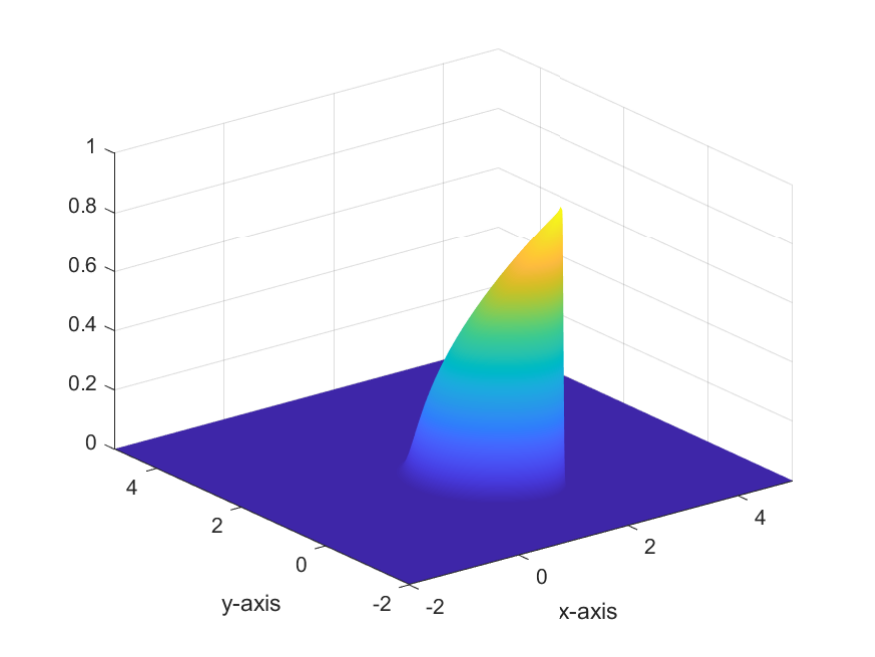}
 }
 \subfigure[With the BP limiter]
 {
  %\label{Fig. sub. 2}
  \includegraphics[width=0.45\textwidth]{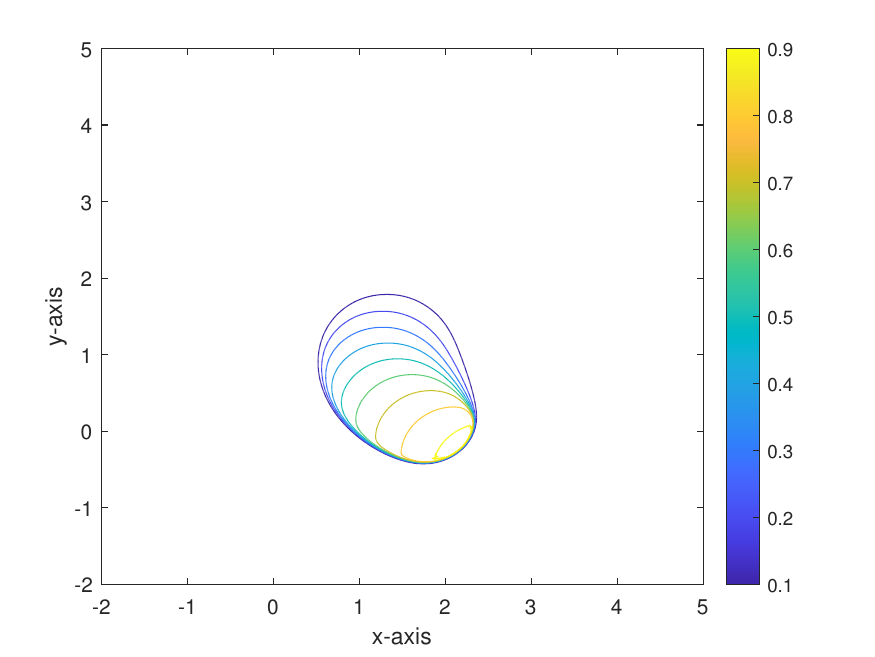}
 }
 \caption{Comparison  of the three-dimension views and contour plots obtained without (a-b) and with (c-d) the BP limiter.}
 \label{fig:2d:ploy:u:BP}
\end{figure}

\paragraph{Example 7}  
We consider the equation \cite{karlsen1997operator}
\begin{equation}
u_{t} +\left(u+(u-0. 25)^{3} \right)_{x} -\left(u+u^{2} \right)_{y} =\nu\left(u_{x x} +u_{y y} \right), \\
\quad
(x, y, t) \in[-2, 5] \times[-2, 5] \times[0, 1]. \notag
\end{equation}
with initial data given by
\begin{align*}
u_{0} (x, y)= \begin{cases} 1, & \text {for} (x-0. 25)^{2} +(y-2. 25)^{2} <0. 5 \\ 0, & \text {otherwise. } \end{cases}
\end{align*}
We set boundary values to zero, which is consistent with the initial data $u_{0}$. Fig. \ref{fig:2d:ploy:u:BP} displays the comparison of three-dimensional views and contour plots between the numerical solution obtained without/with the BP limiter at $T=1$
  with $\nu=5 \times 10^{-3},N=600,\tau= \frac{h}{{6\max{|f'|}}}$, where $\max{|f'|}=  2.6875$.
Combining Fig. \ref{fig:ploy:nonlim:surf} and \ref{fig:ploy:nonlim:contour}, it clearly shows that negative values and non-physical numerical oscillations appear if the BP technique is not applied, while the scheme
 with the BP limiter leads to accurate positive solutions.

% %==============================================
% \begin{figure} [!htbp]
%  \centering %Í¼Æ¬È«¾Ö¾ÓÖÐ
%  \subfigure
%  {
%   %\label{Fig. sub. 1}
%   \includegraphics[width=0.45\textwidth]{figure/2D_ploy_u_1500_1000.eps}
%  }
%  \subfigure
%  {
%   %\label{Fig. sub. 2}
%   \includegraphics[width=0.45\textwidth]{figure/2D_ploy_u_300_300.eps}
%  }
%  \caption{The comparison of the reference solution (left) and obtained numerical solution (right).}
%  %\label{eg:2d:refu}
% \end{figure}

% \begin{figure} [!htbp]
%  \centering
%  \includegraphics[width=10cm,height=6cm]{figure/2D_ploy_contour_non.eps}
%  \caption{
% Contour plot of computed solution without limiter. 
% }
%  \label{fig9}
% \end{figure}

% \begin{figure}[!htbp]
% \centering
% \subfigure{
%  \includegraphics[width=5.5cm]{figure/2D_ploy_T0.eps}
%  %\caption{fig:eg1:diff_u:mass_err}
%}
% \quad
% \subfigure{
%  \includegraphics[width=5.5cm]{figure/2D_ploy_T02.eps}
%}
% \quad
% \subfigure{
%  \includegraphics[width=5.5cm]{figure/2D_ploy_T04.eps}
%}
% \quad
% \subfigure{
%  \includegraphics[width=5.5cm]{figure/2D_ploy_T06.eps}
%}
% \quad
% \subfigure{
%  \includegraphics[width=5.5cm]{figure/2D_ploy_T08.eps}
%}
% \quad
% \subfigure{
%  \includegraphics[width=5.5cm]{figure/2D_ploy_T1.eps}
%}
% \caption{The surface plots of the computed solutions at different time with $\nu = 0.05, N_t = 250, N = 300 $. Top: t = 0; 0.2; Middle: t = 0.4; 0.6; Bottom: t = 0.8; 1.0.
%}
% \label{plus7}
% \end{figure}

%==============================================

\begin{figure} [!ht]
\centering 
\subfigure[Without the  BP  limiter  ]
{
%\label{Fig. sub. 1}
\includegraphics[width=0.4\textwidth]{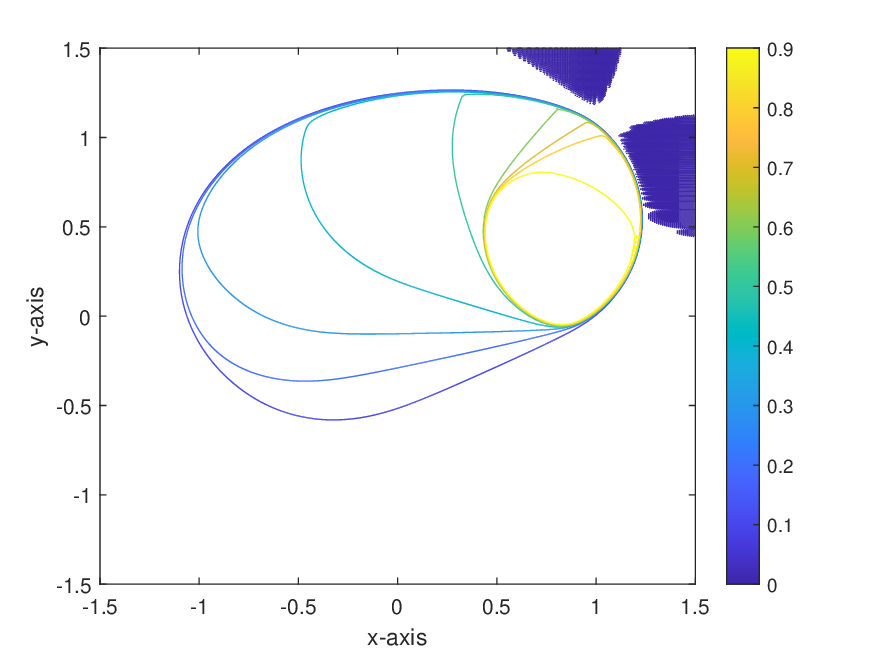}
}
\subfigure[With the BP limiter  ]
{
%\label{Fig. sub. 2}
\includegraphics[width=0.4\textwidth]{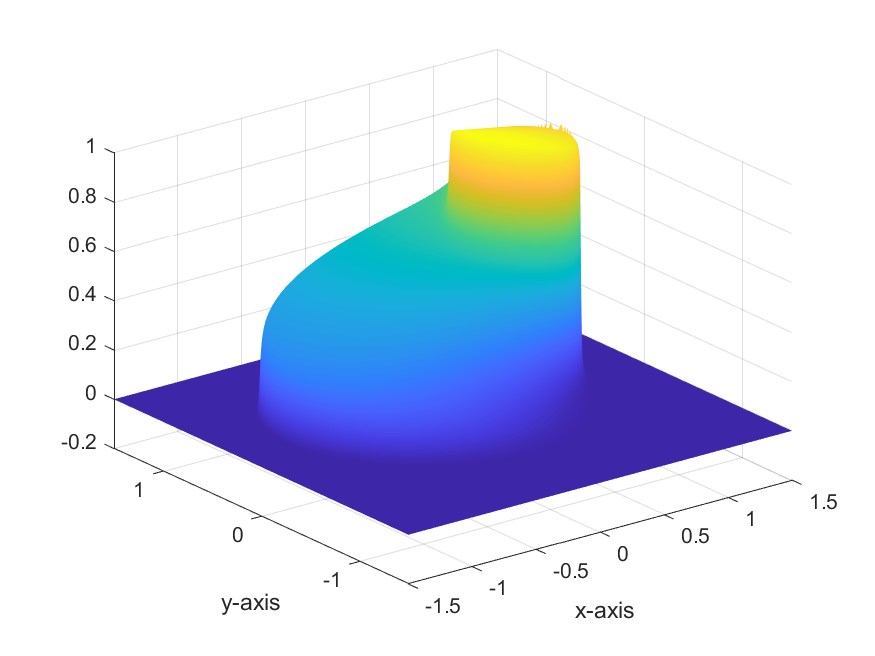}
}
\subfigure[Without the  BP  limiter ]
{
%\label{Fig. sub. 1}
\includegraphics[width=0.4\textwidth]{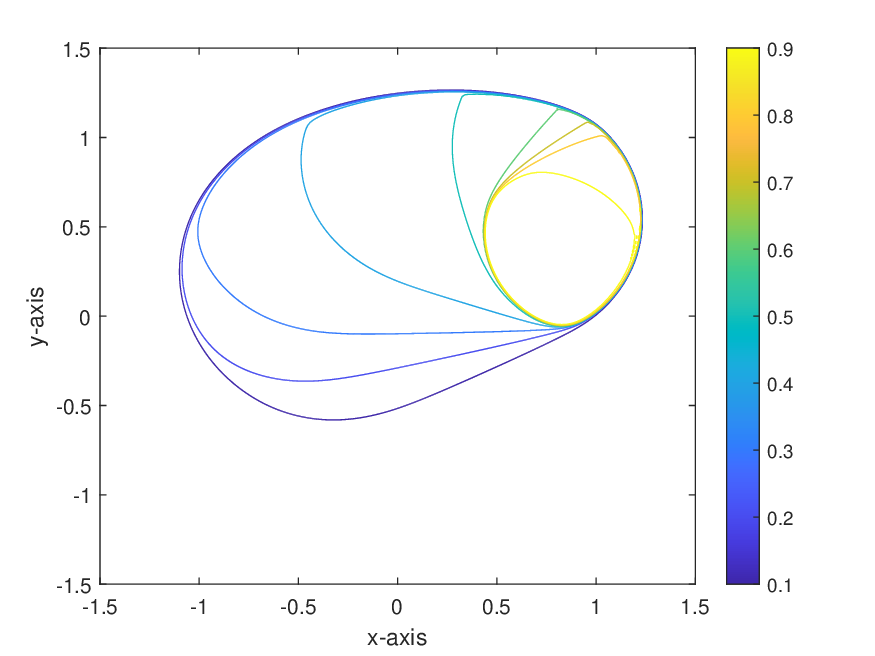}
}
\subfigure[With the BP limiter ]
{
%\label{Fig. sub. 2}
\includegraphics[width=0.4\textwidth]{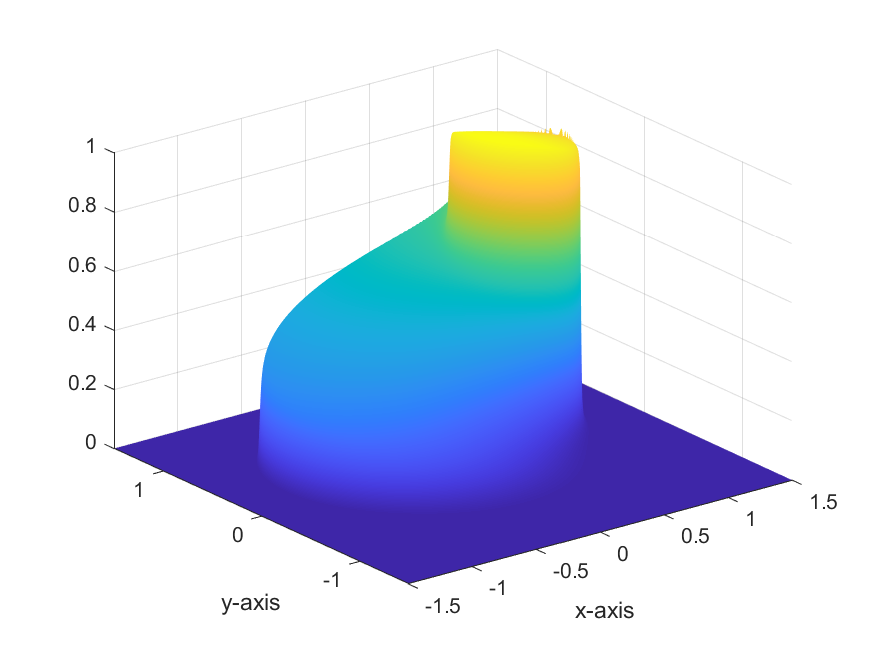}
}
\caption{Comparisons of three-dimension views and contour plots obtained by HOC-ADI-Splitting scheme without (a-b) / with (c-d) BP limiter  using $\nu=4 \times 10^{-3}$.}
\label{2d:frac:u}
\end{figure}

\textbf{Example 8} % \cite{ACOSTA20101397}
%(Two-Dimensional Buckley¨CLeverett Equation)
We now present an example where we generate approximate solutions to the equation \cite{karlsen1997operator,He2022AFM}
\begin{equation*}
u_{t} +f(u)_{x} +g(u)_{y} =\gamma\left(u_{x x} +u_{y y} \right),
\end{equation*}
with initial data
$$
u_{0} (x, y)= \begin{cases} 1, & \text {for } x^{2} +y^{2} <0. 5 \\ 0, & \text {otherwise} \end{cases}
$$
and
$
g(u) =\frac{u^{2} }{u^{2} +(1-u)^{2} },
$
$
f(u) =g(u)\left(1-5(1-u)^{2} \right).
$
The flux functions $f$ and $g$ are both "S-shaped" with $f(0)=g(0)=0$ and $f(1)=g(1)=1$. This problem is motivated from two-phase flow in porous media with a gravitation pull in the $x$-direction. The computational domain is $[-3,3] \times[-3,3]$ and boundary values are again put equal to zero.
%The results are presented in Fig. \ref{fig11} at $T =0. 5$.
We take $N=600, \nu = 4 \times 10^{-3}, \tau =  \frac{h}{6\max{|f'|}}$ and $T=0.5$, where $\max{|f'|}\approx 3.13$.
We displayed the three-dimensional views and plots of the numerical solution without/with the BP limiter  in Fig. \ref{2d:frac:u}. No exact solution to this problem is available, but if compared with the numerical solutions reported in \cite{karlsen1997operator} and \cite{He2022AFM}, our solutions seem to converge to the correct solution. Thus, our scheme with the limiter  provides a highly satisfactory approximation to this model with a nonlinear, degenerate diffusion. 

%In summary, the effectivity of
%the present method is verified by good coincide of this contour with those ones in \cite{karlsen1997operator} \cite{He2022AFM}.
% %
% \begin{table} [!htbp]
% %h£ºhear£¬t£ºtop£¬b£ºbottom£¬p£ºpage£¬ÏÂÒ»Ò³¡£
% \centering
% \caption{Verification of the limiter's bounding effect of problem 8 with $\nu =0. 01, \tau= \frac{h}{6 \mathrm{max} u}$ at $T =0. 5$.
% }
% \label{tab_5}
% \begin{tabular}{|c|cc|cc|}
% %¼¸ÁÐ¾ÍÐ´¼¸¸öc£¬±íÊ¾ÄÚÈÝ¾ÓÖÐÐ´£¬ÄÄÀïÐèÒª·Ö¸îÏß¾ÍÔÚÄÄÀï¼Ó¡°|¡±
%  %\toprule
%   %×îÉÏÃæµÄºáÏß
% %ÌîÐ´±í¸ñÄÚÈÝ£¬½¨ÒéÄÚÈÝ¶àµÄÊ±ºòÓÃ{} & {} & {}
%  % \midrule
%  %ÐÐ¼äÐèÒª·Ö¸îÏß¾Í¼ÓÕâÌõÓï¾ä£¬²»ÐèÒª¾Í²»¼Ó
%  %Ö»ÒªÊÇÐÐ¼äµÄÄÚÈÝ£¬Ìî³äÍêÖ®ºó¶¼Òª¼ÓÉÏ¡°\\¡±
% %\bottomrule
%  %×îÏÂÃæµÄºáÏß
% \hline & \multicolumn{2}{c|}{\text {The scheme without limiter} } & \multicolumn{2}{c|}{\text {The scheme with limiter } } \\
% \hline$N$ & $\mathrm{max(U)-max(u)}$& $\mathrm{min(u)-min(U)}$ & $\mathrm{max(U)-max(u)}$ & $\mathrm{min(U)-min(u)}$\\
% \hline 400 & $  1. 5 E -3$ & $ -8. 6 E -3$ &
% $ 1. 7 E -3$ & $ 6. 4654 E -157 $\\
% \hline 500 & $4. 6E -3$ & $ -4. 3583 E -5$ &
% $4. 4E -3$ & $ 2. 3203 E -188$\\
% \hline 600 & $4. 7E -3$ & $ -1. 6001 E -54$ &
% $4. 5E -1$ & $ 1. 1693 E -213$\\
% \hline 700 & $ 4. 7 E -1$ & $-12. 0146 E -102$ &
% $ 4. 6 E -1$ & $ 9. 4572 E -235$\\
% \hline
% \end{tabular}
% \end{table}

% No exact solution to problem (6. 9)¨C(6. 12) is available, but
% if compared with the numerical solutions reported in [20], our solutions seem to converge
% to the correct entropy solution.

% Fig. 5 shows ¡°log-log plots¡± of the error
% and numerical convergence rates in ?x for " = 0:1, 0:01 and CFL numbers 1, 2,
% and 4.

% One of the well-known difficulties in solving the Burgers equation with a small viscosity coefficient is the steepening of the solution curves with time and the development of a shock-like
% discontinuity. Here, we simulate this example using MECM3 until t = 1. 0 with different ?? from 10?1
% to 10?4 and also illustrate the numerical solutions for small ?? = 10?5 at different time levels to show
% the adaptability of MECM3 for a small viscosity coefficient number. The results are illustrated in
% Fig. 12, which is shown that MECM3 exhibits good performance under small viscosity coefficients.

 \section{Conclusions}
In this paper, we have demonstrated that fourth order accurate compact finite difference schemes combined with Strang splitting method for nonlinear convection diffusion problems. A series of theoretical analyses and numerical experiments demonstrate that this improved method effectively optimizes the original algorithm. By decoupling the convection and diffusion terms, the Strang method allows each term to be handled more efficiently, reducing computational constraints and costs. The ADI method's application to high-dimensional problems ensures that the complexity is manageable, making the method more practical and efficient for real-world applications.
Our schemes preserve bound and mass by utilizing the BP limiter proposed by \cite{li2018A}.
To extend the scheme to two-dimensional, we employ ADI method to further reduce the computational cost.
Using the properties of the limiter and M-matrix, we give the proof of bound-preserving, mass conservation and stability.
% In fact, this paper only gives sufficient conditions of our schemes to preserve bound, but not necessary conditions.
% After a large number of numerical experiments, we found that the actual calculation only needs to satisfy the CFL condition from the convection term, which further weakens the network ratio condition.
Numerical results suggest well performance by the proposed schemes.

%%%%%%%%%%%%%%%%%%%%%%%%%%%%%%%%%%%%%%%%%%%%%%%%%
%\bibliographystyle{plain}
%% bib
%\bibliographystyle{spmpsci}
%\bibliography{ref}

\bibliographystyle{abbrvnat}
%\bibliography{reference}

%%\bibliographystyle{plain}
%\bibliographystyle{elsarticle-num}
%%\bibliographystyle{unsrt}
%%\bibliographystyle{abbrv}
%%\bibliographystyle{spmpsci}  % mathematics and physical sciences
%%\bibliographystyle{spmpsci}
%\bibliography{ref}

\begin{thebibliography}{99}
 
\bibitem{acosta2010mollification}
  Acosta C D, Mejia C E. A mollification based operator splitting method
  for convection diffusion equations.
  Computers \& mathematics with applications, 2010, 59(4): 1397--1408.

\bibitem{ali1992}
Ali A H A, Gardner G A, Gardner L R T. A collocation solution for Burgers equation using cubic B-spline finite
elements. Computer Methods in Applied Mechanics and Engineering, 1992, 100:325--337.

\bibitem{berikelashvili2007convergence}
  Berikelashvili G, Gupta M M, Mirianashvili M.
  Convergence of fourth order compact difference schemes for
  three-dimensional convection diffusion equations.
  SIAM Journal on Numerical Analysis, 2007, 45(1): 443--455.

\bibitem{bertolazzi2005second}
 Bertolazzi E, Manzini G. A second-order maximum principle preserving finite
 volume method for steady convection diffusion problems.
 SIAM Journal on Numerical Analysis, 2005, 43(5): 2172--2199.

\bibitem{chu1998three}
 Chu P C, Fan C. A three-point combined compact difference scheme.
 Journal of Computational Physics, 1998, 140(2): 370--399.

\bibitem{douglas1962alternating}
  Douglas J. Alternating direction methods for three space variables.
  Numerische Mathematik, 1962, 4: 41--63.

\bibitem{einkemmer2016overcoming}
  Einkemmer L, Ostermann A. Overcoming order reduction in diffusion-reaction splitting.
  Part 2: Oblique boundary conditions.
  SIAM Journal on Scientific Computing, 2016, 38(6): A3741--A3757.
 
\bibitem{fattah1985dispersion}
 Fattah Q N, Hoopes J A. Dispersion in anisotropic, homogeneous, porous media.
 Journal of Hydraulic Engineering, 1985, 111(5): 810--827.

\bibitem{genty2011maximum}
  Genty A, Le Potier C. Maximum and minimum principles for radionuclide transport
 calculations in geological radioactive waste repository: comparison between
 a mixed hybrid finite element method and finite volume element discretizations.
 Transport in porous media, 2011, 88(1): 65--85.

\bibitem{gottlieb2009high}
 Gottlieb S, Ketcheson D I, Shu C W. High order strong stability preserving time discretizations.
 Journal of Scientific Computing, 2009, 38(3): 251--289.
 
%\bibitem{gustafsson2013time} 
% Gustafsson B, Kreiss H O, Oliger J. Time Dependent Problems and Difference Methods.
% New York: Wiley, 2013.

\bibitem{He2022AFM}
  He Q, Du W, Shi F, et al.
 A fast method for solving time-dependent nonlinear convection diffusion problems.
 Electronic Research Archive, 2022, 30(6): 2165--2182.

\bibitem{hundsdorfer2003numerical}
  Hundsdorfer W H, Verwer J G, Hundsdorfer W H.
  Numerical solution of time-dependent advection-diffusion-reaction equations. Berlin: Springer, 2003.

\bibitem{JIA2011387}
  Jia H, Li K. A third accurate operator splitting method.
  Mathematical and computer modelling, 2011, 53(1-2): 387--396.

\bibitem{2002A}
  Kalita J C, Dalal D C, Dass A K. A class of higher order compact schemes for the
  unsteady two-dimensional convection diffusion equation with variable convection coefficients.
  International Journal for Numerical Methods in Fluids, 2002, 38(12): 1111--1131.

\bibitem{karaa2004high}
  Karaa S, Zhang J. High order ADI method for solving unsteady convection diffusion problems.
  Journal of Computational Physics, 2004, 198(1): 1--9.

\bibitem{karlsen1997operator}
  Karlsen K H, Risebro N H. An operator splitting method for nonlinear convection diffusion equations.
 Numerische Mathematik, 1997, 77(3): 365--382.

%\bibitem{lapidus2011numerical}
%  Lapidus L, Pinder G F. Numerical Solution of Partial Differential Equations in Science and Engineering.
%  John Wiley \& Sons, 2011.

\bibitem{lee2019energy}
 Lee S, Shin J. Energy stable compact scheme for Cahn-Hilliard equation with periodic boundary condition.
 Computers \& Mathematics with Applications, 2019, 77(1): 189--198.
 
\bibitem{lele1992compact}
  Lele S K. Compact finite difference schemes with spectral-like resolution.
  Journal of computational physics, 1992, 103(1): 16--42.

\bibitem{li2018A}
  Li H, Xie S, Zhang X. A high order accurate bound-preserving compact
  finite difference scheme for scalar convection diffusion equations.
  SIAM Journal on Numerical Analysis, 2018, 56(6): 3308--3345.

\bibitem{mohebbi2010high}
  Mohebbi A, Dehghan M. High-order compact solution of the one-dimensional
  heat and advection-diffusion equations.
  Applied mathematical modelling, 2010, 34(10): 3071--3084.
 
\bibitem{parlange1980water}
 Parlange J Y. Water transport in soils.
 Annual review of fluid mechanics, 1980, 12(1): 77--102.

\bibitem{peaceman2000fundamentals}
  Peaceman D W. Fundamentals of Numerical Reservoir Simulation. Elsevier, 2000.

\bibitem{peaceman1955numerical}
  Peaceman D W, Rachford, Jr H H. The numerical solution of parabolic and elliptic differential equations.
  Journal of the Society for industrial and Applied Mathematics, 1955, 3(1): 28--41.

\bibitem{plemmons1977m}
  Plemmons R J. M-matrix characterizations. I-nonsingular M-matrices.
  Linear Algebra and its applications, 1977, 18(2): 175--188.

\bibitem{poole1974survey}
  Poole G, Boullion T. A survey on M-matrices.
  SIAM review, 1974, 16(4): 419--427.

\bibitem{rui2010mass}
 Rui H, Tabata M. A mass-conservative characteristic finite element scheme
 for convection diffusion problems.
 Journal of Scientific Computing, 2010, 43: 416--432.

\bibitem{shen2021discrete}
  Shen J, Zhang X. Discrete maximum principle of a high order finite difference scheme for a generalized {Allen-Cahn} equation. 
Communications in Mathematical Sciences, 2022, 20(5): 1409--1436. 

\bibitem{shu1988efficient}
  Shu C W, Osher S. Efficient implementation of essentially non-oscillatory shock-capturing schemes.
  Journal of computational physics, 1988, 77(2): 439--471.

\bibitem{srinivasan2018positivity}
  Srinivasan S, Poggie J, Zhang X. A positivity-preserving high order
 discontinuous Galerkin scheme for convection diffusion equations.
 Journal of Computational Physics, 2018, 366: 120--143.

\bibitem{strang1968}
 Strang G. On the construction and comparison of difference schemes.
 SIAM journal on numerical analysis, 1968, 5(3): 506--517.

%\bibitem{Sunbook}
%  Sun Z Z. Numerical Methods of Partial Differential Equations (2nd version, in Chinese).
%  Beijing: Science Press, 2012.

\bibitem{tian2011rational}
  Tian Z F. A rational high-order compact ADI method for unsteady convection diffusion equations.
  Computer Physics Communications, 2011, 182(3): 649--662.

\bibitem{tian2007fourth}
 Tian Z F, Ge Y B. A fourth-order compact ADI method for solving
 two-dimensional unsteady convection diffusion problems.
 Journal of Computational and Applied Mathematics, 2007, 198(1): 268--286.

\bibitem{van2019positivity}
  van der Vegt J J W, Xia Y, Xu Y. Positivity preserving limiters for time-implicit
  higher order accurate discontinuous Galerkin discretizations.
  SIAM journal on scientific computing, 2019, 41(3): A2037--A2063.

\bibitem{2015shiyan}
  Wang H, Shu C W, Zhang Q. Stability and error estimates of local discontinuous Galerkin
 methods with implicit-explicit time-marching for advection-diffusion problems.
 SIAM Journal on Numerical Analysis, 2015, 53(1): 206--227.

\bibitem{wang2005new}
  Wang T. New characteristic difference method with adaptive mesh for
 one-dimensional unsteady convection-dominated diffusion equations.
 International Journal of Computer Mathematics, 2005, 82(10): 1247--1260.

\bibitem{wang2015fourth}
  Wang Y M. Fourth-order compact finite difference methods and monotone
  iterative algorithms for quasi-linear elliptic boundary value problems.
  SIAM Journal on Numerical Analysis, 2015, 53(2): 1032--1057.

\bibitem{xie2008numerical}
  Xie S S, Heo S, Kim S, Woo G, Yi S. Numerical solution of one-dimensional Burgers equation
 using reproducing kernel function. J. Comput. Appl. Math., 2008, 214(2): 417--434.

\bibitem{xiong2015high}
  Xiong T, Qiu J M, Xu Z. High order maximum-principle-preserving
  discontinuous Galerkin method for convection diffusion equations.
  SIAM Journal on Scientific Computing, 2015, 37(2): A583--A608.

\bibitem{zhang2021numerical}
  Zhang L, Ge Y. Numerical solution of nonlinear advection diffusion
 reaction equation using high-order compact difference method.
 Applied Numerical Mathematics, 2021, 166: 127--145.

\bibitem{zhang2022positivity}
  Zhang L, Ge Y, Wang Z. Positivity-preserving high-order
 compact difference method for the Keller-Segel chemotaxis model.
 Mathematical Biosciences and Engineering, 2022, 19(7): 6764--6794.

\bibitem{zhang2012maximum}
 Zhang X, Liu Y, Shu C W. Maximum-principle-satisfying high order finite volume
 weighted essentially nonoscillatory schemes for convection diffusion equations.
 SIAM Journal on Scientific Computing, 2012, 34(2): A627--A658.

\bibitem{zhuang2014positivity}
  Zhuang C, Zeng R. A positivity-preserving scheme for the simulation of streamer discharges
 in non-attaching and attaching gases.
 Communications in Computational Physics, 2014, 15(1): 153--178.
 
\bibitem{zlatev1984implementation}
 Zlatev Z, Berkowicz R, Prahm L P.
 Implementation of a variable stepsize variable formula method in the time-integration part of a
 code for treatment of long-range transport of air pollutants.
 Journal of Computational Physics, 1984, 55(2): 278-301.

\bibitem{kaya2001explicit}
   Kaya  D. 
 An explicit solution of coupled Burgers equations by decomposition method.
 International Journal of Mathematics and Mathematical Sciences, 2001, 27: 675-680.

\bibitem{FU2017jsc}
  Fu K, Liang D. 
 The time second order mass conservative characteristic FDM for advection–diffusion equations in high dimensions. 
 Journal of Scientific Computing, 2017, 73: 26-49.

\bibitem{FU2019sisc}
  Fu K, Liang D. 
 A mass-conservative temporal second order and spatial fourth order characteristic finite volume method for atmospheric pollution advection diffusion problems. 
 SIAM Journal on Scientific Computing, 2019, 41(6): B1178-B1210.

% \bibitem{qin2023}
%  Qin D, Fu K, Liang D. 
%  Positivity preserving temporal second-order spatial fourth-order conservative characteristic methods for convection dominated diffusion equations. Computers & Mathematics with Applications, 2023, 149: 190-202.

\bibitem{liao2010maximum}
    Liao H L, Sun Z Z.
     Maximum norm error bounds of {ADI} and compact {ADI} methods for solving parabolic equations.
     Numerical Methods for Partial Differential Equations, 2010,26:37-60.

\bibitem{shen2016maximum}
    Shen J, Tang T, Yang J.
    On the maximum principle preserving schemes for the generalized {Allen--Cahn} equation. 
    Communications in Mathematical Sciences, 2016, 14(6): 1517-1534.


\bibitem{bookADR} 
     Hundsdorfer W, Verwer J. 
     Numerical Solution of  Time-Dependent Advection-Diffusion-Reaction Equations.
     Heidelberg: Springer, 2003.

\bibitem{wang1d} 
    Wang Y, Xie S S, Fu H F. 
    Efficient, linearized high-order compact difference schemes for nonlinear parabolic equations {I: One-dimensional} problem.
    Numerical Methods for Partial Differential Equations, 2023, 39(2): 1529-1557.

\bibitem{Zhang2002HOC}
   Zhang J, Sun H, Zhao J J. 
   High order compact scheme with multigrid local mesh refinement procedure for convection diffusion problems. 
   Computer Methods in Applied Mechanics and Engineering, 2002, 191(41-42): 4661-4674.

\bibitem{qiao2008}
   Ito K, Qiao Z.
    A high order compact {MAC} finite difference scheme for the Stokes equations: augmented variable approach. 
   Journal of Computational Physics, 2008, 227(17): 8177-8190.

\bibitem{shu1994TVB}
  Cockburn  B, Shu  C W.
Nonlinearly stable compact schemes for shock calculations. 
SIAM Journal on Numerical Analysis, 1994, 31(3): 607–627.

\bibitem{li2003TVB}
Li H, Zhang X X. 
A high order accurate bound-preserving compact finite difference scheme for two-dimensional incompressible flow. 
Communications on Applied Mathematics and Computation, 2024, 6(1): 113-141.

\bibitem{qin2023}
Qin D, Fu K, Liang D. 
Positivity preserving temporal second-order spatial fourth-order conservative characteristic methods for convection dominated diffusion equations.
Computers \& Mathematics with Applications, 2023, 149: 190-202.

\end{thebibliography}
\end{document}